\documentclass[11pt,oneside,a4paper,mathscr]{amsart}
\usepackage{amsmath,amssymb}
\usepackage{tikz-cd}
\usepackage{adjustbox}
\usepackage{amscd,euscript}
\usepackage[frame,cmtip,curve,arrow,matrix,line,graph]{xy}
\usepackage{hyperref}
\usepackage{bm}
\usepackage{here}
\usepackage{enumitem}
\hypersetup{hidelinks}
\usepackage{xcolor}

\newcommand{\conjeq}{%
  \mathrel{%
    \mathop{%
      \vcenter{\hbox{\scalebox{1.35}{$\simeq$}}}%
    }\limits_{\mathrm{conj.}}%
  }%
}

\tikzset{
  mybox/.style={
    draw,
    rounded corners=3pt,
    fill=gray!8,
    inner sep=5pt,
    align=center
  },
  resultbox/.style={
    draw,
    rounded corners=3pt,
    fill=gray!14,
    inner sep=5pt,
    align=center
  }
}
\title[The Dolbeault geometric Langlands-type conjecture]{The Dolbeault geometric Langlands-type conjecture for smooth surfaces}
\author{Yukinobu Toda}

\newtheorem{thm}{Theorem}[section]

\newtheorem{cor}[thm]{Corollary}
\newtheorem{prop}[thm]{Proposition}

\newtheorem{conj}[thm]{Conjecture}
\newtheorem{conj2}[thm]{`Conjecture'}
\newtheorem{lemma}[thm]{Lemma}

\theoremstyle{definition}
\newtheorem{defn}[thm]{Definition}

\newtheorem{thm*}[thm]{Theorem$^*$}
\newtheorem{remark}[thm]{Remark}

\newtheorem{example}[thm]{Example}

\newtheorem{assum}[thm]{Assumption}

\newcommand{\comment}[1]{}

\renewcommand{\leq}{\leqslant}
\renewcommand{\geq}{\geqslant}

\newcommand{\rH}{\mathrm{H}}
\newcommand{\rW}{\mathrm{W}}

\newcommand{\Y}{\mathcal{Y}}

\newcommand{\LL}{\operatorname{L}}
\newcommand{\Pic}{\operatorname{Pic}}

\newcommand{\X}{\mathcal{X}}

\newcommand{\cH}{\mathcal{H}}

\newcommand{\cL}{\mathcal{L}}

\newcommand{\cC}{\mathcal{C}}

\newcommand{\cP}{\mathcal{P}}
\newcommand{\cE}{\mathcal{E}}

\newcommand{\cO}{\mathcal{O}}
\newcommand{\cD}{\mathcal{D}}
\newcommand{\ev}{\mathrm{ev}}
\newcommand{\bgm}{B\mathbb{G}_m}

\newcommand{\diasquare}{\ar@{}[rd]|\square}

\newcommand{\cHom}{\mathcal{H}om}
\newcommand{\rank}{\operatorname{rank}}

\newcommand{\rB}{\mathrm{B}}
\newcommand{\rC}{\mathrm{C}}

\newcommand{\Coh}{\operatorname{Coh}}

\newcommand{\id}{\operatorname{id}}
\newcommand{\Ext}{\operatorname{Ext}}
\newcommand{\Ind}{\operatorname{Ind}}
\newcommand{\IndCoh}{\operatorname{IndCoh}}
\newcommand{\IndL}{\operatorname{IndL}}
\newcommand{\Hom}{\operatorname{Hom}}

\newcommand{\Spec}{\operatorname{Spec}}

\newcommand{\GL}{\operatorname{GL}}

\newcommand{\ch}{\operatorname{ch}}

\renewcommand{\Im}{\operatorname{Im}}

\newcommand{\cl}{\mathrm{cl}}

\newcommand{\inclusion}{\ar@<-0.3ex>@{^{(}->}[r]}
\newcommand{\iinclusion}{\ar@<-0.3ex>@{^{(}->}[rr]}
\newcommand{\linclusion}{\ar@<0.3ex>@{_{(}->}[l]}
\newcommand{\dinclusion}{\ar@<-0.3ex>@{^{(}->}[d]}
\newcommand{\uinclusion}{\ar@<-0.3ex>@{^{(}->}[u]}

\newcommand{\ssslash}{/\!\!/}

\newcommand{\wt}{\mathrm{wt}}
\newcommand{\Hig}{\mathrm{Higgs}}

\newcommand{\Bun}{\mathrm{Bun}}

\makeatletter
\newcommand{\colim@}[2]{%
  \vtop{\m@th\ialign{##\cr
    \hfil$#1\operator@font colim$\hfil\cr
    \noalign{\nointerlineskip\kern1.5\ex@}#2\cr
    \noalign{\nointerlineskip\kern-\ex@}\cr}}%
}
\newcommand{\colim}{%
  \mathop{\mathpalette\colim@{}}\nmlimits@
}
\makeatother

\usetikzlibrary{positioning,fit,calc}
\tikzstyle{block}=[draw=black, text width=1cm, minimum height=2cm, align=center] 
\tikzstyle{block2}=[draw=black, text width=2cm, minimum height=1cm, align=center] 
\tikzstyle{block3}=[draw=black, text width=2cm, minimum height=1cm, align=center] 

\makeatletter

\@addtoreset{equation}{section}	
\makeatother

\begin{document}
\begin{abstract}
Using a general theory of limit categories for symmetric stacks, we construct their semiorthogonal decompositions for moduli stacks of pure one-dimensional coherent sheaves on smooth surfaces with nef anti-canonical divisors. The components of these semiorthogonal decompositions are quasi-BPS categories associated with the corresponding local Calabi--Yau threefolds.

By comparing these semiorthogonal decompositions with those of dg-categories of coherent sheaves on moduli stacks of semistable sheaves, we propose a Dolbeault geometric Langlands-type conjecture for surfaces, generalizing the Dolbeault geometric Langlands conjecture for $\GL_r$-Higgs bundles. Our conjecture includes a correspondence for Higgs bundles twisted by sufficiently ample line bundles on curves.

We discuss examples and construct Hecke actions on limit categories and, under certain conditions, on quasi-BPS categories. For symplectic surfaces, these actions are expected to provide a categorical lift of the corresponding Hecke operators on BPS cohomology constructed by Davison--Hennecart--Kinjo--Schiffmann--Vasserot.
\end{abstract}
\maketitle
\setcounter{tocdepth}{2}
\tableofcontents

\section{Introduction}
\subsection{Overview}
The purpose of this paper is to formulate a \textit{Dolbeault geometric Langlands-type conjecture} for smooth quasi-projective surfaces as an equivalence between certain dg-categories.
It states the following (below $k$ is an algebraically closed field of 
characteristic zero): 
\begin{conj}\label{intro:conj0}
Let $S$ be a smooth quasi-projective surface over $k$ such that $-K_S$ is nef. For a non-zero effective curve class $\beta$ on $S$ and $\chi,w \in \mathbb{Z}$, 
let 
\begin{align*}\mathfrak{M}_S(\beta, \chi)\to \rB\end{align*}
be the derived moduli stack of 
compactly supported pure one-dimensional sheaves $E$ on $S$ 
with $([E], \chi(E))=(\beta, \chi)$, with derived Hilbert--Chow map to $\rB$. 
Let $\sigma$ be a Bridgeland stability condition on $\Coh_{\leq 1}^{\heartsuit}(S)$, given by a complexified ample class. 
Then there is a $\rB$-linear equivalence 
    \begin{align}\label{equiv:CohL}
        \IndCoh(\mathfrak{M}_S(\beta, w)^{\sigma\text{-ss}})_{-\chi_{\beta}}
        \simeq \IndL(\mathfrak{M}_S(\beta, \chi))_{w_{\beta}}^{\sigma}.
    \end{align}
\end{conj}

The left-hand side is the dg-category of ind-coherent sheaves of the moduli stack of semistable sheaves. 
The key ingredient is the notion of (ind-)\textit{limit categories} 
for the moduli stack of \textit{all} pure one-dimensional 
sheaves, 
which appear in the right-hand side of (\ref{equiv:CohL}). 
The stack $\mathfrak{M}_S(\beta, \chi)$ belongs to the class of stacks called \textit{symmetric stacks};
a general theory of limit categories for such stacks is developed in~\cite{TodaLimits}. 
We refer to Conjecture~\ref{intro:conj}, Conjecture~\ref{conj:CY3} for more details on the notation. 

In the case of $S=\mathrm{Tot}_C(\Omega_C)$ for a smooth projective curve $C$, the above conjecture is equivalent to the \textit{Dolbeault geometric Langlands conjecture} for $\GL_r$-Higgs bundles on curves proposed in~\cite{PTlim}. 
Our conjecture gives rise to new and interesting open problems, many of which are more tractable than the original Dolbeault geometric Langlands correspondence. For example, it includes a version of the Dolbeault geometric Langlands conjecture for $L$-twisted Higgs bundles, where $L$ is a sufficiently ample line bundle on a smooth projective curve $C$. Using the results in this paper together with~\cite{TodaLimits}, we will prove this conjecture for $\GL_2$ over the entire Hitchin base in~\cite{PTtwist}.

Conjecture~\ref{intro:conj} is motivated by a prospective Dolbeault geometric Langlands-type duality for a 
(yet to be defined) \textit{categorical Donaldson--Thomas theory} on Calabi--Yau threefolds~\cite{T, HHR1, HHR2}, together with the phenomenon of categorical $\chi$-independence. Although our constructions are carried out on a smooth surface $S$, the underlying geometric intuition comes from the associated local Calabi--Yau threefold $X=\mathrm{Tot}_S(\omega_S)$, see Subsection~\ref{subsec:intromotiv}.

\subsection{Background of limit categories and Dolbeault geometric Langlands conjecture}
The notion of 
limit categories was introduced in our joint work with Tudor P\u{a}durariu~\cite{PTlim} for derived stacks with self-dual cotangent complexes. We used it to formulate a version of the Dolbeault geometric Langlands conjecture for Higgs bundles by Donagi--Pantev~\cite{DPlecture}. 
In the case of $\GL_r$-Higgs bundles on curves, it is of the form 
\begin{align}\label{intro:DGLC}
    \IndCoh(\Hig_{\GL_r}(w)^{\mathrm{ss}})_{-\chi} \simeq \IndL(\Hig_{\GL_r}(\chi))_w
\end{align}
where $\Hig_{\GL_r}(\chi)$ stands for the component corresponding to Higgs bundles of degree $\chi$, 
and $(-)_w$ means the central weight $w$ part. 
The right-hand side is the limit category, which may be regarded as an effective version of classical limits of the dg-category of D-modules on $\Bun_{\GL_r}$. 
An equivalence (\ref{intro:DGLC}) is regarded as a classical limit of 
the de Rham geometric Langlands correspondence~\cite{BD0, AG, GLC1, GLC2, GLC3, GLC4, GLC5}.

Moreover, both sides of (\ref{intro:DGLC}) are compactly generated and the subcategories of compact objects admit semiorthogonal decompositions into \textit{quasi-BPS categories}~\cite{PThiggs},
which are a categorification of BPS invariants in Donaldson--Thomas theory~\cite{DM, BDIKP}. 
A striking feature of the above correspondence is that it exchanges the degree and weight: 
\begin{align*}(\chi, w) \mapsto (w, -\chi).
\end{align*}
Combined with the semiorthogonal decomposition into 
quasi-BPS categories, it yields a categorical analogue of the \textit{$\chi$-independence phenomenon}~\cite{MT, TodGV, 
DMJS, KinjoKoseki, DHKSV, GWZ}, 
one of the curious phenomena in enumerative geometry. 

By the spectral correspondence, a $\GL_r$-Higgs bundle
can be regarded as a compactly supported pure one-dimensional
coherent sheaf on the smooth surface $\mathrm{Tot}_C(\Omega_C)$.
In this paper, we extend the above picture from $\mathrm{Tot}_C(\Omega_C)$ to a more general
smooth quasi-projective surface $S$ with nef anti-canonical divisor.
Unlike the case of $\mathrm{Tot}_C(\Omega_C)$, the moduli stack
of compactly supported pure one-dimensional sheaves on $S$
may not have a self-dual cotangent complex.
Nevertheless, it belongs to the class of \textit{symmetric stacks}
considered in~\cite{TodaLimits}.

In~\cite{TodaLimits}, 
for each symmetric stack $\mathfrak{M}$ and $\delta \in \mathrm{Pic}(\mathfrak{M})_{\mathbb{Q}}^{\mathrm{CL}}$ (see~Subsection~\ref{subsec:symmetric} for the notation), we 
introduced the 
dg-category 
\begin{align*}
    \IndL(\mathfrak{M})_{\delta} \subset \IndCoh(\mathfrak{M})
\end{align*}
whose locally compact objects are given by the weight condition 
with respect to (regularizations of) maps $\nu \colon \bgm \to \mathfrak{M}$, see Subsection~\ref{subsec:limit} for details. 
The limit category is defined to be its subcategory of compact objects 
\begin{align}\label{intro:LM}
\LL(\mathfrak{M})_{\delta}
:=(\IndL(\mathfrak{M})_{\delta})^{\mathrm{cp}}.
\end{align}

In the present paper, we apply
\cite[Theorem~4.24]{TodaLimits} to the limit categories
of moduli stacks of pure one-dimensional sheaves
on smooth surfaces, and describe their semiorthogonal decompositions 
into quasi-BPS categories. We then compare them with those for 
coherent sheaves on the semistable locus established in~\cite{BPT}. 

\subsection{Limit categories of one-dimensional sheaves on surfaces}
Let $S$ be a smooth quasi-projective surface over $k$. 
For a curve class $\beta$ and an integer $\chi \in \mathbb{Z}$,
we denote by
\begin{align*}
    \mathfrak{M}_S(\beta, \chi) \to \rB
\end{align*}
the derived moduli stack of compactly supported pure
one-dimensional coherent sheaves $E$ on $S$ with
$([E], \chi(E))=(\beta, \chi)$, together with the
(derived) Hilbert--Chow morphism.
Note that $\mathfrak{M}_S(\beta, \chi)$ is not
quasi-compact in general, as it may have infinitely many
Harder--Narasimhan strata.


On the other hand, consider a 
Bridgeland stability condition~\cite{Brs1}
determined by a complexified ample class of $S$
\begin{align*}
    \sigma=(Z, \Coh_{\leq 1}^{\heartsuit}(S)), \ Z \colon K(\Coh_{\leq 1}(S)) \to \mathbb{C}.
\end{align*}
Here $\Coh_{\leq 1}^{\heartsuit}(S)$ is the abelian category of coherent sheaves on 
$S$ with at most one-dimensional supports. 
We will use the stability condition in two different ways: 
\begin{enumerate}[leftmargin=1.5em]
    \item As usual, we use $\sigma$ to define the open substack of 
    semistable objects 
   \begin{align}\label{intro:openss}
    \mathfrak{M}_S(\beta, \chi)^{\sigma\text{-ss}}\subset \mathfrak{M}_S(\beta, \chi),
\end{align}
and the compactly generated dg-category of ind-coherent sheaves~\cite{MR3136100, MR3037900}
\begin{align*}
    \IndCoh(\mathfrak{M}_S(\beta, \chi)^{\sigma\text{-ss}}).
\end{align*}
\item Together with an additional datum $w\in \mathbb{Z}$, we use $\sigma$ to define
\begin{align}\label{correct:delta}
   \delta= \delta_N(D_{\sigma, \beta, w}) \in \mathrm{Pic}(\mathfrak{M}_S(\beta, \chi))_{\mathbb{Q}}^{\mathrm{CL}},
\end{align}
and the correct twisting
of the ind-limit category
\begin{align}\label{intro:indl}
    \IndL(\mathfrak{M}_S(\beta, \chi))_{w}^{\sigma}:=
    \IndL(\mathfrak{M}_S(\beta, \chi))_{\delta_N(D_{\sigma, \beta, w})}.
\end{align}
\end{enumerate}

The construction of (\ref{correct:delta}) is based on a computation of the $\mathbb{G}_m$-weights of line bundles 
arising from the Deligne pairing; see Proposition~\ref{prop:line}.
As an application of \cite[Theorem~4.24]{TodaLimits}, we obtain the following result:  
\begin{thm}\emph{(Theorem~\ref{thm:sodLS})}\label{intro:thmS}
Suppose that $-K_S$ is nef. 
Then (\ref{intro:indl}) is compactly generated
by the subcategory of compact objects 
\begin{align}\label{intro:sod:2}
    \LL(\mathfrak{M}_S(\beta, \chi))_{w_{\beta}}^{\sigma}
    :=(\IndL(\mathfrak{M}_S(\beta, \chi))_{w_{\beta}}^{\sigma})^{\mathrm{cp}},
\end{align}
which admits the semiorthogonal decomposition of the form
\begin{align}\label{intro:sod:3}
    \LL(\mathfrak{M}_S(\beta, \chi))_{w_{\beta}}^{\sigma}
    =\left\langle \bigotimes_{i=1}^k\mathbb{T}^{+}_S(\beta_i, \chi_i)_{w_{i, \beta_i}}^{\sigma}\right\rangle.
\end{align}
Here the right-hand side ranges over the partitions 
\begin{align*}(\beta, \chi, w)=(\beta_1, \chi_1, w_1)+\cdots+(\beta_k, \chi_k, w_k)
\end{align*}
such that 
\begin{align}\label{intro:mu}
    \mu^{\sigma}(\beta, w)=\mu^{\sigma}(\beta_i, w_i), \ 
    \mu^{\sigma}(\beta_1, \chi_1)>\cdots>\mu^{\sigma}(\beta_k, \chi_k).
\end{align}
Here $\mu^{\sigma}(-) \in \mathbb{Q}$ is the slope with respect to $\sigma$. 
The transformation $(w, w_i) \mapsto (w_{\beta}, w_{i, \beta_i})$ is given by
\begin{align}\label{intro:trans}
    (-)_{\beta}=(-)+\frac{1}{2}\beta(K_S+\beta).
\end{align}
The semiorthogonal decomposition~(\ref{intro:sod:3}) is ordered by decreasing values of 
 \begin{align}\label{intro:obeta}
        o(\{(\beta_i, \chi_i)\}_{1\leq i\leq k}):=\sum_{i=1}^k \operatorname{Im}Z(\beta_i, \chi_i) \cdot\mu^{\sigma}(\beta_i, \chi_i)^2.
    \end{align}
\end{thm}

The semiorthogonal summands in the above theorem are the (positive)
\textit{quasi-BPS categories}. Here we distinguish the following two versions, see Definition~\ref{def:qbpspm} for details:
\begin{align*}
    &\mathbb{T}^{+}_S(\beta, \chi)_{w}^{\sigma} :=\LL(\mathfrak{M}_S(\beta, \chi)^{\sigma\text{-ss}})_{\delta_N(D_{\sigma, \beta, w})}, \\
&\mathbb{T}^{-}_S(\beta, \chi)_{w}^{\sigma} :=\LL(\mathfrak{M}_S(\beta, \chi)^{\sigma\text{-ss}})_{\delta_N(-D_{\sigma, \beta, -w})}. 
\end{align*}
We write $\mathbb T_S=\mathbb T_S^+$ when there is no sign.
The two versions agree for generic stability conditions by Lemma~\ref{lem:generic}.
The quasi-BPS categories for generic stability conditions were introduced and studied for K3 surfaces~\cite{PTK3} and Higgs bundles~\cite{PThiggs}; they categorify BPS invariants for the non-compact Calabi--Yau 3-fold
\begin{align*}X=\mathrm{Tot}_S(\omega_S).
\end{align*}
\subsection{The Dolbeault geometric Langlands-type conjecture for smooth surfaces}\label{subsec:introconj}
The semiorthogonal decomposition of~\cite{BPT}, recalled in Theorem~\ref{thm:bpt},
applied to the moduli stack of semistable objects~(\ref{intro:openss}) yields
\begin{align}\label{intro:sod:1}
        \Coh(\mathfrak{M}_S(\beta, w)^{\sigma\text{-ss}})_{-\chi_{\beta}}
        =\left\langle \bigotimes_{i=1}^k \mathbb{T}^{-}_S(\beta_i, w_i)_{-\chi_{i, \beta_i}}^{\sigma}
        \right\rangle.
    \end{align}
    Here the right-hand side ranges over the partitions satisfying the same condition~(\ref{intro:mu}), 
    and the transformation $(\chi, \chi_i)\mapsto (\chi_{\beta}, \chi_{i, \beta_i})$ is given by~(\ref{intro:trans}). The order of the 
    semiorthogonal decomposition is again determined by decreasing values of~(\ref{intro:obeta}). 

\begin{figure}[!hbp]\caption{Conceptual overview}\label{figure:T}
\[
\begin{tikzcd}[
  ampersand replacement=\&,
  column sep=1.2em,
  row sep=2.7em,
  cells={nodes={inner sep=2pt}}
]
\&
|[draw, rounded corners, inner sep=5pt]|
{\begin{tabular}{c}
$(S,\sigma)$\\[-1pt]
{\scriptsize surface with}\\[-2pt]
{\scriptsize stability condition}
\end{tabular}}
\arrow[dl, dashed]
\arrow[dr, dashed]
\&
\\[0.3em]
|[draw, rounded corners, inner sep=5pt]|
{\begin{tabular}{c}
$\Coh\!\bigl(\mathfrak M_S(\beta,w)^{\sigma\text{-ss}}\bigr)_{-\chi_{\beta}}$\\[1pt]
{\scriptsize coherent sheaves on}\\[-2pt]
{\scriptsize semistable moduli stacks}
\end{tabular}}
\arrow[d, equal]
\&
\conjeq
\arrow[d, Rightarrow]
\&
|[draw, rounded corners, inner sep=5pt]|
{\begin{tabular}{c}
$\LL(\mathfrak M_S(\beta,\chi))_{w_{\beta}}^{\sigma}$\\[1pt]
{\scriptsize limit category on}\\[-2pt]
{\scriptsize the full moduli stack}
\end{tabular}}
\arrow[d, equal]
\\
|[draw, rounded corners, inner sep=5pt]|
{\begin{tabular}{c}
$\left\langle
  \bigotimes_{i=1}^k
  \mathbb{T}^{-}_S(\beta_i,w_i)_{-\chi_{i, \beta_i}}^{\sigma}
\right\rangle$\\[1pt]
{\scriptsize combinatorics of}\\[-2pt]
{\scriptsize weight polytope}
\end{tabular}}
\&
\conjeq
\arrow[d, Rightarrow]
\&
|[draw, rounded corners, inner sep=5pt]|
{\begin{tabular}{c}
$\left\langle
  \bigotimes_{i=1}^k
  \mathbb{T}^{+}_S(\beta_i,\chi_i)_{w_{i, \beta_i}}^{\sigma}
\right\rangle$\\[1pt]
{\scriptsize Harder--Narasimhan}\\[-2pt]
{\scriptsize stratification}
\end{tabular}}
\\
\&
|[draw, rounded corners, inner sep=5pt]|
{\begin{tabular}{c}
$\mathbb{T}^{-}_S(\beta,w)_{-\chi_{\beta}}^{\sigma}
 \conjeq
  \mathbb{T}^{+}_S(\beta,\chi)_{w_{\beta}}^{\sigma}$\\[1pt]
{\scriptsize duality of}\\[-2pt]
{\scriptsize quasi-BPS categories}
\end{tabular}}
\arrow[d, Rightarrow]
\&
\\[-0.2em]
\&
\textit{Categorical $\chi$-independence}
\&
\end{tikzcd}
\]
\end{figure}
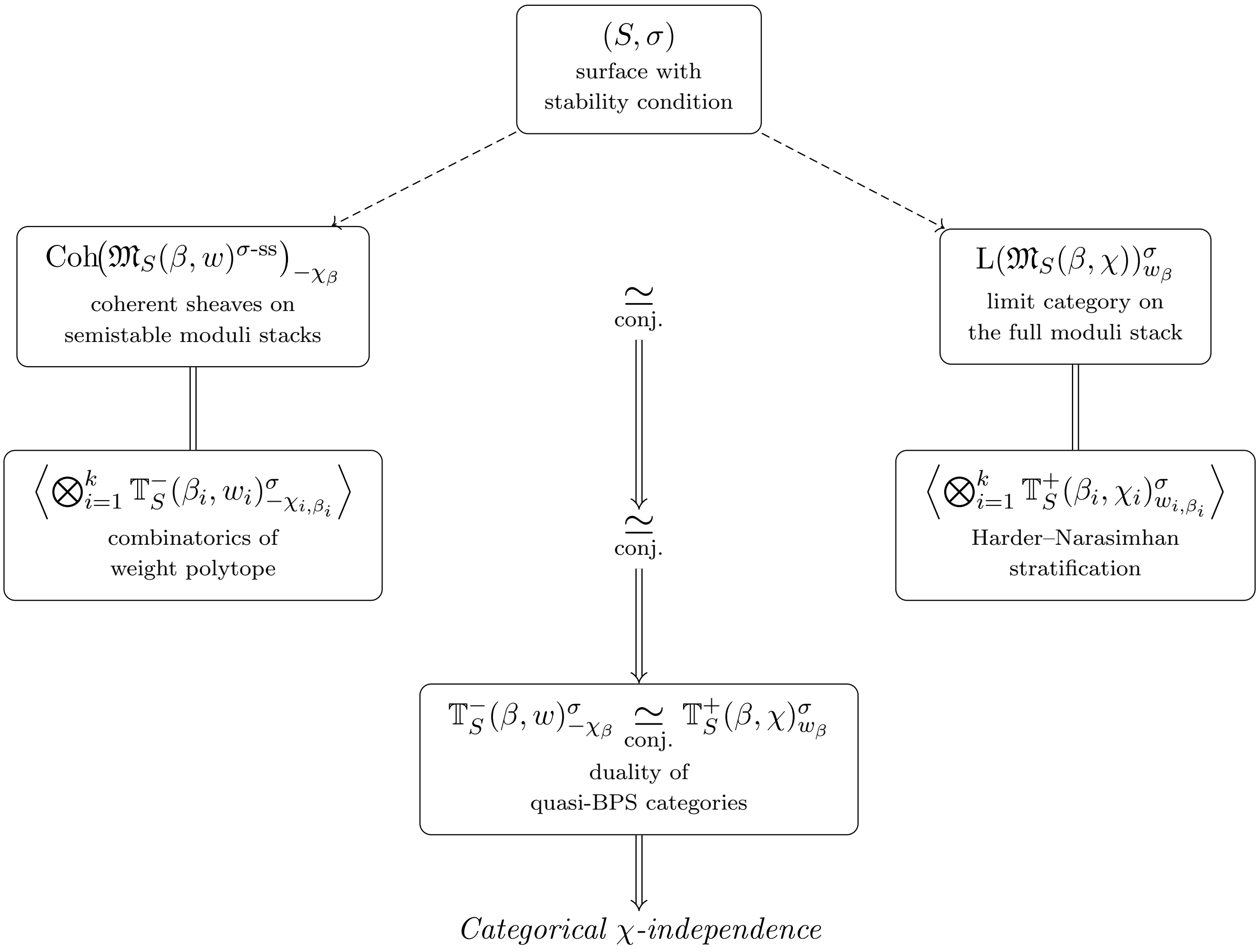

Based on the above comparison of the semiorthogonal decompositions~(\ref{intro:sod:3}) and~(\ref{intro:sod:1}), we propose the following conjecture:
\begin{conj}\label{intro:conj}
Suppose that $-K_S$ is nef. 
    There is a $\rB$-linear equivalence 
    \begin{align}\label{equiv:CohL2}
        \Coh(\mathfrak{M}_S(\beta, w)^{\sigma\text{-ss}})_{-\chi_{\beta}}
        \simeq \LL(\mathfrak{M}_S(\beta, \chi))_{w_{\beta}}^{\sigma}.
    \end{align}
    Moreover, the above equivalence restricts to the following duality between the quasi-BPS categories: 
    \begin{align}\label{intro:equiv:T}
        \mathbb{T}^{-}_S(\beta, w)_{-\chi_{\beta}}^{\sigma} \simeq \mathbb{T}^{+}_S(\beta, \chi)_{w_{\beta}}^{\sigma}.
    \end{align}
\end{conj}
Note that, because of the compact generation of both sides of (\ref{equiv:CohL}), an equivalence (\ref{equiv:CohL}) is equivalent to (\ref{equiv:CohL2}). 
The conceptual overview of the above conjecture is described in Figure~\ref{figure:T}.

\subsection{Examples}\label{intro:exam}
In Section~\ref{sec:examples}, we will discuss Conjecture~\ref{intro:conj} in several cases; 
\begin{enumerate}[leftmargin=1.5em]
\item Over the locus $\rB^{\circ} \subset \rB$ corresponding to integral curves, Conjecture~\ref{intro:conj} follows from the autoduality of 
compactified Jacobians of integral planar curves~\cite{Ardual}, see Proposition~\ref{prop:integral}.
    \item $S=\mathrm{Tot}_C(L)$ for a line bundle $L$ on a smooth projective curve $C$
    with $L=\Omega_C$ or $\deg L>2g(C)-2$; in the former case, 
    Conjecture~\ref{intro:conj} recovers the \textit{Dolbeault geometric Langlands conjecture} for $G=\GL_r$ proposed in~\cite{PTlim}. The latter case gives its generalization to 
    $L$-twisted $\GL_r$-Higgs bundles, which we prove in~\cite{PTtwist} for $\GL_2$ under the stronger assumption $\deg L>2g(C)$. 
    \item $S$ is a K3 surface: in this case, Conjecture~\ref{intro:conj} gives a version of 
    Dolbeault geometric Langlands-type conjecture for \textit{Beauville--Mukai systems}, which 
    are (when there are no strictly semistable sheaves) Lagrangian fibrations of irreducible holomorphic symplectic manifolds. 
    \item Other cases, such as $S=\mathbb{P}^2$, blow-up of $\mathbb{P}^1 \times \mathbb{A}^1$. 
\end{enumerate}

\subsection{Hecke operators}
Using the functoriality results for symmetric stacks 
established in~\cite{TodaLimits}, we construct categorical Hecke operators on limit categories. Combining this with Theorem~\ref{intro:thmS}, we also obtain Hecke operators on quasi-BPS categories, which, for symplectic surfaces, are expected to lift the Hecke operators on BPS cohomology in~\cite{DHKSV}.

Let $\mathfrak{M}_S^0(d)$ be the derived moduli stack of zero-dimensional sheaves on $S$ of length $d$. The direct sum 
\begin{align*}
    \mathbb{H}:=\bigoplus_{d\in \mathbb{Z}_{\geq 0}} \LL(\mathfrak{M}_S^0(d))_0
\end{align*}
admits a monoidal structure given by the categorical Hall product. 
We prove the following result, which generalizes the construction of Hecke operators for $\GL_r$-Higgs bundles in~\cite[Section~9]{PTlim}: 
\begin{thm}\emph{(Theorem~\ref{thm:hecke1})}\label{intro:thmhecke}
Suppose that $-K_S$ is nef. For fixed $\beta$ and a $\mathbb{Q}$-divisor class $D$ on $S$, 
   the direct sum 
    \begin{align*}
        \bigoplus_{\chi\in \mathbb{Z}} \LL(\mathfrak{M}_{S}(\beta, \chi))_{\delta_N(D)}
    \end{align*}
    admits right and left actions of $\mathbb{H}$.
\end{thm}

As a corollary, we have the following: 
\begin{cor}\emph{(Corollary~\ref{cor:hecke2})}\label{intro:corhecke}
Suppose that $-K_S$ is nef and $(\beta, w)$ is $\sigma$-primitive. 
Then the direct sum of quasi-BPS categories
\begin{align*}
\bigoplus_{\chi\in \mathbb{Z}}\mathbb{T}^{+}_S(\beta, \chi)_{w_{\beta}}^{\sigma}
\end{align*}
admits right and left actions of $\mathbb{H}$. 
\end{cor}

The above Hecke actions are expected to be compatible with Wilson 
actions on $\Coh(\mathfrak{M}_S(\beta, w)^{\sigma\text{-ss}})$ under the equivalences in Conjecture~\ref{intro:conj}, see Conjecture~\ref{conj:WH}.
If $S$ is a symplectic surface, Corollary~\ref{cor:hecke2} may be regarded as a categorical lift of the Hecke operators on BPS cohomology in~\cite{DHKSV}; see Remark~\ref{rmk:heckelift}.

\subsection{Motivation from categorical Donaldson--Thomas theory}\label{subsec:intromotiv}
Conjecture~\ref{intro:conj} is motivated by categorical Donaldson--Thomas theory and categorical $\chi$-independence phenomena. 

Let $X$ be a smooth Calabi--Yau 3-fold with a Bridgeland stability condition 
\begin{align*}
    \sigma=(Z, \Coh_{\leq 1}^{\heartsuit}(X)), \ Z \colon K(\Coh_{\leq 1}(X)) \to \mathbb{C}.
\end{align*}
Let $\mathfrak{M}_X$ be the derived moduli stack of compactly supported pure one-dimensional 
coherent sheaves on $X$, and we denote by 
\begin{align*}
 \mathfrak{M}_X^{\sigma\text{-ss}} \subset \mathfrak{M}_X
\end{align*}
its open substack of $\sigma$-semistable sheaves. 
We expect the following \textit{Dolbeault geometric 
Langlands-type `conjecture'}, which extends Conjecture~\ref{intro:conj} to Calabi--Yau 3-folds:  
\begin{conj2}\label{conj2:M}
    There is an equivalence of the form 
    \begin{align*}
        \mathrm{IndDT}(\mathfrak{M}_X^{\sigma\text{-ss}}) \simeq \IndL(\mathfrak{M}_X)^{\sigma}.
    \end{align*}
    Moreover, the 
    above equivalence restricts to an equivalence of quasi-BPS categories. 
\end{conj2}
Here we use quotation marks around `conjecture' because neither side 
nor the `quasi-BPS categories' are yet defined 
in general. Roughly speaking, the two sides should be as follows: 
\begin{itemize}[leftmargin=1.5em]
    \item The dg-category $\mathrm{IndDT}(\mathfrak{M}_X^{\sigma\text{-ss}})$ should be 
obtained by gluing locally defined categories of matrix factorizations and should 
categorify the (cohomological) Donaldson--Thomas invariants counting one-dimensional 
semistable sheaves. 
We refer to~\cite{HHR1, HHR2} for recent progress in this direction. 
\item The dg-category $\IndL(\mathfrak{M}_X)^{\sigma}$
should be obtained as a subcategory of 
the DT category $\mathrm{IndDT}(\mathfrak{M}_X)$, defined without imposing a 
stability condition but instead by a 
certain weight condition with respect to maps $\bgm \to \mathfrak{M}_X$.   
\item We expect that both sides are compactly generated and that the subcategories of compact objects admit semiorthogonal decompositions in terms of quasi-BPS categories for Calabi--Yau 3-folds, analogous to~(\ref{intro:sod:1}) and~(\ref{intro:sod:3}). An equivalence is expected to preserve 
these semiorthogonal decompositions and may yield a duality between quasi-BPS categories. 
\end{itemize}

By categorical dimensional reduction (or Koszul duality)~\cite{I, MR2982435, Hirano, T}, 
Conjecture~\ref{intro:conj} is expected to give a rigorous formulation 
of `Conjecture~\ref{conj2:M}' for the non-compact Calabi--Yau 3-fold  
\begin{align*}X=\mathrm{Tot}_S(\omega_S)
\end{align*}
when $-K_S$ is nef. For a more general surface $S$ (without the assumption that $-K_S$ is nef), we expect a formulation of the above `conjecture' for the local surface $X$ using the singular support quotient construction as in~\cite{T}. 
In Section~\ref{sec:DLCY3}, we formulate an expected equivalence 
of the form 
\begin{align}\notag
  \mathrm{IndDT}^{\mathbb{G}_m}(\mathfrak{M}_X(\beta, w)^{\sigma_X\text{-ss}})_{-\chi_{\beta}}
  \simeq 
    \IndL^{\mathbb{G}_m}(\mathfrak{M}_X(\beta, \chi))_{w_{\beta}}^{\sigma_X},
\end{align}
see Conjecture~\ref{conj:local-surfaces} for details. 

The conjectural duality of quasi-BPS categories should be relevant to 
\textit{categorical $\chi$-independence}. 
The $\chi$-independence phenomenon is a
curious feature of the enumerative geometry of Calabi--Yau 3-folds: 
the BPS invariants (or Gopakumar--Vafa invariants) counting one-dimensional semistable sheaves $E$
with $([E], \chi(E))=(\beta, \chi)$ are conjectured to be independent of $\chi$. This conjecture is important for providing a mathematical foundation for Gopakumar--Vafa invariants; see~\cite[Conjecture~3.20]{MT} and~\cite[Conjecture~1.2]{TodGV} for its formulation, and~\cite{DMJS, KinjoKoseki, DHKSV, GWZ} for recent progress.

We expect a categorical analogue of $\chi$-independence, 
and the conjectural duality between quasi-BPS categories may be relevant to it. Indeed, the equivalence~(\ref{intro:equiv:T}) swaps the Euler characteristic and the weight, and hence interchanges $\chi$-independence with \textit{weight-independence}; in some situations 
the latter is automatic, as in the case of Higgs bundles~\cite{PThiggs}.
For K3 surfaces, categorical $\chi$-independence also yields an understanding of twisted non-commutative crepant resolutions of symplectic singularities; 
see~\cite{PTK3} and Subsection~\ref{subsec:K3}.

On the other hand, we do not expect to obtain categorical $\chi$-independence via Hecke operators as in~\cite{DHKSV}; instead, the 
Hecke operators in Theorem~\ref{intro:thmhecke} are expected to play an important role in the reduction of Conjecture~\ref{intro:conj} to a Whittaker normalization, as in~\cite{TodaGL2, TodaA, PTtwist}.

\subsection{Acknowledgements}
The author would like to thank Tudor P{\u a}durariu, Tatsuki Kinjo,
Francesco Sala, Duiliu Diaconescu, Mauro Porta, Philasang Yoo,
and David Fang for helpful discussions.

The final stages of this work were carried out while the author was visiting the Morningside Center of Mathematics, Chinese Academy of Sciences, as a Distinguished Visiting Professor from October 18 to November 18, 2026. The author thanks the Morningside Center of Mathematics for its hospitality.

The author is supported by the World Premier International Research
Center Initiative (WPI Initiative), MEXT, Japan, the Inamori Research
Institute for Science, and JSPS KAKENHI Grant Number JP24H00180.

\subsection*{AI disclosure}
ChatGPT-6 Astra was used to explore a proof of Lemma~\ref{lem:HN}, to check the proofs in the paper, and to assist with editing. The author takes full responsibility for the correctness and content of the paper.

\subsection{Notation and conventions}\label{subsec:notation0}
\subsubsection*{Conventions on stacks}
All (derived) stacks are defined over an algebraically closed field $k$ of characteristic zero.
We write $\X(k)$ for the $k$-valued points of a stack $\X$, and
$\X^{\mathrm{cl}}$ for its classical truncation. 
A stack $\X$ is called \textit{QCA} if it is quasi-compact with affine
geometric stabilizers and its classical inertia stack is of finite
presentation over $\X^{\mathrm{cl}}$; see~\cite[Definition~1.1.8]{MR3037900}.
It is \textit{locally QCA} if it admits a Zariski open cover by QCA substacks.

For a morphism $f\colon\X\to\Y$, we denote its relative cotangent and
tangent complexes by $\mathbb{L}_f$ and $\mathbb{T}_f$.
The morphism is \textit{quasi-smooth} if $\mathbb{L}_f$ is perfect of
cohomological amplitude $[-1,1]$; we then set
$\omega_f:=\det\mathbb{L}_f$ and $\dim f:=\operatorname{rank}\mathbb{L}_f$.
All fiber products are derived. For a vector space $V$, we write
$V[-1]:=\Spec\mathrm{Sym}^*(V^{\vee}[1])$.

\subsubsection*{Categories of sheaves}\label{subsec:notation}
For a dg-category $\cC$, let $\cC^{\mathrm{cp}}$ be its subcategory of
compact objects and $\Ind(\cC)$ its ind-completion.
We call $\cC$ \textit{compactly generated} if
$\cC=\Ind(\cC^{\mathrm{cp}})$.
For compactly generated dg-categories $\cC_1,\cC_2$, we denote their tensor
product by $\cC_1\otimes\cC_2$ and set
\begin{align*}
 \cC_1^{\mathrm{cp}}\otimes\cC_2^{\mathrm{cp}}
 :=(\cC_1\otimes\cC_2)^{\mathrm{cp}}.
\end{align*}
The latter is split-generated by exterior products $c_1\boxtimes c_2$
with $c_i\in\cC_i^{\mathrm{cp}}$.

For a derived stack $\X$, we write $\Coh(\X)$ and $\IndCoh(\X)$ for the
dg-categories of coherent and ind-coherent sheaves, following~\cite{MR3701352},
and $\Coh^{\heartsuit}(\X)$ for the heart of the standard t-structure.
For a QCA stack, $\IndCoh(\X)$ is compactly generated with compact objects
$\Coh(\X)$~\cite{MR3037900}.
For a quasi-smooth stack $\mathfrak{M}$, its derived dual is
\begin{align*}
 \mathbb{D}_{\mathfrak{M}}(-):=\cHom(-,\mathcal{O}_{\mathfrak{M}})
 \colon\Coh(\mathfrak{M})^{\mathrm{op}}\stackrel{\sim}{\to}
 \Coh(\mathfrak{M}).
\end{align*}

\subsubsection*{\texorpdfstring{$\mathbb{G}_m$-weights}{Gm-weights}}\label{subsec:Gmwt}
For a $\mathbb{G}_m$-gerbe $\X\to X$ with Brauer class $\alpha$, we have
\begin{align*}
 \Coh(\X)=\bigoplus_{w\in\mathbb{Z}}\Coh(\X)_w,
 \quad \Coh(\X)_w\simeq\Coh(X,\alpha^w),
\end{align*}
where the subscript $w$ denotes weight $w$ and the right-hand side is the
category of $\alpha^w$-twisted coherent sheaves~\cite{MR2309155}.
For $A=\bigoplus_w A_w\in\Coh(\bgm)$, we set
\begin{align*}
 A^{>0}:=\bigoplus_{w>0}A_w,\quad
 A^{<0}:=\bigoplus_{w<0}A_w,\quad
 \wt(A):=\{w\in\mathbb{Z}:A_w\neq0\}.
\end{align*}
We identify $H^2(\bgm,\mathbb{Q})=\mathbb{Q}$ and use
\begin{align*}
 c_1\colon K(\bgm)_{\mathbb{Q}}\to\mathbb{Q},
 \quad c_1(A)=\sum_{w\in \mathbb{Z}} w\cdot\dim A_w,
\end{align*}
where $\dim A_w$ denotes the Euler characteristic.

\subsubsection*{Filtered and graded objects}
We recall the notion of stacks of \textit{filtered and graded objects}, following~\cite{Halpinstab,HalpK32}.

Let $\Theta:=\mathbb{A}^1/\mathbb{G}_m$, where $\mathbb{G}_m$ acts on $\mathbb{A}^1$ by weight one. For a derived stack $\mathfrak{M}$, 
the stacks of filtered and graded objects are defined by 
\begin{align*}
    \mathrm{Filt}(\mathfrak{M}):=\mathrm{Map}(\Theta, \mathfrak{M}), \ 
    \mathrm{Grad}(\mathfrak{M}):=\mathrm{Map}(\bgm, \mathfrak{M}).
\end{align*}
There are natural evaluation maps 
\begin{align}\label{dia:ev}
  \xymatrix{
   \mathrm{Grad}(\mathfrak{M})  \ar@<-0.5ex>[r]&  \ar@<-0.5ex>[l]_-{\ev_0} \mathrm{Filt}(\mathfrak{M}) \ar[r]^-{\ev_1} & \mathfrak{M}.
   }
\end{align}
A map $\nu \colon \bgm \to \mathfrak{M}$ is called a \textit{graded point}, and $f \colon \Theta \to \mathfrak{M}$ is called a \textit{filtration} of $f(1)$. For a filtration $f$, the map  $f_0:=\ev_0(f) \colon \bgm \to \mathfrak{M}$ is called its \textit{associated graded}. A filtration $f$ is called \textit{non-degenerate} if the induced map on stabilizer groups of $f_0$ is non-trivial. 

A \textit{$\Theta$-stratum} of $\mathfrak{M}$ is an open and closed substack $\mathfrak{S} \subset \mathrm{Filt}(\mathfrak{M})$
such that 
$\ev_1|_{\mathfrak{S}} \colon \mathfrak{S} \to \mathfrak{M}$
is a closed immersion. 
Its \textit{center} is an open and closed substack $\mathfrak{Z} \subset \mathrm{Grad}(\mathfrak{M})$
such that $\ev_0(\mathfrak{S})=\mathfrak{Z}$. 
There is also the notion of \textit{$\Theta$-stratification}, 
see~\cite{Halpinstab}, Subsection~\ref{subsec:thetast}.

\subsubsection*{Semiorthogonal decomposition}
For a pretriangulated dg-category $\cD$, a totally ordered set $I$  and a family of
pretriangulated full subcategories $\cD_i$ for $i\in I$, we write
\begin{align*}
 \cD=\langle\cD_i:i\in I\rangle
\end{align*}
if the $\cD_i$ generate $\cD$ and $\Hom(\cD_i,\cD_j)=0$ for $i>j$.

\section{Limit categories for symmetric stacks}\label{sec:limit}

In this section, we recall the general results on limit categories for
symmetric stacks from~\cite{TodaLimits} which will be used below.  We
include the definitions in order to fix the conventions for the weight
conditions, and state the functoriality, semiorthogonal decomposition and
K\"unneth theorems in the form needed for moduli stacks of sheaves on
surfaces.

\subsection{Symmetric stacks and the component lattice}\label{subsec:symmetric}

Let $\mathfrak{M}$ be a quasi-smooth derived stack over $k$.  It is called
\textit{symmetric} if, for every $x\in\mathfrak{M}(k)$, the
virtual representation of $G_x:=\operatorname{Aut}(x)$ given by the tangent
complex is self-dual:
\begin{align}\label{def:symmetric-review}
 [\mathbb{T}_{\mathfrak{M}}|_x]
 = [\mathbb{T}_{\mathfrak{M}}|_x^{\vee}]
 \quad\mbox{in }K(BG_x).
\end{align}
Note that this 
condition is weaker than the requirement that
$\mathbb{L}_{\mathfrak{M}}$ is self-dual.  In particular, zero-shifted symplectic stacks
and cotangent stacks are symmetric.  More importantly for this paper, the
derived moduli stack of compactly supported one-dimensional sheaves on a
smooth surface is symmetric; see Proposition~\ref{prop:sym}.

We recall the component lattice of~\cite{BHIK}.  For a lattice $\Lambda$, let
$T_{\Lambda}$ be the torus with cocharacter lattice $\Lambda$.  The formal
lattice
\begin{align*}
 \operatorname{CL}(\mathfrak{M})(\Lambda)
 :=\pi_0\operatorname{Map}(BT_{\Lambda},\mathfrak{M})
\end{align*}
is called the \textit{component lattice} of $\mathfrak{M}$.  We write
$\operatorname{CL}(\mathfrak{M})_{\mathbb{Q}}$ for its rational version.
There is a natural map
\begin{align}\label{map:deltaP-review}
 \delta_P\colon \operatorname{Pic}(\mathfrak{M})_{\mathbb{Q}}
 \to
 \operatorname{CL}(\mathfrak{M})_{\mathbb{Q}}^{\vee},
 \quad
 \delta_P(\mathcal{L})(\nu)=c_1(\nu^*\mathcal{L}).
\end{align}
Its image is denoted by
$\operatorname{Pic}(\mathfrak{M})_{\mathbb{Q}}^{\operatorname{CL}}$.
Although most parameters used below belong to this image, the definition of
a limit category makes sense for any
$\delta\in\operatorname{CL}(\mathfrak{M})_{\mathbb{Q}}^{\vee}$.

\subsection{The weight condition and limit categories}\label{subsec:limit}
We recall limit categories for symmetric stacks; 
see~\cite[Subsection~3.2]{TodaLimits} for more details. 

Let $\nu\colon \bgm\to\mathfrak{M}$ be a graded point.  A
\textit{regularization} of $\nu$ is a finite-dimensional
$\mathbb{G}_m$-representation $V$ together with a quasi-smooth map
\begin{align*}
 \nu_V^{\mathrm{reg}}\colon
 V[-1]/\mathbb{G}_m\to\mathfrak{M}
\end{align*}
whose restriction along the zero section is $\nu$.  A regularization always
exists.  Let
$\iota\colon V[-1]/\mathbb{G}_m\to\bgm$ be the projection.

Suppose first that $\mathfrak{M}$ is quasi-smooth, QCA and symmetric.  For
$\delta\in\operatorname{CL}(\mathfrak{M})_{\mathbb{Q}}^{\vee}$, the \textit{limit
category}
\begin{align}\label{def:limit-symmetric}
 \LL(\mathfrak{M})_{\delta}\subset\Coh(\mathfrak{M})
\end{align}
consists of objects $\mathcal{E}$ such that, for every
$\nu$ and its regularization as above, the following condition is 
satisfied: 
\begin{align}\label{wt:limit-symmetric}
& \wt\!\left(\iota_*(\nu_V^{\mathrm{reg}})^*\mathcal{E}\right)
 \subset {} \\
 \notag&\left[
 \frac{1}{2}c_1(\nu^*\mathbb{L}_{\mathfrak{M}}^{<0})
       +c_1((V^{\vee})^{<0}),
 \frac{1}{2}c_1(\nu^*\mathbb{L}_{\mathfrak{M}}^{>0})
       +c_1((V^{\vee})^{>0})
 \right]
 +\delta(\nu).
\end{align}
The condition is independent of the choice of regularization.  When
$\mathcal{E}$ is perfect, it is equivalent to the simpler condition
\begin{align}\label{wt:limit-perfect}
 \wt(\nu^*\mathcal{E})\subset 
 \left[
 \frac{1}{2}c_1(\nu^*\mathbb{L}_{\mathfrak{M}}^{<0}),
 \frac{1}{2}c_1(\nu^*\mathbb{L}_{\mathfrak{M}}^{>0})
 \right]+\delta(\nu).
\end{align}
If moreover $\mathfrak{M}$ admits a good moduli space, it is enough to use
graded points whose images are closed.

For a quasi-smooth, locally QCA and symmetric stack, not necessarily
quasi-compact, one sets
\begin{align}\label{def:IndL-review}
 \IndL(\mathfrak{M})_{\delta}
 :=\lim_{\mathfrak{U}\subset\mathfrak{M}}
 \Ind(\LL(\mathfrak{U})_{\delta}),
\end{align}
where $\mathfrak{U}$ runs through the QCA open substacks and the transition
functors are restrictions.  We use the notation
\begin{align}\label{def:L-Ltilde-review}
 \LL(\mathfrak{M})_{\delta}
 :=(\IndL(\mathfrak{M})_{\delta})^{\mathrm{cp}},\quad
 \widetilde{\LL}(\mathfrak{M})_{\delta}
 :=\lim_{\mathfrak{U}\subset\mathfrak{M}}
 \LL(\mathfrak{U})_{\delta}.
\end{align}
Thus $\widetilde{\LL}$ is the subcategory of \textit{locally} compact objects, whereas $\LL$ is the subcategory of \textit{globally} compact objects in
$\IndL$.  These two categories need not agree for a non-quasi-compact
stack.

\subsection{Functoriality}

We recall the functoriality result proved in~\cite{TodaLimits}.
Consider a correspondence
\begin{align}\label{dia:limit-correspondence}
 \mathfrak{M}_1\stackrel{q}{\leftarrow}
 \mathfrak{N}\stackrel{p}{\to}\mathfrak{M}_2,
\end{align}
where $\mathfrak{M}_i$ are quasi-smooth, QCA and symmetric,
$\mathfrak{N}$ is quasi-smooth, $q$ is quasi-smooth and $p$ is projective.
Suppose that for every graded point
$g\colon\bgm\to\mathfrak{N}$ one has
\begin{align}\label{wt:limit-correspondence}
 c_1(g^*\mathbb{L}_p^{>0})=c_1(g^*\mathbb{L}_q^{<0}),
 \quad
 c_1(g^*\mathbb{L}_p^{<0})=c_1(g^*\mathbb{L}_q^{>0}).
\end{align}
The following result is proved in~\cite{TodaLimits}.

\begin{prop}\emph{(\cite[Proposition~3.13]{TodaLimits})}\label{prop:limit-correspondence}
If $\delta_i\in\operatorname{CL}(\mathfrak{M}_i)_{\mathbb{Q}}^{\vee}$
satisfy
\begin{align}\label{delta:limit-correspondence}
 q^*\delta_1=p^*\delta_2+\frac{1}{2}\delta_P(\det\mathbb{L}_p),
\end{align}
then the functor $p_{*}q^*$ restricts to
\begin{align*}
 p_*q^*\colon
 \LL(\mathfrak{M}_1)_{\delta_1}\to
 \LL(\mathfrak{M}_2)_{\delta_2}.
\end{align*}
The analogous statement holds for locally QCA stacks with $\LL$ replaced
by $\widetilde{\LL}$.
\end{prop}
\begin{remark}
There are two variants which will be used implicitly.  Replacing (\ref{delta:limit-correspondence})
with \begin{align}\notag
 q^*\delta_1=p^*\delta_2-\frac{1}{2}\delta_P(\det\mathbb{L}_p),
\end{align}
gives a functor 
\begin{align*}
 p_*(q^*(-) \otimes \omega_q)\colon
 \LL(\mathfrak{M}_1)_{\delta_1}\to
 \LL(\mathfrak{M}_2)_{\delta_2}.
\end{align*} If $q^*\delta_1=p^*\delta_2$ and
$\det\mathbb{L}_p^{1/2}$ exists, then
\begin{align*}
 p_*\bigl(q^*(-)\otimes\det\mathbb{L}_p^{1/2}\bigr)
 \colon\LL(\mathfrak{M}_1)_{\delta_1}\rightarrow
 \LL(\mathfrak{M}_2)_{\delta_2}.
\end{align*}
If both $p$ and $q$ are quasi-smooth and projective, the functor in
Proposition~\ref{prop:limit-correspondence} admits a right adjoint, see
\cite[Corollary~3.16]{TodaLimits}:
\begin{align}\label{adjoint:limit-correspondence}
 q_*p^!\colon
 \LL(\mathfrak{M}_2)_{\delta_2}\rightarrow
 \LL(\mathfrak{M}_1)_{\delta_1},
 \quad
 p^!(-)=p^*(-)\otimes\omega_p[\dim p].
\end{align}
\end{remark}

\subsection{Semiorthogonal decompositions from \texorpdfstring{$\Theta$}{Theta}-stratification}
\label{subsec:thetast}
Let $\ell\in\operatorname{CL}(\mathfrak{M})_{\mathbb{Q}}^{\vee}$ and let
$b$ be a positive definite rational quadratic norm on the component
lattice; see~\cite[Subsection~2.5]{TodaLimits}.  These data define the numerical invariant
\begin{align*}
 \mu(f)=-\frac{\ell(f_0)}{\sqrt{b(f_0)}}
\end{align*}
for a non-degenerate filtration $f\colon\Theta\to\mathfrak{M}$.  Suppose
that it determines a locally finite $\Theta$-stratification
(cf.~\cite[Definition~2.3]{TodaLimits})
\begin{align}\label{theta:limit-review}
 \mathfrak{M}
 =\mathfrak{M}^{\ell\text{-ss}}
 \sqcup\bigsqcup_{\mu\in I_{>0}}\mathfrak{S}_{\mu},
\end{align}
indexed by a locally finite set
$I\subset\mathbb{R}_{\geq 0}$ containing $0$, 
whose cumulative unions form an open cover of $\mathfrak{M}$. 
This means that 
\begin{enumerate}
\item For $\mu \preceq \mu'$, we have 
    $\mathfrak{M}_{\preceq \mu} \subset \mathfrak{M}_{\preceq \mu'}$, and we have the open cover 
    \begin{align*}
 \mathfrak{M}=\bigcup_{\mu\in I}\mathfrak{M}_{\preceq\mu}.
    \end{align*}
    \item There is a $\Theta$-stratum 
    $p_{\mu} \colon \mathfrak{S}_{\mu} \hookrightarrow \mathfrak{M}_{\preceq \mu}$
    with center $q_{\mu} \colon \mathfrak{S}_{\mu} \to \mathfrak{Z}_{\mu}$
such that 
\begin{align*}
    \mathfrak{M}_{\prec \mu}:=\bigcup_{\mu'\prec \mu} \mathfrak{M}_{\preceq \mu'}=\mathfrak{M}_{\preceq \mu} \setminus \mathfrak{S}_{\mu}.
\end{align*}
    \item The open substack $\mathfrak{M}_{\preceq 0}$ is
    precisely the locus of $\ell$-semistable points.
    We denote it by $\mathfrak{M}^{\ell\text{-ss}}$.
    \item For every 
    $c\in I\setminus\{0\}$ and point 
    $f\in\mathfrak{S}_c(k)$ regarded as a map
    $f\colon\Theta\to\mathfrak{M}$, the graded point $f_0$
    is non-degenerate and $\mu(f)=c$.
    Moreover, for every map $g\colon\Theta\to\mathfrak{M}$
    with non-degenerate $g_0$ and $g(1)\cong f(1)$, we have
    $\mu(g)\leq c$. If $\mu(g)=c$, there exist positive integers
    $m,n$ such that
  $
        f\circ(-)^m \cong g\circ(-)^n$, 
        where $(-)^r \colon \Theta \to \Theta$ is the $r$-th power map. 
\end{enumerate}

For each $\mu \in I_{>0}$, we have the 
derived stack 
with $\Theta$-stratum
\begin{align}\label{dtheta}
    (\mathfrak{M}_{\preceq \mu}, \mathfrak{S}_{\mu}, \mathfrak{Z}_{\mu}).
\end{align}
We assume the following condition
(cf.~\cite[Assumption~4.1]{TodaLimits}):
\begin{assum}\label{assum:dag}
For each $\mu \in I_{>0}$, the derived stack with $\Theta$-stratum 
(\ref{dtheta}) satisfies the following conditions: 
\begin{enumerate}
\item 
The derived stack 
$\mathfrak{M}_{\preceq \mu}$ is quasi-smooth, QCA, locally of finite presentation, has a separated diagonal, and is symmetric; 
\item 
For any $z\in \mathfrak{Z}_{\mu}(k)$ corresponding to 
$\nu \colon \bgm \to \mathfrak{M}_{\preceq\mu}$, 
the $\mathbb{G}_m$-weights of $\cH^{-1}(\nu^* \mathbb{L}_{\mathfrak{M}})$ are non-negative;
\item The center $\mathfrak{Z}_{\mu}$ of the $\Theta$-stratum $\mathfrak{S}_{\mu}$ admits 
a good moduli space.
\end{enumerate}
\end{assum}
The following is the main result of~\cite{TodaLimits}:
\begin{thm}[{\cite[Theorem~4.24]{TodaLimits}}]
\label{thm:limit-sod-general}
Under Assumption~\ref{assum:dag}, $\IndL(\mathfrak{M})_{\delta}$ is compactly
generated.  Its compact objects admit a semiorthogonal decomposition
\begin{align}\label{sod:limit-general}
 \LL(\mathfrak{M})_{\delta}
 =\left\langle
 \Upsilon_{\mu}\LL(\mathfrak{Z}_{\mu})_{\delta_{\mu}},
 j_!\LL(\mathfrak{M}^{\ell\text{-ss}})_{\delta}
 :\mu\in I_{>0}
 \right\rangle,
\end{align}
ordered by decreasing $\mu$, where
\begin{align}\label{delta:center-review}
 \delta_{\mu}
 =\left.
 \bigl(\delta+\delta_P(\det\mathbb{L}_{p_{\mu}}^{1/2})\bigr)
 \right|_{\mathfrak{Z}_{\mu}}.
\end{align}
Here $j\colon\mathfrak{M}^{\ell\text{-ss}}\hookrightarrow\mathfrak{M}$ is the open immersion.
The functors $\Upsilon_{\mu}$ and $j_!$ are fully faithful, and $j_!$ is
left adjoint to the restriction functor 
$j^* \colon \LL(\mathfrak{M})_{\delta} \to \LL(\mathfrak{M}^{\ell\text{-ss}})_{\delta}$.
\end{thm}

The following consequence is important when checking that a coherent object
satisfying the weight condition is compact.  

\begin{cor}[{\cite[Corollary~4.25]{TodaLimits}}]
\label{cor:limit-compactness}
We have 
$\LL(\mathfrak{M})_{\delta}\subset
\widetilde{\LL}(\mathfrak{M})_{\delta}$, and an object
$\mathcal{E}\in\widetilde{\LL}(\mathfrak{M})_{\delta}$ belongs to
$\LL(\mathfrak{M})_{\delta}$ if and only if there exists $\mu$ such that
\begin{align}\label{orth:limit-compactness}
 \Hom(\mathcal{E},A)=0
\end{align}
for every $A\in\widetilde{\LL}(\mathfrak{M})_{\delta}$ supported on
$\mathfrak{S}_{>\mu}:=\mathfrak{M}\setminus
\mathfrak{M}_{\preceq\mu}$.
\end{cor}

For comparison, when a quasi-smooth QCA symmetric stack $\mathfrak{M}$ has
a separated diagonal and a good moduli space, 
it is proved in~\cite{BPT} 
(recalled in~\cite[Theorem~3.11]{TodaLimits})
that we have the 
semiorthogonal decomposition
\begin{align}\label{sod:coh-limit-review}
 \Coh(\mathfrak{M})
 =\left\langle
 \ast_{\lambda}\LL(\mathfrak{M}^{\lambda})_{\delta_{\lambda}}
 :\lambda\in
 |\operatorname{CL}(\mathfrak{M})_{\mathbb{Q}}|
 \right\rangle.
\end{align}
Here $\mathfrak{M}^{\lambda}\subset \mathrm{Grad}(\mathfrak{M})$ is the corresponding component, $\ast_{\lambda}$ is the Hall induction functor, and the
parameters $\delta_{\lambda}$ are determined by $\delta$, the relative
cotangent complex of $\ev_1$, and a choice of positive quadratic
norm.  
The semiorthogonal decomposition (\ref{sod:coh-limit-review}) is obtained 
from combinatorics of weight polytopes, while (\ref{sod:limit-general}) is obtained 
from the geometry of $\Theta$-stratifications.

\section{Moduli stacks of one-dimensional sheaves on surfaces}\label{sec:surfaces}
In this section, we apply Theorem~\ref{thm:limit-sod-general} to moduli 
stacks of pure one-dimensional sheaves on smooth quasi-projective surfaces under conditions ensuring that \cite[Assumption~4.1]{TodaLimits} holds. 
By comparing the resulting semiorthogonal decomposition with that in~\cite{BPT}, 
we propose a certain equivalence of categories, which generalizes the Dolbeault geometric Langlands conjecture for $\GL_r$-Higgs bundles in~\cite{PTlim}. 
\subsection{Some notation on surfaces}\label{subsec:notationS}
Let $S$ be a smooth quasi-projective surface over $k$. 
Here we fix notation. 
\subsubsection*{The groups of one-cycles}
Let $Z^1(S)$ be the abelian group of one-cycles on $S$, 
and $Z_1^{c}(S) \subset Z^1(S)$ the subgroup of compactly supported 
one-cycles. We have the intersection pairing 
\begin{align}\label{pair:Z}
  Z_1^c(S)\times   Z^1(S) \to \mathbb{Z}, \ 
    (C, D) \mapsto C \cdot D.
\end{align}
We define 
\begin{align*}
    N^1(S):=Z^1(S)/\equiv, \ N_1^c(S) :=Z_1^c(S)/\equiv
\end{align*}
where $\equiv$ is the numerical equivalence relation with 
respect to the pairing (\ref{pair:Z}). 
The pairing (\ref{pair:Z}) descends to the pairing 
\begin{align*}
  N_1^c(S) \times N^1(S)\to \mathbb{Z} 
\end{align*}
which induces the perfect pairing 
\begin{align*}
    N_1^c(S)_{\mathbb{Q}} \times  N^1(S)_{\mathbb{Q}}\to \mathbb{Q}. 
\end{align*}

The inclusion $Z_1^c(S) \subset Z^1(S)$ induces the map 
\begin{align}\label{map:N1S}
    N_1^c(S) \to N^1(S).
\end{align}
By composition, we also obtain the pairing 
\begin{align}\label{pair:N1}
    N_1^c(S)_{\mathbb{Q}} \times  N_1^c(S)_{\mathbb{Q}}
    \to   N_1^c(S)_{\mathbb{Q}} \times  N^1(S)_{\mathbb{Q}}
    \to \mathbb{Q}. 
\end{align}

\subsubsection*{The categories of one-dimensional sheaves}
We denote by 
\begin{align*}
    \Coh_{\leq 1}(S) \subset \Coh(S)
\end{align*}
the subcategory of compactly supported objects $E$ with $\dim \mathrm{Supp}(E) \leq 1$. 
The standard t-structure on $\Coh(S)$ restricts to $\Coh_{\leq 1}(S)$, with heart
\begin{align*}
    \Coh_{\leq 1}^{\heartsuit}(S) :=\Coh^{\heartsuit}(S) \cap \Coh_{\leq 1}(S).
\end{align*}
\subsubsection*{Maps from the K-groups}
For $E \in \Coh_{\leq 1}(S)$, we denote by $[E] \in N_1^c(S)$ the fundamental one-cycle 
associated with $E$. 
We have the following maps
\begin{align}\label{map:cl}
  &  \cl \colon K(\Coh_{\leq 1}(S)) \to \Lambda:=N_1^c(S) \oplus \mathbb{Z}, \ 
    E \mapsto ([E], \chi(E)), \\
 \notag  &  \cl \colon K(S) \to \Gamma:=\mathbb{Z}\oplus N^1(S), \ 
    E \mapsto (\rank(E), \det(E)).
\end{align}

For $E_1 \in \Coh(S)$ and $E_2 \in \Coh_{\leq 1}(S)$, their \textit{Euler pairing} is 
\begin{align*}
    \chi(E_1, E_2):=\sum_{i\in \mathbb{Z}} (-1)^i \dim \Hom^i(E_1, E_2).
\end{align*}
It descends to the pairing 
\begin{align*}
    \chi \colon \Gamma \times \Lambda \to \mathbb{Z}, \ 
    ((r, \beta_1), (\beta_2, \chi)) \mapsto r\chi-\beta_1 \beta_2.
\end{align*}
We have the commutative diagram 
\begin{align*}
    \xymatrix{
K(\Coh_{\leq 1}(S)) \ar[r] \ar[d]_-{\cl} & K(S) \ar[d]^-{\cl} \\
\Lambda \ar[r] & \Gamma
    }
\end{align*}
where the bottom arrow is $(\beta, m) \mapsto (0, \beta)$ via the map (\ref{map:N1S}). 
By composing $\chi$ with $\Lambda \to \Gamma$, we also have the pairing 
\begin{align*}
    \chi \colon \Lambda \times \Lambda \to \mathbb{Z}, \ ((\beta_1, \chi_1), (\beta_2, \chi_2)) \mapsto 
    -\beta_1 \beta_2.
\end{align*}

\subsubsection*{Degrees of sheaves on curves}
For a compactly supported effective divisor $C \subset S$ 
and a rank-one object $E \in \Coh(C)$, we define 
\begin{align}\label{degCE}
    \deg_C(E):=\chi(E)-\chi(\mathcal{O}_C).
\end{align}
By setting $\beta=[C] \in N_1^c(S)$, we can write it as
\begin{align}\label{deg:chi}
    \deg_C(E)=\chi(E)+\frac{1}{2}\beta (K_S+\beta).
\end{align}
Here the second $\beta$ is regarded as an element of $N^1(S)$ by the map (\ref{map:N1S}). 
Below we often use the transformation rule based on the formula (\ref{deg:chi}): 
\begin{align}\label{chi:trans}
    \chi \mapsto \chi_{\beta}:=\chi+\frac{1}{2}\beta (K_S+\beta)
\end{align}

\subsection{The derived stack of one-dimensional sheaves on surfaces}\label{subsec:onedim}
We consider the moduli stacks 
\begin{align*}
\mathfrak{M}_S \subset 
    \mathfrak{M}_S^{\leq 1}. 
\end{align*}
\begin{itemize}[leftmargin=1.5em]
\item The stack $\mathfrak{M}_S^{\leq 1}$ is the 
derived moduli stack of
$E \in \Coh^{\heartsuit}_{\leq 1}(S)$;
\item the substack 
$\mathfrak{M}_S \subset \mathfrak{M}_S^{\leq 1}$ is an open substack of pure one-dimensional sheaves $E\in \Coh_{\leq 1}^{\heartsuit}(S)$, i.e. there is no non-zero zero-dimensional subsheaf of $E$. 
\end{itemize}
We have the following proposition: 
\begin{prop}\label{prop:sym}
    The derived stack $\mathfrak{M}_S^{\leq 1}$ is symmetric. 
\end{prop}
\begin{proof}
    Take $E\in \Coh^{\heartsuit}_{\leq 1}(S)$, and let $x\in \mathfrak{M}_S^{\leq 1}(k)$ be the corresponding point. 
    Then we have 
    \begin{align*}
        \mathbb{T}_{\mathfrak{M}_S^{\leq 1}}|_{x}=\Hom^*(E, E)[1].
    \end{align*}
    By Serre duality, we have 
    \begin{align*}
          \mathbb{T}_{\mathfrak{M}_S^{\leq 1}}|_{x}^{\vee}\cong
          \Hom^{*}(E, E\otimes \omega_S)[1]
    \end{align*}
    as an object of $\Coh(BG_x)$. 

    If $\dim \mathrm{Supp}(E)=0$, then $E\otimes \omega_S \cong E$
    and $ \mathbb{T}_{\mathfrak{M}_S^{\leq 1}}|_{x}$ is a symmetric 
    virtual $G_x$-representation. Therefore we may assume that 
    $\dim \mathrm{Supp}(E)=1$. 
    Since $S$ is quasi-projective, we can write $\omega_S$ as 
    \begin{align*}
        \omega_S=\mathcal{O}_S(H_1-H_2)
    \end{align*}
    where $H_i \subset S$ are smooth divisors
    which intersect $\mathrm{Supp}(E)$ transversely, and 
    $E$ is pure one-dimensional in neighborhoods of $H_i$.
    Since we have 
    \begin{align*}
    [\omega_S]=[\mathcal{O}_S]-[\mathcal{O}_{H_2}]+[\omega_S \otimes \mathcal{O}_{H_1}] \in K(S),
    \end{align*}
    noting that $(E\otimes \omega_S)|_{H_1}\cong E|_{H_1}$ (as they are zero-dimensional sheaves), 
    we have the equality in $K(BG_x)$:
    \begin{align*}
        [\mathbb{T}_{\mathfrak{M}_S^{\leq 1}}|_{x}^{\vee}]
        =[\mathbb{T}_{\mathfrak{M}_S^{\leq 1}}|_{x}]
        +[\Hom_{H_2}^*(E|_{H_2}, E|_{H_2})]
        -[\Hom_{H_1}^*(E|_{H_1}, E|_{H_1})].
    \end{align*}

The action of $G_x=\mathrm{Aut}(E)$ on $\Hom^*(E|_{H_i}, E|_{H_i})$
factors through $\mathrm{Aut}(E|_{H_i})$. 
Therefore, by decomposing $E|_{H_i}$ according to its zero-dimensional supports, 
the proposition is reduced to the local statement in Lemma~\ref{lem:A1} below.
\end{proof}

\begin{lemma}\label{lem:A1}
 Let $Q\in \Coh^{\heartsuit}(\mathbb{A}^1)$ be a zero-dimensional 
 sheaf supported at $0\in \mathbb{A}^1$. Then 
 \begin{align*}
     [\Hom^*(Q, Q)]=0 \in K(B\mathrm{Aut}(Q)).
 \end{align*}
\end{lemma}
\begin{proof}
    Let $n$ be the length of $Q$. The moduli stack of sheaves of length $n$ on 
    $\mathbb{A}^1$, together with its good moduli space, is given by 
    \begin{align*}
        \pi\colon\mathcal{Y}=\mathfrak{gl}_n/\GL_n \to Y=\mathfrak{gl}_n \ssslash \GL_n. 
    \end{align*}
    The sheaf $Q$ corresponds to a point $y\in \pi^{-1}(0)$. 
    Let $y_0 \in \pi^{-1}(0)$ be the closed point, corresponding to 
    $\mathcal{O}_0^{\oplus n}$. Then we have 
    $G_y \subset G_{y_0}=\GL_n$ and 
    \begin{align*}
        [\mathbb{T}_{\mathcal{Y}}|_{y_0}]=[\mathfrak{gl}_n]-[\mathfrak{gl}_n]=0 \in K(B\GL_n).
    \end{align*}
    Let $i \colon BG_y \to BG_{y_0}$. 
    It follows that 
    \begin{align*}
        [\Hom^*(Q, Q)[1]]=[\mathbb{T}_{\mathcal{Y}}|_{y}]=i^*[\mathbb{T}_{\mathcal{Y}}|_{y_0}]=0 \in K(BG_y).
    \end{align*}
\end{proof}

In later sections, we will consider maps 
$\nu \colon \bgm \to \mathfrak{M}_S$. 
Such a map corresponds to a direct sum 
\begin{align}\label{dsum}
    F=F_1 \oplus \cdots \oplus F_k, \ F_i \in \Coh_{\leq 1}^{\heartsuit}(S)
\end{align}
such that $F_i$ is pure and has $\mathbb{G}_m$-weight $\nu_i$. 
Equivalently, it corresponds to the $\mathbb{Z}$-graded object 
\begin{align}\label{dsum2}
    \bigoplus_{m\in \mathbb{Z}}E_m, \ E_m \in \Coh^{\heartsuit}_{\leq 1}(S)
\end{align}
such that $E_m$ is pure and $E_m=0$ for $\lvert m \rvert \gg 0$. 
The objects (\ref{dsum}) and (\ref{dsum2}) are related by 
setting $E_m=\bigoplus_{i:\,\nu_i=m}F_i$ for $m\in\mathbb{Z}$. 

\subsection{Derived Hilbert--Chow map}
In this subsection, we construct a (derived) Hilbert--Chow map 
on $\mathfrak{M}_S$, which generalizes the Hitchin fibration for 
the moduli stacks of 
$\GL_r$-Higgs bundles. 
Namely we construct a map 
\begin{align}\label{map:h}h \colon \mathfrak{M}_S \to \rB
\end{align}
where $\rB$ is a \textit{derived scheme} classifying effective divisors $C \subset S$, with tangent complex $\Gamma(N_{C/S})$, 
where $N_{C/S}=\mathcal{O}_C(C)$ is the normal bundle of $C$ in $S$.
We construct such $\rB$ as a derived moduli stack of pairs $(L, s)$
where $L$ is a line bundle on $S$ with a non-zero section $s$. 

We first construct a map (\ref{map:h}) for a compactification of $S$. 
We take a smooth projective compactification of $S$
\begin{align*}
    S \subset \overline{S}. 
\end{align*}
Let $\Pic(\overline{S})$ be the derived moduli stack of line bundles on $\overline{S}$. There is a 
determinant morphism 
\begin{align}\label{map:det}
    \det \colon \mathfrak{M}_{\overline{S}} \to \mathrm{Pic}(\overline{S}), \ 
    F \mapsto \det F.
\end{align}
We have the universal line bundle 
\begin{align}\label{line:L}
    \mathcal{L} \to \overline{S} \times \mathrm{Pic}(\overline{S})
\end{align}
and consider the derived stack 
\begin{align*}
    \Spec \mathrm{Sym}(((p_{\mathrm{Pic}})_{*}\mathcal{L})^{\vee}) \to \mathrm{Pic}(\overline{S}).
\end{align*}
Here $p_{\mathrm{Pic}} \colon \overline{S} \times \mathrm{Pic}(\overline{S}) \to \mathrm{Pic}(\overline{S})$ is the projection. 
The above stack classifies pairs 
\begin{align*}(L, s), \ 
L \to \overline{S}, \ s \colon \mathcal{O}_{\overline{S}} \to L
\end{align*}
where $L$ is a line bundle and $s$ is a section of $L$.
We denote by 
\begin{align*}
  \overline{\rB} \subset   \Spec \mathrm{Sym}(((p_{\mathrm{Pic}})_{*}\mathcal{L})^{\vee})
\end{align*}
the open substack corresponding to $(L, s)$ such that $s\neq 0$. 
\begin{lemma}\label{lem:B}
The stack $\overline{\rB}$ is a quasi-smooth derived scheme, 
and each component of $\overline{\rB}^{\mathrm{cl}}$ is projective. Moreover 
there is a natural bijection between $\overline{\rB}(k)$
and the set of effective divisors on $\overline{S}$.
\end{lemma}
\begin{proof}
We first show the bijection between $\overline{\rB}(k)$ and 
effective divisors on $\overline{S}$. 
A pair $(L, s)$ corresponds to an effective divisor 
by taking the zero-locus of $s$; conversely, an effective divisor $C\subset \overline{S}$ determines 
$(\mathcal{O}_{\overline{S}}(C), s)$ where $s \colon \mathcal{O}_{\overline{S}} \hookrightarrow \mathcal{O}_{\overline{S}}(C)$ is the natural map. 

We next show that $\overline{\rB}$ is quasi-smooth. Let $C\subset \overline{S}$
correspond to $(L, s)$. Then the tangent complex of $\overline{\rB}$ is given by 
\begin{align*}
    \mathbb{T}_{\overline{\rB}}|_{(L, s)}=\Gamma((\mathcal{O}_{\overline{S}}\to L), L)=\Gamma(\mathcal{O}_C(C)).
\end{align*}
Since $\Gamma(\mathcal{O}_C(C))$
is concentrated in degrees $[0, 1]$, $\overline{\rB}$ is quasi-smooth. 

Finally we show that each component of $\overline{\rB}^{\mathrm{cl}}$ is projective. 
Let $\mathrm{Pic}_{\circ}(\overline{S}) \subset
\mathrm{Pic}(\overline{S})$ be a connected component, and set
\begin{align*}
    \overline{\rB}_{\circ}
    :=\overline{\rB}
    \times_{\mathrm{Pic}(\overline{S})}
    \mathrm{Pic}_{\circ}(\overline{S}).
\end{align*}
It is enough to show that 
$(\overline{\rB}_{\circ})^{\mathrm{cl}}$
is a projective scheme.

The stack $\mathrm{Pic}_{\circ}(\overline{S})$
is a $\mathbb{G}_m$-gerbe over a derived scheme
$J_{\circ}(\overline{S})$, whose classical truncation
is projective.
For any line bundle $L$ on $\overline{S}$, the action of
$\mathrm{Aut}(L)=\mathbb{G}_m$ on
$H^0(\overline{S},L)\setminus\{0\}$ is free.
It follows that the morphism
\begin{align*}
    (\overline{\rB}_{\circ})^{\mathrm{cl}}
    \to J_{\circ}(\overline{S})^{\mathrm{cl}}
\end{align*}
has fibers $\mathbb{P}(H^0(\overline{S},L))$,
and is a projective bundle over each stratum
on which $h^0(\overline{S},L)$ is constant.
There exists a sufficiently ample
effective divisor $H$ on $\overline{S}$ such that
$
    H^{>0}(\overline{S},L(H))=0
$
for every line bundle $L$ corresponding to a point
of $J_{\circ}(\overline{S})$.
Then 
$(\overline{\rB}_{\circ})^{\mathrm{cl}}$ is realized 
as a closed subscheme of the projective bundle
\begin{align*}
    P_{\circ}\to J_{\circ}(\overline{S})^{\mathrm{cl}}
\end{align*}
whose fiber over $[L]$ is
$\mathbb{P}(H^0(\overline{S},L(H)))$.
Since $J_{\circ}(\overline{S})^{\mathrm{cl}}$
is projective, so is $P_{\circ}$.
Therefore $(\overline{\rB}_{\circ})^{\mathrm{cl}}$
is a projective scheme.
\end{proof}

We also have the following lemma: 
\begin{lemma}\label{lem:lift2}
There is a lift of the map (\ref{map:det}) to $\overline{\rB}$
\begin{align*}
    \xymatrix{
    & \overline{\rB} \ar[d] \\
    \mathfrak{M}_{\overline{S}} \ar[r]^-{\det} \ar[ur]^-{h} & \mathrm{Pic}(\overline{S}).
    }
\end{align*}
\end{lemma}
\begin{proof}
Let $T$ be a derived affine scheme over $k$. 
We construct the morphism $h$ on $T$-valued points. 
Let
\begin{align*}
\mathcal{E}_T \in \Coh(\overline{S}\times T)
\end{align*}
be a $T$-family corresponding to a morphism
$T\to\mathfrak{M}_{\overline{S}}$.
Since $\overline{S}$ is smooth, $\mathcal{E}_T$ is perfect. Moreover, for any
geometric point $t\to T$, the restriction
$\mathcal{E}_t$ is a pure one-dimensional sheaf on the smooth surface
$\overline{S}$. Hence $\mathcal{E}_t$ is Cohen--Macaulay of codimension one,
and therefore has Tor-amplitude contained in $[-1,0]$.
By the fiberwise criterion for the Tor-amplitude of perfect complexes,
$\mathcal{E}_T$ has Tor-amplitude contained in $[-1,0]$ on
$\overline{S}\times T$.

It follows that, Zariski locally on $\overline{S}\times T$, the perfect
complex $\mathcal{E}_T$ is represented by a two-term complex
\begin{align}\label{resol:derived}
V_1 \stackrel{d}{\to} V_0
\end{align}
of vector bundles; see, for example,
\cite[Proposition~2.13(6),(7)]{AG14}.
Since $\mathcal{E}_T$ has virtual rank zero, $V_0$ and $V_1$ have the same
rank. We therefore obtain a section
\begin{align}\label{det:section}
\det(d) \colon
\mathcal{O}_{\overline{S}\times T}
\to
\det V_0\otimes \det V_1^{\vee}
\simeq \det(\mathcal{E}_T).
\end{align}

The local sections (\ref{det:section}) are independent of the choice of the
presentation (\ref{resol:derived}) and glue to a canonical section
\begin{align}\label{section:ET}
s_{\mathcal{E}_T} \colon
\mathcal{O}_{\overline{S}\times T}
\to
\det(\mathcal{E}_T).
\end{align}
Indeed, this is the determinant section associated with a rank-zero perfect
complex of Tor-amplitude $[-1,0]$. Locally it is given by the determinant of
the differential as in (\ref{det:section}); its invariance under a change of
a locally free presentation follows from the determinant construction for
perfect complexes, or equivalently from the $\operatorname{Div}$
construction of Knudsen--Mumford~\cite{KM76}. The construction is compatible
with pull-back: for any morphism of derived affine schemes $T'\to T$, we have
\begin{align*}
\det(\mathcal{E}_{T'})
\simeq
\det(\mathcal{E}_T)|_{\overline{S}\times T'}, \quad
s_{\mathcal{E}_{T'}}
=
s_{\mathcal{E}_T}|_{\overline{S}\times T'}.
\end{align*}

For any geometric point $t\to T$, the section
$s_{\mathcal{E}_T}|_{\overline{S}\times t}$ is non-zero. Indeed, away from
$\mathrm{Supp}(\mathcal{E}_t)$, the morphism $d$ in
(\ref{resol:derived}) is an isomorphism, and hence $\det(d)$ is invertible.
Therefore the pair
$
\left(
\det(\mathcal{E}_T),s_{\mathcal{E}_T}
\right)
$
defines a $T$-point of $\overline{\rB}$. Since this construction is
compatible with base change in $T$, it defines a morphism of derived stacks
\begin{align*}
h\colon
\mathfrak{M}_{\overline{S}}
\to
\overline{\rB}.
\end{align*}
By construction, its composition with
$\overline{\rB}\to\Pic(\overline{S})$ is the determinant morphism
(\ref{map:det}).

Finally, on classical points, if
\begin{align*}
0\to V_1\stackrel{d}{\to}V_0\to E\to 0
\end{align*}
is a locally free resolution of a pure one-dimensional sheaf $E$, the
zero divisor of $\det(d)$ is the zeroth Fitting divisor of $E$. Hence the
morphism $h$ sends $E$ to its Fitting support.
\end{proof}

Let $\cL_{\overline{\rB}}$ be the pull-back of $\cL$ in (\ref{line:L})
to $\overline{S}\times \overline{\rB}$. 
There is a universal morphism 
\begin{align*}
    \overline{s} \colon \mathcal{O}_{\overline{S}\times \overline{\rB}}\to \cL_{\overline{\rB}}.
\end{align*}
We define 
\begin{align*}
    \overline{\mathrm{C}} \hookrightarrow \overline{S}\times \overline{\rB}
\end{align*}
to be the derived zero locus of $\overline{s}$. It is the universal 
effective divisor on $\overline{S}$ over $\overline{\rB}$. 

We then restrict to one-dimensional sheaves on $S$. 
Let 
\begin{align*}
    \rB \subset \overline{\rB}
\end{align*}
be the open subscheme corresponding to compactly supported effective divisors in $S$. 
By construction, 
the map $h$ restricted to the open substack 
$\mathfrak{M}_S \subset \mathfrak{M}_{\overline{S}}$
fits into the Cartesian square 
\begin{align*}
    \xymatrix{
\mathfrak{M}_S \inclusion \ar[d]_-{h} \diasquare & \mathfrak{M}_{\overline{S}} \ar[d]^-{h} \\
\rB \inclusion & \overline{\rB}.
    }
\end{align*}

\begin{lemma}
The above construction of $h \colon \mathfrak{M}_S \to \rB$ is independent of 
 a choice of $S\subset \overline{S}$. 
\end{lemma}
\begin{proof}
    Let $S \subset \overline{S}$ and $S \subset \overline{S}'$
    be two smooth compactifications, and let 
    \begin{align*}
        h \colon \mathfrak{M}_S \to \rB, \ 
         h' \colon \mathfrak{M}_S \to \rB'
    \end{align*}
    be the derived Hilbert--Chow maps constructed from 
    the above two compactifications. We may assume that 
    there is a birational morphism $f \colon \overline{S}' \to \overline{S}$ which is an isomorphism over $S$. 
    Let 
    \begin{align*}
         h \colon \mathfrak{M}_{\overline{S}} \to \overline{\rB}, \ 
         h' \colon \mathfrak{M}_{\overline{S}'} \to \overline{\rB}'
    \end{align*}
    be the derived Hilbert--Chow maps for $\overline{S}$, 
    $\overline{S}'$. Then we have the morphism 
    $f^* \colon \overline{\rB} \to \overline{\rB}'$
    given by $(L, s) \mapsto (f^*L, f^*s)$. 
    It restricts to an equivalence 
    $\rB \stackrel{\sim}{\to} \rB'$, hence the lemma follows. 
\end{proof}

We also define the universal divisor 
\begin{align*}
\mathrm{C}:=\overline{\mathrm{C}}\times_{\overline{\rB}}\rB \to \rB.
\end{align*}
Similarly it is independent of $S\subset \overline{S}$. Moreover the universal sheaf $\cE$ is written as  
\begin{align}\label{ECreg}
    \cE \cong i_{*}\cE_{\mathrm{C}} \in \Coh(S \times \mathfrak{M}_S), \ 
    \cE_{\mathrm{C}} \in \Coh(\mathrm{C}\times_{\rB}\mathfrak{M}_S)
\end{align}
where $i \colon \mathrm{C}\times_{\rB}\mathfrak{M}_S \hookrightarrow S \times \mathfrak{M}_S$ is the closed immersion.

\subsection{Moduli stacks of Bridgeland semistable objects}
We construct a Bridgeland stability condition on
$\Coh_{\leq 1}(S)$ with the standard heart. 
Take a rational element 
\begin{align*}
    B+i H\in N^1(S)_{\mathbb{C}}
\end{align*}
such that $H$ is ample. 
We consider the group homomorphism 
\begin{align*}
    Z_{B, H} \colon K(\Coh_{\leq 1}(S)) \to \mathbb{C}
\end{align*}
given by 
\begin{align*}
    Z_{B, H}(F)=-\int_S e^{-(B+iH)}\ch(F).
\end{align*}
It factors through the map $\cl$ in (\ref{map:cl})
\begin{align*}
Z_{B, H} \colon K(\Coh_{\leq 1}(S)) \stackrel{\cl}{\to} \Lambda \to \mathbb{C}.
\end{align*}
We use the same symbol $Z_{B, H}$ for the 
factored map 
\begin{align}\label{ZBH}
    Z_{B, H} \colon \Lambda \to \mathbb{C}, \ 
    (\beta, \chi) \mapsto -\chi+\left(B-\frac{1}{2}K_S\right)\beta+iH \cdot \beta.
\end{align}

The \textit{slope} of the central charge $Z_{B, H}$ is given by 
\begin{align}\label{slope:mu}
\mu^{\sigma}(\beta, \chi) :=-\frac{\mathrm{Re}Z(\beta, \chi)}{\mathrm{Im}Z(\beta, \chi)}=\frac{\chi +(K_S/2-B)\beta}{H\cdot \beta}\in \mathbb{Q} \cup \{\infty\}.
\end{align}
Here $\mu^{\sigma}(\beta, \chi)=\infty$ if $H \cdot \beta=0$. 
\begin{remark}
From (\ref{ZBH}), 
there is $\alpha \in K(S)_{\mathbb{C}}$
such that 
\begin{align*}
    Z_{B, H}(F)=\chi(\alpha \otimes F).
\end{align*}
Explicitly, we have 
\begin{align}\label{classalpha}
    \alpha=-[\mathcal{O}_S]+[G_1]+[G_2]i
\end{align}
where $G_i$ are $\mathbb{Q}$-linear combinations of one-dimensional sheaves with fundamental one-cycles $B-K_S/2$, $H$, respectively. 
\end{remark}

The pair 
\begin{align}\label{sigma:stab}
    \sigma_{B, H}=(Z_{B, H}, \Coh_{\leq 1}^{\heartsuit}(S))
\end{align}
determines a Bridgeland stability condition~\cite{Brs1} on $\Coh_{\leq 1}(S)$. 
For $\sigma=\sigma_{B, H}$, we denote by 
\begin{align*}
    \mathfrak{M}_S^{\sigma\text{-ss}} \subset \mathfrak{M}_S
\end{align*}
the derived open substack of $\sigma$-semistable objects. 
We have the decompositions into open and closed substacks 
\begin{align*}
    \mathfrak{M}_S^{\sigma\text{-ss}}=\bigsqcup_{(\beta, \chi)\in \Lambda}\mathfrak{M}_S(\beta, \chi)^{\sigma\text{-ss}}, \ \mathfrak{M}_S=\bigsqcup_{(\beta, \chi)\in \Lambda}\mathfrak{M}_S(\beta, \chi)
\end{align*}
where each component corresponds to one-dimensional sheaves $E$ with $[E] =\beta$
in $N_1^c(S)$ and $\chi(E)=\chi$. 

Since any non-zero coherent sheaf $E$ has scalar automorphisms 
$\mathbb{G}_m \subset \mathrm{Aut}(E)$, we have the 
decomposition 
\begin{align*}
    \Coh(\mathfrak{M}_S(\beta, \chi)^{\sigma\text{-ss}})=\bigoplus_{w\in \mathbb{Z}}\Coh(\mathfrak{M}_S(\beta, \chi)^{\sigma\text{-ss}})_w
\end{align*}
where each summand corresponds to $\mathbb{G}_m$-weight $w$, 
called \textit{central weight} $w$ part. 

\subsection{Line bundles on the moduli stack}
Recall the universal one-dimensional sheaf (\ref{ECreg})
\begin{align}\label{univ:E}
    \cE \in \Coh(S\times \mathfrak{M}_S). 
\end{align}
For an element $\gamma \in K(S)$,
we define 
\begin{align*}
    \cL_{\gamma} := \det p_{\mathfrak{M}*}(p_S^* \gamma \otimes \cE) \in \mathrm{Pic}(\mathfrak{M}_S).
\end{align*}
Here $p_{\mathfrak{M}} \colon S\times \mathfrak{M}_S \to \mathfrak{M}_S$ and $p_S \colon S\times \mathfrak{M}_S \to S$ are the projections. 
Its image under the map 
\begin{align*}
    \delta_P \colon \mathrm{Pic}(\mathfrak{M}_S) \to \mathrm{Pic}(\mathfrak{M}_S)^{\mathrm{CL}}
\end{align*}
is described as follows. Let $\nu \colon \bgm \to \mathfrak{M}_S$ be a map corresponding to 
the direct sum 
\begin{align*}
    \mathrm{pt} \mapsto M=M_1\oplus \cdots \oplus M_k
\end{align*}
with $\mathbb{G}_m$-weights $(\nu_1, \ldots, \nu_k)$. Then we have 
\begin{align*}
\delta_P(\cL_{\gamma})(\nu)=
    c_1(\nu^* \cL_{\gamma})=\sum_{i=1}^k \chi(\gamma \otimes M_i)\cdot \nu_i.
\end{align*}
Since it only depends on the numerical class of $\gamma$, it descends to
the additive map 
\begin{align*}
\delta_N \colon 
\Gamma=\mathbb{Z} \oplus N^1(S) \to \mathrm{Pic}(\mathfrak{M}_S)^{\mathrm{CL}}, \ (r, D)\mapsto \delta_P(\cL_{\gamma})
\end{align*}
where $\cl(\gamma)=(r, D)$. 
It extends to the map 
\begin{align}\label{deltaN}
    \delta_N \colon \Gamma_{\mathbb{Q}} \to \mathrm{Pic}(\mathfrak{M}_S)_{\mathbb{Q}}^{\mathrm{CL}}.
\end{align}
Note that $\delta_N(r, D)$ has central weight $D\beta +r \chi$ on $\mathfrak{M}_S(\beta, \chi)$.
In what follows, we only consider the case of $r=0$, and set 
\begin{align*}
    \delta_N(D):=\delta_N(0, D) \in \mathrm{Pic}(\mathfrak{M}_S)_{\mathbb{Q}}^{\mathrm{CL}}.
\end{align*}
\begin{remark}
    Although the construction of $\delta_N$ allows a more general parameter in $\Gamma_{\mathbb{Q}}$, in this paper we restrict to the rank-zero part $N^1(S)_{\mathbb{Q}} \subset \Gamma_{\mathbb{Q}}$, which is sufficient for the purpose of this paper. However, allowing additional parameters in $\Gamma_{\mathbb{Q}}$ may be relevant to the study of a 
    version of quantum geometric Langlands for Higgs bundles and 
    a generalization of the Hecke correspondence in Section~\ref{sec:hecke}.
\end{remark}

For a non-zero effective class $\beta \in N_1^c(S)$ and $w\in \mathbb{Q}$, we consider the affine subspace
\begin{align*}
    N^1(S)_{\mathbb{Q}}^{\beta\cdot (-)=w} \subset N^1(S)_{\mathbb{Q}}
\end{align*}
consisting of elements $D$ such that $D\cdot \beta=w$. 
For a choice of an ample class $H\in N^1(S)_{\mathbb{Q}}$, 
since $H\cdot \beta>0$, the above embedding admits a projection 
\begin{align}\label{proj:A}
    p_{H, \beta, w} \colon N^1(S)_{\mathbb{Q}} \to   N^1(S)_{\mathbb{Q}}^{\beta\cdot (-)=w}, \ 
    D\mapsto D+\frac{w-D\cdot \beta}{H\cdot \beta}H.
\end{align}

For a stability condition $\sigma=\sigma_{B, H}$ as in (\ref{sigma:stab}), we set 
\begin{align}\label{def:Dsigma}
    D_{\sigma, \beta, w}:=p_{H, \beta, w}\left(\frac{\beta}{2}+B  \right) \in N^1(S)_{\mathbb{Q}}^{\beta \cdot (-)=w}.
\end{align}
It is explicitly written as 
     \begin{align*}
        D_{\sigma, \beta, w}=\frac{w}{H\cdot \beta}H-\frac{\beta^2}{2H\cdot \beta}H-\frac{B\cdot \beta}{H\cdot \beta}H+\frac{\beta}{2}+B.
        \end{align*}
We note that, for the transformation $w\mapsto w_{\beta}$ as in (\ref{chi:trans}), we have 
\begin{align}\label{formula:D}
    D_{\sigma, \beta, w_{\beta}}=\mu^{\sigma}(\beta, w)H+\frac{\beta}{2}+B.
\end{align}
Here $\mu^{\sigma}(\beta, w) \in \mathbb{Q}$ is the slope of the central charge (\ref{slope:mu}). 

Later we will use the element 
\begin{align}\label{twist:D}
    \delta_N( D_{\sigma, \beta, w_{\beta}}) \in \mathrm{Pic}(\mathfrak{M}_S)_{\mathbb{Q}}^{\mathrm{CL}}
\end{align}
for the relevant twist of the limit category depending on a choice of a stability condition $\sigma$. 
Here we regard an element of $N^1(S)_{\mathbb{Q}}$ as an element of $\Gamma_{\mathbb{Q}}$
by $D \mapsto (0, D)$. 
        
\subsection{Line bundles via the Deligne pairing}
We denote by 
\begin{align*}
    \mathfrak{M}_S^{\mathrm{reg}} \subset \mathfrak{M}_S
\end{align*}
the open substack corresponding to line bundles on an effective divisor on $S$, called the \textit{regular locus}. 
Then 
\begin{align*}
\cE_{\mathrm{C}}^{\mathrm{reg}}:=\cE_{\mathrm{C}}|_{\mathrm{C}\times_{\rB}\mathfrak{M}_S^{\mathrm{reg}}} \in \Coh(\mathrm{C}\times_{\rB}\mathfrak{M}_S^{\mathrm{reg}})
\end{align*}
is a line bundle, where $\cE_{\mathrm{C}}$ is given in (\ref{ECreg}). 
Below we construct a line bundle on $\mathfrak{M}_S^{\mathrm{reg}}\times_{\rB} \mathfrak{M}_S$
via \textit{Deligne pairing}. 

Consider the following diagram 
\begin{align*}
\xymatrix{
& \ar[ld]_-{p_{12}} \mathfrak{M}_S^{\mathrm{reg}}\times_{\rB}\rC \times_{\rB} \mathfrak{M}_S \ar[rd]^-{p_{23}} \ar[d]_-{p_{13}} & \\
\mathfrak{M}_S^{\mathrm{reg}}\times_{\rB}\rC & 
\mathfrak{M}_S^{\mathrm{reg}}\times_{\rB} \mathfrak{M}_S & \rC \times_{\rB} \mathfrak{M}_S
}
\end{align*}
and define the following line bundle on $\mathfrak{M}_S^{\mathrm{reg}}\times_{\rB} \mathfrak{M}_S$
\begin{align*}
    \cP^{\mathrm{reg}}:=
    &\det p_{13*}(p_{12}^* \cE^{\mathrm{reg}}_{\rC}\otimes p_{23}^* \cE_{\rC}) \\
    &\otimes 
    (\det p_{13*}p_{12}^* \cE_{\rC}^{\mathrm{reg}})^{-1} \otimes 
     (\det p_{13*}p_{23}^* \cE_{\rC})^{-1}
     \otimes \det p_{13*}\mathcal{O}.
\end{align*}
Note that the above determinants make sense since 
the objects before taking the determinants are perfect complexes. 

For $b\in \rB$, we also consider the regular locus in the fiber 
\begin{align*}
h^{-1}(b)^{\mathrm{reg}}:=h^{-1}(b)\cap \mathfrak{M}_S^{\mathrm{reg}}.
\end{align*}
Let $C\subset S$ be a compactly supported effective divisor with fundamental 
class $\beta \in N_1^c(S)$ and $b\in \rB$ the corresponding point. 
The line bundle $\cP^{\mathrm{reg}}$ over $b$ gives the 
line bundle 
\begin{align*}
    \cP_b^{\mathrm{reg}} \to h^{-1}(b)^{\mathrm{reg}} \times h^{-1}(b)
\end{align*}
whose fiber at $(F_1, F_2)$ is given by the \textit{Deligne pairing} 
\begin{align*}
    \cP_{(F_1, F_2)}=\det \Gamma(F_1 \otimes F_2)\otimes \det \Gamma(F_1)^{-1}\otimes \det \Gamma(F_2)^{-1}\otimes \det \Gamma(\cO_C).
\end{align*}
For a line bundle $L\to C$, we have the associated line bundle on $h^{-1}(b)$
\begin{align*}
    \cP_L:=(\iota_L \times \id)^*\cP_b^{\mathrm{reg}} \to h^{-1}(b),
\end{align*}
where $\iota_L \colon \mathrm{pt} \to h^{-1}(b)^{\mathrm{reg}}$
is the map corresponding to $L$. 

Below we regard $E\in \Coh^{\heartsuit}(C)$ as an object of 
$\Coh_{\leq 1}^{\heartsuit}(S)$ by the push-forward along 
$C \hookrightarrow S$. 
We have the following proposition, which relates the semistability of $L$ 
to the weight condition for the line bundle $\cP_L$ to lie in the limit category. 
This is a generalization of~\cite[Proposition~2.24]{TodaA}, and motivates the Dolbeault geometric Langlands-type 
conjecture in the next subsection: 
\begin{prop}\label{prop:line}
Let $\sigma=\sigma_{B, H}$ for a rational element $B+iH \in N^1(S)_{\mathbb{C}}$
such that $H$ is ample. 
Then 
a line bundle $L$ on $C$ with $\chi(L)=w$ is $\sigma$-semistable if and only if, for any map 
$\nu \colon \bgm \to h^{-1}(b)$, we have 
\begin{align}\label{wt:L}
    \wt(\nu^* \cP_L) \subset \frac{1}{2}\left[c_1(\nu^* \mathbb{L}_{h^{-1}(b)}^{<0}), c_1(\nu^* \mathbb{L}_{h^{-1}(b)}^{>0})  \right]+\delta_N(D_{\sigma, \beta, w_{\beta}})(\nu).
\end{align} 
Here we have set $\beta=[C]$, $w\mapsto w_{\beta}$ is given by (\ref{chi:trans}) so that 
$\deg_C(L)=w_{\beta}$, and $\delta_N(D_{\sigma, \beta, w_{\beta}})$ is given in (\ref{twist:D}).
\end{prop}
\begin{proof}
Let $\nu$ be a map 
   $\nu \colon \bgm \to h^{-1}(b)$ corresponding to the direct sum 
   \begin{align*}
      \mathrm{pt} \mapsto M=M_1 \oplus \cdots \oplus M_k, \ M_i \in \Coh^{\heartsuit}(C)
   \end{align*}
   with $\mathbb{G}_m$-weights 
   \begin{align*}(\nu_1, \ldots, \nu_k), \ \nu_1>\cdots>\nu_k.
   \end{align*}
   Let $C_i \subset S$ be the Fitting support of $M_i$, and set 
   \begin{align*}l_i=\deg_{C_i}(L|_{C_i}), \ \beta_i=[C_i] \in N_1^c(S).
   \end{align*}
   By the definition of $\cP_L$, we have 
\begin{align*}
    \nu^* \cP_L=\bigotimes_{i=1}^k\left( \det \Gamma(L\otimes M_i)\otimes \det \Gamma(M_i)^{-1}\right)
    \otimes \det \Gamma(L)^{-1}\otimes \det \Gamma(\cO_C).
\end{align*}
  Therefore we have the identity  
   \begin{align}\label{id:wt}
       c_1(\nu^* \cP_L)-\frac{w_{\beta}}{H\cdot \beta}c_1(\nu^* \cL_H)&=\sum_{i=1}^k l_i \cdot \nu_i-\frac{w_{\beta}}{H\cdot \beta}\sum_{i=1}^k (H\cdot \beta_i)\nu_i.
       \end{align}

       On the other hand, since we have 
       \begin{align*}
           \deg_C(L)=l_1+\cdots+l_k=w_{\beta}, 
       \end{align*} 
  we have the identity
   \begin{align*}
       \widetilde{l}_1+\cdots+\widetilde{l}_k=0, \ 
       \widetilde{l}_i:=l_i-w_{\beta}\frac{H \cdot \beta_i}{H\cdot \beta}.
   \end{align*}
    For a subset $I\subset \{1, \ldots, k\}$, we set  
   \begin{align*}
   l_I:=\sum_{i\in I}l_i,\quad \widetilde{l}_I:=\sum_{i\in I} \widetilde{l}_i, \quad \beta_I:=\sum_{i\in I} \beta_i.
   \end{align*}
   Setting $I_i:=\{1, \ldots, i\}$, we can rewrite the right-hand side of
   (\ref{id:wt}) as
   \begin{align*}\sum_{i=1}^k \widetilde{l}_i \cdot \nu_i =\sum_{i=1}^{k-1}\widetilde{l}_{I_i}(\nu_i-\nu_{i+1}).
   \end{align*}
   
   We also have 
   \begin{align}\label{c1:nu}
       c_1(\nu^* \mathbb{L}_{h^{-1}(b)}^{<0})=\sum_{i>j} \beta_i \beta_j (\nu_i-\nu_j) =
       -\sum_{i=1}^{k-1}\beta_{I_i} \beta_{I_i^{\circ}}(\nu_i-\nu_{i+1}).
   \end{align}
Here $I^{\circ} \subset \{1, \ldots, k\}$ is the complement of $I$. 
We also have the identity 
\begin{align}\label{id:Dsigma}
    D_{\sigma, \beta, w_{\beta}}-\frac{w_{\beta}}{H\cdot \beta}H=\frac{1}{2}\beta+B-\frac{\beta^2}{2H\cdot \beta}H-\frac{B\cdot \beta}{H\cdot \beta}H.
\end{align}
Then we have 
\begin{align}\label{c1:beta}
    &c_1(\nu^* \cL_{D_{\sigma, \beta, w_{\beta}}-\frac{w_{\beta}}{H\cdot\beta}H}) \\
 \notag   &=
    \sum_{i=1}^{k-1} \left(\frac{1}{2}\beta+B-\frac{\beta^2}{2H\cdot \beta}H-\frac{B\cdot \beta}{H\cdot \beta}H  \right)\beta_{I_i} \cdot (\nu_i-\nu_{i+1}).
\end{align}

   Suppose that $L$ is $\sigma$-semistable. For any non-empty proper subset $I\subset \{1, \ldots, k\}$
   and its complement $I^{\circ}$, we have the exact sequence 
   \begin{align*}
   0\to L(-C_I)|_{C_{I^{\circ}}} \to L \to L|_{C_I} \to 0.
\end{align*}
Here $C_I=\sum_{i\in I}C_i$, so that $\beta_I=[C_I]\in N_1^c(S)$. 
By the $\sigma$-semistability, we have the inequality 
\begin{align*}
    \frac{\chi(L)+(K_S/2-B)\cdot \beta}{H\cdot \beta} \leq \frac{\chi(L|_{C_I})+(K_S/2-B)\cdot \beta_I}{H \cdot \beta_I}. 
\end{align*}
   We have the identity 
   \begin{align*}
       \chi(L|_{C_I})=l_I+\chi(\cO_{C_I}), \ \chi(\cO_{C_I})=-\frac{1}{2}\beta_I(K_S+\beta_I).
   \end{align*}
   Therefore we obtain the inequality 
   \begin{align*}
         \frac{w_{\beta}-B\cdot \beta-\beta^2/2}{H\cdot \beta} \leq \frac{l_I-B\cdot \beta_I-\beta_I^2/2}{H \cdot \beta_I}.
   \end{align*}
   It follows that 
   \begin{align*}
       \widetilde{l}_I &\geq \frac{H\cdot \beta_I}{H\cdot \beta}\left(w_{\beta}-\beta\left(B+\frac{1}{2}\beta\right)  \right)
       +\beta_I\left(B+\frac{1}{2}\beta_I\right)-w_{\beta}\frac{H\cdot \beta_I}{H\cdot \beta} \\
       &=-\frac{H\cdot \beta_I}{H\cdot \beta}\beta\left(B+\frac{1}{2}\beta \right)  +\beta_I\left(B+\frac{1}{2}\beta_I\right) \\
       &=-\frac{1}{2}\beta_{I} \beta_{I^{\circ}}+\left(\frac{1}{2}\beta+B-\frac{\beta^2}{2H\cdot \beta}H-\frac{B\cdot \beta}{H\cdot \beta}H  \right)\beta_{I}.
   \end{align*}
   By (\ref{c1:beta}) and (\ref{c1:nu}), it follows that the lower bound in (\ref{wt:L}) holds. 
   The upper bound is proved similarly. 

Conversely, suppose that $L$ is not $\sigma$-semistable. 
Then there is an exact sequence in $\Coh^{\heartsuit}(C)$
\begin{align*}
    0\to E_2 \to L \to E_1 \to 0
\end{align*}
such that, by setting $\beta_i=[E_i] \in N_1^c(S)$, we have  
\begin{align}\label{ineq:chiLE}
    \frac{\chi(L)+(K_S/2-B)\cdot \beta}{H\cdot \beta}>
    \frac{\chi(E_1)+(K_S/2-B)\cdot \beta_1}{H\cdot \beta_1}.
\end{align}

By taking the saturation of $E_2$, we may assume that $E_1$ is pure one-dimensional. Since we have the surjection \begin{align*}\mathcal{O}_C \twoheadrightarrow E_1 \otimes L^{-1}
\end{align*}
we have 
$E_1 \cong L\otimes \mathcal{O}_{C_1}$ for a decomposition 
$C=C_1+C_2$ into effective divisors such that $\beta_i=[C_i]$.  
Setting $l_1=\deg_{C_1}(L|_{C_1})$, we deduce from (\ref{ineq:chiLE}) that  
\begin{align*}
    \widetilde{l}_1=l_1-w_{\beta}\frac{H\cdot \beta_1}{H\cdot \beta}
   < -\frac{1}{2}\beta_{1} \beta_{2}+\left(\frac{1}{2}\beta+B-\frac{\beta^2}{2H\cdot \beta}H-\frac{B\cdot \beta}{H\cdot \beta}H  \right)\beta_{1}.
\end{align*}
Then for any map $\nu \colon \bgm \to h^{-1}(b)$
corresponding to the direct sum $\mathcal{O}_{C_1} \oplus \mathcal{O}_{C_2}$ with $\mathbb{G}_m$-weights $(\nu_1, \nu_2)$
with $\nu_1>\nu_2$, we have 
\begin{align*}
    &c_1(\nu^*\cP_L)-\frac{w_{\beta}}{H\cdot \beta}c_1(\nu^* \cL_H)\\
    &=\widetilde{l}_1(\nu_1-\nu_2) 
    <\frac{1}{2}c_1(\nu^* \mathbb{L}_{h^{-1}(b)}^{<0})
    +c_1(\nu^* \cL_{D_{\sigma, \beta, w_{\beta}}-\frac{w_{\beta}}{H\cdot \beta}H}).
\end{align*}
Therefore the condition (\ref{wt:L}) is not satisfied. 
\end{proof}

\subsection{The Dolbeault geometric Langlands-type conjecture}
Let $S$ be a smooth quasi-projective surface over $k$. 
In this subsection, we propose a conjecture 
relating the dg-categories of coherent sheaves on moduli stacks of 
semistable objects to the limit categories of the moduli 
stacks of pure one-dimensional sheaves without a stability condition. 

Recall the divisor class $D_{\sigma, \beta, w}$ as in (\ref{def:Dsigma}). 
We give the following definition: 
\begin{defn}
    For a non-zero effective class $\beta\in N_1^c(S)$, $\chi, w\in \mathbb{Z}$,
and $\sigma=\sigma_{B, H}$, we define 
\begin{align*}
&\IndL(\mathfrak{M}_S(\beta, \chi))_{w}^{\sigma}:=\IndL(\mathfrak{M}_S(\beta, \chi))_{\delta_N(D_{\sigma, \beta, w})}. 
\end{align*}
The subcategory of compact objects is denoted by 
\begin{align*}
 \LL(\mathfrak{M}_S(\beta, \chi))_{w}^{\sigma}:=\LL(\mathfrak{M}_S(\beta, \chi))_{\delta_N(D_{\sigma, \beta, w})}.
\end{align*}
\end{defn}

We propose the following conjecture, which generalizes the \textit{Dolbeault geometric Langlands conjecture} for $\GL_r$-Higgs bundles in~\cite{PTlim}:
\begin{conj}\label{conj:CY3}
Suppose that $-K_S$ is nef. Then 
    there is a $\rB$-linear equivalence 
    \begin{align}\label{conj:equiv}
    \IndCoh(\mathfrak{M}_S(\beta, w)^{\sigma\text{-ss}})_{-\chi_{\beta}}
    \simeq \IndL(\mathfrak{M}_S(\beta, \chi))_{w_{\beta}}^{\sigma}.
    \end{align}
     Here $(\chi, w)\mapsto (\chi_{\beta}, w_{\beta})$ is the transformation rule (\ref{chi:trans}). 
    \end{conj}

\begin{remark}\label{rmk:equiv:compact}
The left-hand side of (\ref{conj:equiv}) is compactly generated by~\cite{MR3037900}, and the right-hand side is 
also compactly generated by Theorem~\ref{thm:sodLS} below.
    Therefore an equivalence (\ref{conj:equiv}) is equivalent to 
    an equivalence 
    \begin{align*}
        \Coh(\mathfrak{M}_S(\beta, w)^{\sigma\text{-ss}})_{-\chi_{\beta}}
    \simeq \LL(\mathfrak{M}_S(\beta, \chi))_{w_{\beta}}^{\sigma}.
    \end{align*}
\end{remark}

   \begin{remark}
       The assumption that $-K_S$ is nef is relevant in checking Assumption~\ref{assum:dag} (ii). A conjectural equivalence without this assumption will be discussed in Section~\ref{sec:DLCY3}.
   \end{remark}

\subsection{Quasi-BPS categories}\label{subsec:qbps}
Here we define the \textit{quasi-BPS category}, 
which generalizes the construction for K3 surfaces~\cite{PTK3}
and for Higgs bundles~\cite{PThiggs}:
\begin{defn}\label{def:qbpspm}
For a non-zero effective class $\beta\in N_1^c(S)$, $\chi, w\in \mathbb{Z}$, 
and $\sigma=\sigma_{B, H}$, 
the positive and negative \textit{quasi-BPS categories} are defined by 
\begin{align*}
&\mathbb{T}^{+}_S(\beta, \chi)_{w}^{\sigma} :=\LL(\mathfrak{M}_S(\beta, \chi)^{\sigma\text{-ss}})_{\delta_N(D_{\sigma, \beta, w})}, \\
&\mathbb{T}^{-}_S(\beta, \chi)_{w}^{\sigma} :=\LL(\mathfrak{M}_S(\beta, \chi)^{\sigma\text{-ss}})_{\delta_N(-D_{\sigma, \beta, -w})}.
\end{align*}
\end{defn}
Note that 
both parameters have central weight $w$. We use the notation
$\mathbb{T}_S(\beta,\chi)_w^\sigma:=\mathbb{T}^{+}_S(\beta,\chi)_w^\sigma$
when there is no sign. 

The quasi-BPS category was introduced in~\cite{PTK3} for K3 surfaces, 
and in~\cite{PThiggs} for Higgs bundles. In those works, 
we defined these categories for generic stability conditions. 
In our situation, the genericity of stability conditions is 
defined as follows: 
\begin{defn}\label{def:generic}
For $(\beta, \chi) \in \Lambda$ with $\beta$ effective, a stability condition $\sigma=\sigma_{B, H}$ is \textit{generic} with respect to $(\beta, \chi)$ if, for any decomposition in $\Lambda$
\begin{align*}(\beta, \chi)=(\beta_1, \chi_1)+(\beta_2, \chi_2), \ 
\mu^{\sigma}(\beta_i, \chi_i)=\mu^{\sigma}(\beta, \chi)
\end{align*}
with non-zero effective curve classes $\beta_i$, we have 
$(\beta_i, \chi_i) \in \mathbb{Q}_{>0}(\beta, \chi)$.  
\end{defn}

For a generic stability condition, the following lemma simplifies the definition of the quasi-BPS categories: 
\begin{lemma}\label{lem:generic}
Suppose that $\sigma$ is generic with respect to $(\beta, \chi) \in \Lambda$. Then for $\gamma=(r, D) \in \Gamma_{\mathbb{Q}}$, the dg-category 
\begin{align*}
    \LL(\mathfrak{M}_S(\beta, \chi)^{\sigma\text{-ss}})_{\delta_N(\gamma)}
\end{align*}
    only depends on $\gamma$ through 
    the central weight of $\delta_N(\gamma)$, that is, 
    $r \chi+D\beta$.
    In particular, we have 
\begin{align*}
    \mathbb{T}^{+}_S(\beta, \chi)_w^{\sigma}=\mathbb{T}^{-}_S(\beta, \chi)_w^{\sigma} 
    =\LL(\mathfrak{M}_S(\beta, \chi)^{\sigma\text{-ss}})_{\delta_N(\gamma)}
\end{align*}
for any choice of $\gamma=(r, D)\in \Gamma_{\mathbb{Q}}$ such that $w=r\chi+D\beta$.
\end{lemma}
\begin{proof}
    It is enough to show that 
    \begin{align*}
        \delta_N(\gamma) \in \mathrm{Pic}(\mathfrak{M}_S(\beta, \chi)^{\sigma\text{-ss}})_{\mathbb{Q}}^{\mathrm{CL}}
    \end{align*}
    depends only on its central weight if $\sigma$ is generic. 
    Let $\nu \colon \bgm \to \mathfrak{M}_S(\beta, \chi)^{\sigma\text{-ss}}$ be a map, corresponding to the direct 
    sum 
    \begin{align*}
       F= F_1 \oplus \cdots \oplus F_k, \ F_i \in \Coh^{\heartsuit}_{\leq 1}(S)
    \end{align*}
    with $\mathbb{G}_m$-weights $(\nu_1, \ldots, \nu_k)$ 
    such that $\cl(F_i)=(\beta_i, \chi_i)$.
    Since the above direct sum is $\sigma$-semistable, we have 
    $\mu^{\sigma}(\beta_i, \chi_i)=\mu^{\sigma}(\beta, \chi)$. 
    By the genericity of $\sigma$, we have $(\beta_i, \chi_i)=t_i(\beta, \chi)$ for $t_i \in \mathbb{Q}_{>0}$. 
By writing $\gamma=\cl(\gamma')$ for $\gamma' \in K(S)_{\mathbb{Q}}$, we have 
    \begin{align*}
        \delta_N(\gamma)(\nu)=\sum_{i=1}^k \chi(\gamma'\otimes F_i) \nu_i =\sum_{i=1}^k t_i \chi(\gamma'\otimes F)\nu_i
    \end{align*}
    and $\chi(\gamma'\otimes F)=r \chi+D\beta$ is the central weight of $\delta_N(\gamma)$. 
\end{proof}

We also give the following definition: 
\begin{defn}\label{def:primitive}
We say that $(\beta, \chi) \in \Lambda$ with $\beta$ effective is \textit{$\sigma$-primitive} 
if there is no decomposition 
\begin{align*}
    (\beta, \chi)=(\beta_1, \chi_1)+(\beta_2, \chi_2), \ 
    \mu^{\sigma}(\beta_i, \chi_i)=\mu^{\sigma}(\beta, \chi)
\end{align*}
    such that $\beta_i$ are non-zero and effective. 
\end{defn}
If $(\beta, \chi) \in \Lambda$ is $\sigma$-primitive, 
then there is no strictly $\sigma$-semistable sheaf 
$E \in \Coh^{\heartsuit}_{\leq 1}(S)$ with $\cl(E)=(\beta, \chi)$. 
In this case, we have 
\begin{align*}
    \mathbb{T}^{+}_S(\beta, \chi)_w^{\sigma}=\mathbb{T}^{-}_S(\beta, \chi)_w^{\sigma}
    =\Coh(\mathfrak{M}_S(\beta, \chi)^{\sigma\text{-st}})_w
    \simeq \Coh(M_S(\beta, \chi)^{\sigma\text{-st}}, \alpha^w)
\end{align*}
where $\mathfrak{M}_S(\beta, \chi)^{\sigma\text{-st}}$
is a $\mathbb{G}_m$-gerbe over $M_S(\beta, \chi)^{\sigma\text{-st}}$ with the corresponding Brauer class $\alpha$.

\subsection{Semiorthogonal decompositions into quasi-BPS categories}\label{subsec:sodqbps}
The following result, which generalizes~\cite[Theorem~1.1]{PTK3} and~\cite[Theorem~1.1]{PThiggs}, is obtained by applying the semiorthogonal decomposition~(\ref{sod:coh-limit-review}) and is essentially proved in~\cite[Section~4.4.7]{BPT}:

\begin{thm}\emph{(\cite{BPT})}\label{thm:bpt}
Let $S$ be a smooth quasi-projective surface, $(\beta, w)\in \Lambda$, $\chi \in \mathbb{Z}$ and $\sigma=\sigma_{B, H}$. 
    There is a semiorthogonal decomposition 
    \begin{align}\label{sod:1}
        \Coh(\mathfrak{M}_S(\beta, w)^{\sigma\text{-ss}})_{-\chi_{\beta}}
        =\left\langle \bigotimes_{i=1}^k \mathbb{T}^{-}_S(\beta_i, w_i)_{-\chi_{i, \beta_i}}^{\sigma}
        \right\rangle.
    \end{align}
    Here the right-hand side ranges over the partitions 
    \begin{align*}(\beta, \chi, w)=(\beta_1, \chi_1, w_1)+\cdots+(\beta_k, \chi_k, w_k)
    \end{align*}
    such that 
    \begin{align*}
        \mu^{\sigma}(\beta, w)=\mu^{\sigma}(\beta_i, w_i), \ 
        \mu^{\sigma}(\beta_1, \chi_1)>\cdots>\mu^{\sigma}(\beta_k, \chi_k), 
    \end{align*}
    and $(\chi, \chi_i) \mapsto (\chi_{\beta}, \chi_{i, \beta_i})$
    is the transformation rule (\ref{chi:trans}). 
    The semiorthogonal decomposition is ordered by decreasing values of 
    \begin{align}\label{obeta}
        o(\{(\beta_i, \chi_i)\}_{1\leq i\leq k}):=\sum_{i=1}^k (H\cdot \beta_i) \mu^{\sigma}(\beta_i, \chi_i)^2.
    \end{align}
\end{thm}
\begin{proof}
The application of the semiorthogonal decomposition~(\ref{sod:coh-limit-review}) in the above setting
is given in~\cite[Section~4.4.7]{BPT}; for each 
\begin{align*}\delta \in \mathrm{CL}(\mathfrak{M}_S(\beta, w))_{\mathbb{Q}}^{\vee}, \ 
r=\operatorname{Im}Z_{B, H} \colon  \Lambda \to \mathbb{Q}, 
\end{align*}
we have the semiorthogonal decomposition 
    \begin{align}\label{sod:LM}
         \Coh(\mathfrak{M}_S(\beta, w)^{\sigma\text{-ss}})
        =\left\langle \bigotimes_{i=1}^k \LL(\mathfrak{M}_S(\beta_i, w_i)^{\sigma\text{-ss}})_{\delta_i+\mu_i\ell_i} 
      : \mu_1<\cdots<\mu_k  \right\rangle.
    \end{align}
    Here the right-hand side ranges over the partitions 
    \begin{align*}(\beta, w)=(\beta_1, w_1)+\cdots+(\beta_k, w_k), \ \mu^{\sigma}(\beta_i, w_i)=\mu^{\sigma}(\beta, w).
    \end{align*}
    The elements $\delta_i, \ell_i \in  \mathrm{CL}(\mathfrak{M}_S(\beta_i, w_i)^{\sigma\text{-ss}})_{\mathbb{Q}}^{\vee}$ are given as follows. 
     Let $\nu \colon \bgm \to \mathfrak{M}_S(\beta_i, w_i)^{\sigma\text{-ss}}$ be a map
    \begin{align*}
        \mathrm{pt} \mapsto \bigoplus_{m\in \mathbb{Z}}E_{i, m}, \ \cl(E_{i, m})=(\beta_{i, m}, w_{i, m}).
    \end{align*}
    Then $\delta_i$ is given by 
    \begin{align*}
        \delta_i(\nu)&=\delta(\nu)+\frac{1}{2}\sum_{m\in \mathbb{Z}}m\cdot \chi\left((\beta_{i, m}, w_{i, m}), \sum_{j>i}(\beta_j, w_j)-\sum_{j<i}(\beta_j, w_j) \right), \\
        &=\delta(\nu)+\frac{1}{2}\sum_{m\in \mathbb{Z}}m\cdot\beta_{i, m}\left( \sum_{j<i}\beta_j-\sum_{j>i}\beta_j \right),
    \end{align*}
    and $\ell_i$ is given by  
    \begin{align*}
        \ell_i(\nu)&=\sum_{m\in \mathbb{Z}}m \cdot r(\beta_{i, m}, w_{i, m}) =\sum_{m\in \mathbb{Z}} m \cdot H\beta_{i, m}.
    \end{align*}
The order of the semiorthogonal decomposition is given by the decreasing values of 
\begin{align}\label{ord:r}
    \sum_{i=1}^k r(\beta_i, w_i) \mu_i^2 =\sum_{i=1}^k (H \cdot \beta_i) \mu_i^2 \in \mathbb{Q}_{\geq 0}.
\end{align}
We have also used the K\"unneth formula 
of limit categories in~\cite[Theorem~4.27]{TodaLimits}.
    
We take $\delta$ to be 
\begin{align*}
    \delta=\delta_N\left( -\frac{\beta}{2}-B+\frac{\beta^2}{2H\cdot \beta}H+\frac{B\cdot \beta}{H\cdot \beta}H \right).
\end{align*}
    Then we have 
    \begin{align*}
        \delta_i+\mu_i \ell_i=\delta_N\left(-\frac{\beta}{2}-B 
        +\frac{1}{2}\left(\sum_{j<i}\beta_j -\sum_{j>i}\beta_j  \right)-\mu_i' H\right).
    \end{align*}
    Here $\mu_i'$ is given by 
    \begin{align}\label{mui'}
     \mu_i'=-\left(\mu_i   +\frac{\beta^2}{2H\cdot \beta}+\frac{B\cdot \beta}{H\cdot \beta}\right).
    \end{align}
    We can write $\mu_i'=\mu^{\sigma}(\beta_i, \chi_i)$ for some $\chi_i\in \mathbb{Q}$. By $\mu_1<\cdots<\mu_k$, we have 
    \begin{align*}
        \mu^{\sigma}(\beta_1, \chi_1)>\cdots>\mu^{\sigma}(\beta_k, \chi_k).
    \end{align*}
We have the equality 
\begin{align*}
    -\frac{\beta}{2}-B 
        +\frac{1}{2}\left(\sum_{j<i}\beta_j -\sum_{j>i}\beta_j  \right)-\mu_i' H
        =-\mu^{\sigma}(\beta_i, \chi_i)H-\frac{\beta_i}{2}-B-\sum_{j>i}\beta_j.
\end{align*}
Thus $\delta_i+\mu_i\ell_i=
\delta_N(-D_{\sigma,\beta_i,\chi_{i,\beta_i}}-\sum_{j>i}\beta_j)$.
Tensoring by $\mathcal L_{\sum_{j>i}\beta_j}$ and using the definition of
$\mathbb T^{-}_S$, we obtain an equivalence
\begin{align}\label{Lequiv}
\otimes \cL_{\sum_{j>i}\beta_j} \colon 
  \LL(\mathfrak{M}_S(\beta_i, w_i)^{\sigma\text{-ss}})_{\delta_i+\mu_i \ell_i}
    \simeq \mathbb{T}^{-}_S(\beta_i, w_i)_{-\chi_{i, \beta_i}}^{\sigma}.
\end{align}
    
By (\ref{formula:D}), the central weight of $\delta_i+\mu_i\ell_i$ is 
    \begin{align*}
        \left(-\mu^{\sigma}(\beta_i, \chi_i)H-\frac{\beta_i}{2}-B-\sum_{j>i}\beta_j \right)\beta_i=
        -\chi_{i, \beta_i}-\beta_i\left(\sum_{j>i}\beta_j  \right). 
    \end{align*}
    If the dg-categories (\ref{Lequiv}) are non-zero, 
    then the above is an integer. Hence 
    $\chi_{i, \beta_i} \in \mathbb{Z}$, equivalently 
     $\chi_i \in \mathbb{Z}$. Moreover, the condition 
    \begin{align*}
        -\sum_{i=1}^k \left( \chi_{i, \beta_i}+\beta_i\left(\sum_{j>i}\beta_j \right)\right)=-\chi_{\beta}
    \end{align*}
    is equivalent to $\chi=\chi_1+\cdots+\chi_k$. 
    Therefore the result holds by taking the central weight $-\chi_{\beta}$-part of (\ref{sod:LM}).
    The order of the semiorthogonal decomposition is given by (\ref{obeta}), 
    by substituting (\ref{mui'}) into (\ref{ord:r}) and using that $(\beta_1, \chi_1, w_1)+\cdots +(\beta_k, \chi_k, w_k)=(\beta, \chi, w)$.
\end{proof}

\subsection{Semiorthogonal decomposition of the limit category}
We also consider the semiorthogonal decomposition of the limit 
category using Theorem~\ref{thm:limit-sod-general}. 

We first recall the Harder--Narasimhan stratification of $\mathfrak{M}_S(\beta, \chi)$ with respect to the stability condition 
$\sigma=\sigma_{B, H}$ 
in terms of numerical data $(\ell, b)$. 
Let $\alpha \in K(S)_{\mathbb{C}}$ be as in (\ref{classalpha}), and 
take the universal sheaf $\cE$ as in (\ref{univ:E}). 
Let $p_{\mathfrak{M}}$, $p_S$ be the projections from 
$S\times \mathfrak{M}_S(\beta, \chi)$ onto the corresponding factors. 
We then set 
\begin{align*}
    \ell:=\lvert Z_{B, H}(\beta, \chi)\rvert^2
    \det p_{\mathfrak{M}*}\left(\cE \otimes p_S^* \Im \left( \frac{-\alpha}{Z_{B, H}(\beta, \chi)} \right) \right) \in \mathrm{Pic}(\mathfrak{M}_S(\beta, \chi))_{\mathbb{Q}}
\end{align*}
and 
\begin{align*}
    b:=2\ch_2\left(p_{\mathfrak{M}*}(\cE \otimes p_S^* \Im \alpha)  \right) \in H^4(\mathfrak{M}_S(\beta, \chi), \mathbb{Q}).
\end{align*}
As a special case of~\cite[Section~6.4]{Halpinstab}, 
the above $(\ell, b)$ defines the $\Theta$-stratification 
\begin{align}\label{MSmu}
    \mathfrak{M}_S(\beta, \chi)=\bigsqcup_{\mu \in I_{>0}} \mathfrak{S}_{\mu} \sqcup \mathfrak{M}_S(\beta, \chi)^{\ell\text{-ss}}
\end{align}
which coincides with the Harder--Narasimhan stratification with respect to $\sigma$-stability; namely we have 
\begin{align*}
    \mathfrak{M}_S(\beta, \chi)^{\ell\text{-ss}}=\mathfrak{M}_S(\beta, \chi)^{\sigma\text{-ss}},
\end{align*}
and each $\mathfrak{S}_{\mu}$ corresponds to one-dimensional sheaves $E$
with a Harder--Narasimhan filtration with respect to $\sigma$
\begin{align}\label{HN:E2}
    0=E_0\subset E_1 \subset E_2 \subset \cdots \subset E_k=E.
\end{align}
Here each $F_i=E_i/E_{i-1}$ is $\sigma$-semistable
and by setting $\cl(F_i)=(\beta_i, \chi_i)$, 
we have 
\begin{align}\label{mu:ineq}
    \mu^{\sigma}(\beta_1, \chi_1)>\cdots>\mu^{\sigma}(\beta_k, \chi_k)
\end{align}
and we have the following equality (see~\cite[Lemma~6.4.17]{Halpinstab}):
\begin{align}\label{id:mu}
\mu
= (H\cdot\beta)
\left(
  \sum_{i=1}^k
  (H\cdot\beta_i)\,
  \mu^{\sigma}(\beta_i,\chi_i)^2
  -
  (H\cdot\beta)\,
  \mu^{\sigma}(\beta,\chi)^2
\right)^{1/2}.
\end{align}
We call the right-hand side of (\ref{id:mu}) the \textit{instability value} of the Harder--Narasimhan 
filtration (\ref{HN:E2}), and denote it by $\mu_{\mathrm{HN}}$. 
\begin{lemma}\label{lem:locfin}
For a fixed $(\beta, \chi)$, the set $I \subset \mathbb{R}_{\geq 0}$ of possible values 
(\ref{id:mu}) is locally finite. 
\end{lemma}
\begin{proof}
    Since $\beta-\beta_i$ is an effective divisor for each $i$, there are only finitely many possible values of $k\geq 1$ and $\beta_i$ for $1\leq i\leq k$. If (\ref{id:mu}) is bounded above, then the definition of $\mu^{\sigma}(\beta_i, \chi_i)$ also gives only finitely many possible values of $\chi_i$. Therefore $I$ is locally finite. 
\end{proof}

Let $\mathfrak{S}_{\mu} \to \mathfrak{Z}_{\mu}$ be the center of $\mathfrak{S}_{\mu}$. A $k$-point of $\mathfrak{Z}_{\mu}$ consists of 
a direct sum 
\begin{align}\label{sum:F}
    F=F_1 \oplus \cdots \oplus F_k
\end{align}
where $F_i$ are $\sigma$-semistable sheaves as above, with a 
canonical map (up to composition with 
$(-)^n \colon \mathbb{G}_m \to \mathbb{G}_m$) 
\begin{align*}
    \mathbb{G}_m \to \mathrm{Aut}(F), \ 
    t\mapsto (t^{\nu_1}, \ldots, t^{\nu_k}).
\end{align*}
By~\cite[Lemma~6.4.17]{Halpinstab}, we have 
\begin{align}\label{weight:nu}
    (\nu_1, \ldots, \nu_k) \in \mathbb{R}_{>0}(\mu^{\sigma}(\beta_1, \chi_1)-\mu^{\sigma}(\beta, \chi), \ldots, \mu^{\sigma}(\beta_k, \chi_k)-\mu^{\sigma}(\beta, \chi)).
\end{align}

By applying Theorem~\ref{thm:limit-sod-general} to the above situation, we obtain the following result: 
\begin{thm}\label{thm:sodLS}
Suppose that $-K_S$ is nef. 
Then we have the semiorthogonal decomposition
\begin{align}\label{sod:2}
    \LL(\mathfrak{M}_S(\beta, \chi))_{w_{\beta}}^{\sigma}=\left\langle \bigotimes_{i=1}^k\mathbb{T}^{+}_S(\beta_i, \chi_i)_{w_{i, \beta_i}}^{\sigma}\right\rangle.
\end{align}
Here the right-hand side ranges over the partitions 
\begin{align*}(\beta, \chi, w)=(\beta_1, \chi_1, w_1)+\cdots+(\beta_k, \chi_k, w_k)
\end{align*}
such that 
\begin{align*}
    \mu^{\sigma}(\beta, w)=\mu^{\sigma}(\beta_i, w_i), \ 
    \mu^{\sigma}(\beta_1, \chi_1)>\cdots>\mu^{\sigma}(\beta_k, \chi_k).
\end{align*}
The semiorthogonal decomposition is ordered by 
decreasing values of 
\begin{align}\label{order:o}
      o(\{(\beta_i, \chi_i)\}_{1\leq i\leq k})=\sum_{i=1}^k (H\cdot \beta_i) \mu^{\sigma}(\beta_i, \chi_i)^2.
\end{align}
\end{thm}
\begin{proof}
By Lemma~\ref{lem:checkassum} below, for each $\mu \in I_{>0}$, 
the $\Theta$-stratum 
\begin{align}\label{theta:mu}
    \mathfrak{S}_{\mu} \hookrightarrow\mathfrak{M}_S(\beta, \chi)_{\preceq \mu}
\end{align}
satisfies Assumption~\ref{assum:dag}. 
By applying Theorem~\ref{thm:limit-sod-general}, noting (\ref{formula:D}) and
the K\"unneth formula for limit categories~\cite[Theorem~4.27]{TodaLimits}, we obtain the 
semiorthogonal decomposition 
\begin{align}\label{sod:limitLM}
    &\LL(\mathfrak{M}_S(\beta, \chi))_{w_{\beta}}^{\sigma} \\
  \notag  &=
    \left\langle \bigotimes_{i=1}^k \LL(\mathfrak{M}_S(\beta_i, \chi_i)^{\sigma\text{-ss}})_{\delta_N(\mu^{\sigma}(\beta, w)H+\beta/2+B+\frac{1}{2}\left(\sum_{i>j}\beta_j-\sum_{i<j}\beta_j \right))} \right\rangle.
\end{align}
Here the right-hand side ranges over partitions $(\beta, \chi)=(\beta_1, \chi_1)+\cdots+(\beta_k, \chi_k)$ such that 
\begin{align*}
    \mu^{\sigma}(\beta_1, \chi_1)>\cdots>\mu^{\sigma}(\beta_k, \chi_k).
\end{align*}
The term $(\sum_{i>j}\beta_j-\sum_{i<j}\beta_j)/2$
is derived from $\delta_P(\det \mathbb{L}_{p_{\mu}}^{1/2})$
by the formula 
\begin{align*}
    \mathbb{L}_{p_{\mu}}|_{F}=\bigoplus_{i<j}\Hom^*(F_i, F_j)^{\vee}
\end{align*}
where $F$ is a direct sum as in (\ref{sum:F}). 
Indeed, for any map 
\begin{align*}\nu' \colon \bgm \to \mathfrak{Z}_{\mu}, \ 
\mathrm{pt} \mapsto F=F_1\oplus \cdots \oplus F_k
\end{align*}
with $\mathbb{G}_m$-weights 
$(\nu_1', \ldots, \nu_k')$, we have 
\begin{align*}
    c_1((\nu')^* \det \mathbb{L}_{p_{\mu}}^{1/2})=\frac{1}{2}\sum_{i<j}\beta_i \beta_j(\nu_j'-\nu_i') 
    =\frac{1}{2} \delta_N\left( \sum_{i>j}\beta_j-\sum_{i<j}\beta_j \right)(\nu').
\end{align*}

Below, we rewrite the semiorthogonal decomposition (\ref{sod:limitLM}). 
By writing $\mu^{\sigma}(\beta, w)=\mu^{\sigma}(\beta_i, w_i)$ for some $w_i\in \mathbb{Q}$, we have 
\begin{align*}
    \mu^{\sigma}(\beta, w)H+\frac{\beta}{2}+B+\frac{1}{2}\left(\sum_{i>j}\beta_j-\sum_{i<j}\beta_j \right)
    =\mu^{\sigma}(\beta_i, w_i)H+\frac{\beta_i}{2}+B+\sum_{i>j}\beta_j.
\end{align*}
Therefore, by (\ref{formula:D}), we have an equivalence 
\begin{align*}
\otimes\mathcal{L}_{\sum_{i>j}\beta_j} \colon &
\mathbb{T}^{+}_S(\beta_i, \chi_i)_{w_{i, \beta_i}}^{\sigma} \\
&\simeq 
    \LL(\mathfrak{M}_S(\beta_i, \chi_i)^{\sigma\text{-ss}})_{\delta_N(\mu^{\sigma}(\beta, w)H+\beta/2+B+\frac{1}{2}\left(\sum_{i>j}\beta_j-\sum_{i<j}\beta_j \right))}.
\end{align*}
It follows that $w_i \in \mathbb{Z}$ if the above is non-zero. 
The central weight of the right-hand side is 
\begin{align*}
   \left( \mu^{\sigma}(\beta_i, w_i)H+\frac{\beta_i}{2}+B+\sum_{i>j}\beta_j\right)\beta_i=w_{i, \beta_i}+\beta_i\left(\sum_{i>j}\beta_j\right).
\end{align*}
Comparing the central weights of both sides in (\ref{sod:limitLM}), we obtain 
\begin{align*}
    \sum_{i=1}^k w_{i, \beta_i}+\sum_{i>j}\beta_i\beta_j=w_{\beta}
\end{align*}
and hence $w=w_1+\cdots+w_k$. 

 From the construction, the order of the semiorthogonal decomposition is the same as the order of the Harder--Narasimhan stratification. Since it is 
    given by the instability values in (\ref{id:mu}), the order of the semiorthogonal 
    decomposition is given by (\ref{order:o}). 
\end{proof}

We have used the following lemma:
\begin{lemma}\label{lem:checkassum}
Each $\Theta$-stratum (\ref{theta:mu}) satisfies Assumption~\ref{assum:dag}.  
\end{lemma}
\begin{proof}
    For (i), the quasi-smoothness follows from the fact that $S$ is a surface~\cite[Section~3.4.1]{T}. 
The QCA condition follows from Lemma~\ref{lem:locfin} together with the fact that each Harder--Narasimhan stratum is of finite type. The local finite presentation property is a standard fact for moduli stacks of sheaves. The stack $\mathfrak{M}_S(\beta,\chi)_{\preceq\mu}$ has an affine diagonal, hence a separated diagonal, 
and the symmetric condition follows from Proposition~\ref{prop:sym}.

For (ii), let 
$x\in \mathfrak{M}_S(\beta, \chi)(k)$ be the point corresponding to the 
direct sum $F=\oplus_{i=1}^k F_i$ as in (\ref{sum:F}). 
We have 
\begin{align*}
   \mathbb{T}_{\mathfrak{M}}|_{x}=\bigoplus_{i, j}\Hom^*(F_i, F_j)[1].
\end{align*}
For $i<j$, we have 
\begin{align*}
    \Hom(F_i, F_j)=\Hom(F_i, F_j\otimes \omega_S)=0
\end{align*}
by Lemma~\ref{lem:S}. 
Therefore 
\begin{align*}
    \cH^{-1}(\mathbb{L}_{\mathfrak{M}}|_{x})=\bigoplus_{i\geq j}
    \Hom(F_i, F_j\otimes \omega_S).
\end{align*}
For the one-parameter subgroup $\lambda$ with weights $(\nu_1,\ldots,\nu_k)$ defining the stratum, each summand has $\lambda$-weight $\nu_j-\nu_i \geq 0$. Therefore (ii) is satisfied. 

For (iii), we have 
\begin{align*}
    \mathfrak{Z}_{\mu}=\bigsqcup_{(\beta_i, \chi_i)}\times_{i=1}^k \mathfrak{M}_S(\beta_i, \chi_i)^{\sigma\text{-ss}}
\end{align*}
where the right-hand side ranges over the partitions 
\begin{align*}(\beta, \chi)=(\beta_1, \chi_1)+\cdots+(\beta_k, \chi_k)
\end{align*}
such that the conditions (\ref{mu:ineq}) and (\ref{id:mu}) are satisfied.
Note that there are only finitely many such 
decompositions. 
It admits a good moduli space, so (iii) is satisfied. 
\end{proof}
   The assumption that $-K_S$ is nef is used in the following lemma, 
    which was used in the proof of Lemma~\ref{lem:checkassum}:
    \begin{lemma}\label{lem:S}
    Suppose that $-K_S$ is nef. 
    For $B+iH \in N^1(S)_{\mathbb{C}}$ such that $H$ is ample
    and $\sigma=\sigma_{B, H}$, let $F_1, F_2 \in \Coh_{\leq 1}^{\heartsuit}(S)$
    be $\sigma$-semistable with $\mu^{\sigma}(F_1)>\mu^{\sigma}(F_2)$. 
    Then we have 
    \begin{align}\label{vanish}
     \Hom(F_1, F_2)=   \Hom(F_1, F_2 \otimes \omega_S)=0.
    \end{align}
    \end{lemma}
    \begin{proof}
The first vanishing 
holds without the assumption on $-K_S$;
this is a consequence of the semistability of $F_i$ and $\mu^{\sigma}(F_1)>\mu^{\sigma}(F_2)$. 

We show the second vanishing when $-K_S$ is nef. 
Suppose that there is a non-zero morphism 
$\phi \colon F_1 \otimes \omega_S^{-1}\to F_2$. We may assume that both $F_1$ and $F_2$ are pure one-dimensional. 
Then $A=\mathrm{Im}(\phi) \subset F_2$ is a pure one-dimensional sheaf. 
By the $\sigma$-semistability of $F_1, F_2$, 
together with the injection $A \hookrightarrow F_2$
and surjection $F_1 \twoheadrightarrow A\otimes \omega_S$, 
we have the inequalities 
\begin{align*}
    \mu^{\sigma}(F_1) \leq \mu^{\sigma}(A\otimes \omega_S), \ 
    \mu^{\sigma}(A) \leq \mu^{\sigma}(F_2). 
\end{align*}
On the other hand, as $-K_S$ is nef, we have 
\begin{align*}
    \mu^{\sigma}(A\otimes \omega_S)=
    \mu^{\sigma}(A)+\frac{K_S \cdot [A]}{H\cdot [A]}
    \leq \mu^{\sigma}(A).
\end{align*}
Therefore $\mu^{\sigma}(F_1) \leq \mu^{\sigma}(F_2)$, which contradicts $\mu^{\sigma}(F_1)>\mu^{\sigma}(F_2)$. 
    \end{proof}

We expect that the semiorthogonal decomposition (\ref{sod:1})
corresponds to (\ref{sod:2}) under a conjectural equivalence in Conjecture~\ref{conj:CY3}. 
In particular, we conjecture the following: 
\begin{conj}\label{conj:qbps}
    There is an equivalence 
    \begin{align}\label{equiv:qbps}
        \mathbb{T}^{-}_S(\beta, w)_{-\chi_{\beta}}^{\sigma} \simeq \mathbb{T}^{+}_S(\beta, \chi)_{w_{\beta}}^{\sigma}.
    \end{align}
\end{conj}

\subsection{Numerical Picard rank one case}
In Conjectures~\ref{conj:CY3} and~\ref{conj:qbps}, the relevant categories depend on 
a choice of a stability condition $\sigma$ in general. 
In the numerical Picard rank one case, they are independent of $\sigma$. 
Indeed, we have the following lemma: 
\begin{lemma}\label{lem:indep}
   Let $H$ be an ample divisor on $S$. For $a, b \in \mathbb{Q}$ with $a>0$, let $\sigma=\sigma_{bH, aH}$, $\sigma_0=\sigma_{0, H}$.  
   Then we have 
   \begin{align*}
       \mathfrak{M}_S(\beta, \chi)^{\sigma\text{-ss}}=\mathfrak{M}_S(\beta, \chi)^{\sigma_0\text{-ss}}, \ D_{\sigma, \beta, w}=D_{\sigma_0, \beta, w}.
   \end{align*}
\end{lemma}
\begin{proof}
    Both follow from the definitions by direct computation. 
\end{proof}

Suppose that 
\begin{align*}N^1(S)_{\mathbb Q}=\mathbb Q[H]
\end{align*}for an ample divisor $H$ on $S$. 
In this case $D_{\sigma,\beta,-w}=-D_{\sigma,\beta,w}$, so
$\mathbb T^{+}_S=\mathbb T^{-}_S$; we omit the sign.
We write
\begin{align*}
    \mathfrak{M}_S(\beta, \chi)^{\text{ss}}:=\mathfrak{M}_S(\beta, \chi)^{\sigma\text{-ss}}
    \end{align*}
    and 
    \begin{align*}
    \LL(\mathfrak{M}_S(\beta, \chi))_w :=\LL(\mathfrak{M}_S(\beta, \chi))_{w}^{\sigma}, \
    \mathbb{T}_S(\beta, \chi)_w :=\mathbb{T}_S(\beta, \chi)_w^{\sigma}
\end{align*}
since they are independent of $\sigma=\sigma_{bH, aH}$ by Lemma~\ref{lem:indep}.
The equivalence in Conjecture~\ref{conj:CY3} can then be written as 
\begin{align}\label{conj:picone}
\Coh(\mathfrak{M}_S(\beta, w)^{\mathrm{ss}})_{-\chi_{\beta}} \simeq \LL(\mathfrak{M}_S(\beta, \chi))_{w_{\beta}}.  
\end{align}

\section{Examples and discussion}\label{sec:examples}
In this section, we give examples of Conjecture~\ref{conj:CY3} and discuss their
importance and consequences. 
\subsection{Over the locus of integral curves}

Let $\rB^\circ\subset\rB$ be the derived open subscheme
parametrizing integral divisors in $S$, and set
\begin{align}\label{map:hcirc}
h^{\circ} \colon 
\mathfrak M_S^\circ(\beta,\chi)
 :=\mathfrak M_S(\beta,\chi)\times_{\rB}\rB^\circ \to \rB^{\circ}. 
\end{align}
A sheaf parametrized by $\mathfrak M_S^\circ(\beta,\chi)$
is a rank-one torsion-free sheaf on its support (which is an integral curve), and is therefore
$\sigma$-stable. In particular, the above stack is quasi-compact.
Then
\[
\LL(\mathfrak M_S^\circ(\beta,\chi))_w^{\sigma}
 =\Coh(\mathfrak M_S^\circ(\beta,\chi))_w
 =\mathbb T_S(\beta,\chi)_w^{\sigma}\big|_{\rB^\circ}.
\]
Consequently, the restriction of Conjecture~\ref{conj:CY3}
to $\rB^\circ$ reduces to a $\rB^\circ$-linear equivalence
\begin{equation}\label{equiv:integ}
 \Coh(\mathfrak M_S^\circ(\beta,w))_{-\chi_\beta}
 \simeq
 \Coh(\mathfrak M_S^\circ(\beta,\chi))_{w_\beta}.
\end{equation}
\begin{prop}\label{prop:integral}
   There is a $\rB^{\circ}$-linear equivalence
   (\ref{equiv:integ}). 
\end{prop}
\begin{proof}
We may assume that $S$ is projective. 
By taking a closed immersion $S \hookrightarrow \mathbb{P}^N$ for $N\gg 0$, 
the universal curve over $\rB^\circ$
defines a morphism to the derived Hilbert scheme of
$\mathbb P^N$. As in~\cite[Lemma~4.2]{MT}, a neighborhood of
its image is a smooth classical scheme $T$ carrying a locally
versal family of integral planar curves $\cC_T\to T$.
If $\mathfrak J_T(\chi)$ denotes the stack of rank-one torsion-free
sheaves of Euler characteristic $\chi$ on its fibers, then
there is a natural morphism 
\begin{align}\label{construct:map}
\mathfrak J_T(\chi)\times_T\rB^{\circ} \to
 \mathfrak M_S^\circ(\beta,\chi).
\end{align}

We show that (\ref{construct:map}) is an equivalence. 
Since both are isomorphic in classical truncations, 
it is enough to check that the induced map on
$\rB^\circ$-relative cotangent complexes is an equivalence.
Since both sides are $\mathbb{G}_m$-gerbes over the good moduli spaces, and the induced maps on the central stabilizers are isomorphisms, it is enough to show that the map (\ref{construct:map}) is an equivalence over 
the good moduli spaces, which we denote by $u \colon J \to M$. 

For a sheaf $F$ supported on an integral curve $C$, the
induced map on the 
obstruction spaces by the map $M\to \rB^{\circ}$ is 
\[
 \operatorname{Ext}^2_S(F,F)
 \to H^1(C,\mathcal O_C(C)).
\]
By Serre duality, its dual is the natural map 
\[
 H^0(C,\omega_S|_C)
 \to
 \operatorname{Hom}_C(F,F\otimes \omega_S|_C),
\]
which is injective because $F$ is rank-one torsion-free.
Thus $\mathbb L_{M/\rB^\circ}$ is perfect of amplitude
$[-1,0]$, whose rank is 
\[
 \operatorname{rk}\mathbb L_{M/\rB^\circ}
 =1-\chi_S(F,F)-\chi(C,\mathcal O_C(C))
 =p_a(C).
\]

Set $\rB_0=(\rB^\circ)^{\mathrm{cl}}$ and let
$u_0\colon J_0\to M_0$ be the derived base change of $u$
to $\rB_0$. Since $\mathfrak{J}_T(\chi)\to T$ is flat, 
it follows that $J_0$ is classical, and
$u_0$ identifies it with $M_0^{\mathrm{cl}}$. Consequently,
\[
 u_0^*\mathbb L_{M_0/\rB_0}
 \to \mathbb L_{J_0/\rB_0}
\]
is a $\rB_0$-relative perfect obstruction theory on $J_0$. Both complexes have
amplitude $[-1,0]$ with the same virtual rank, so 
this morphism is a quasi-isomorphism. It follows that 
$\mathbb{L}_{J/M}|_{\rB_0}=0$, 
hence $\mathbb{L}_{J/M}=0$, and this implies that 
$u$ is an equivalence. 

The Cohen--Macaulay extension of the relative Poincar\'e 
line bundle~\cite{Ardual} on $\mathfrak{J}_T(w)\times_T \mathfrak{J}_T(\chi)$
gives a Fourier--Mukai equivalence
\begin{align*}
    \Coh(\mathfrak{J}_T(w))_{-\chi_{\beta}} \stackrel{\sim}{\to} 
      \Coh(\mathfrak{J}_T(\chi))_{w_{\beta}},
\end{align*}
see
\cite[Proposition~2.28]{TodaGL2}.
By taking the derived base change of the above Fourier--Mukai kernel
and its
inverse along $\rB^\circ\to T$, 
and using an equivalence (\ref{construct:map}), 
we obtain an equivalence~\eqref{equiv:integ}.
\end{proof}

The above proposition shows that Conjectures~\ref{conj:CY3} and
\ref{conj:qbps} hold over 
$\rB^{\circ} \subset \rB$. 
The locus $\rB^{\circ}$ corresponds to the \textit{elliptic locus}
in the case of moduli stacks of Higgs bundles. 
However, beyond $\rB^\circ$, these conjectures are substantially harder, as discussed below.
\subsection{Higgs bundles}
\subsubsection*{Moduli stacks of Higgs bundles}
Let $C$ be a smooth projective curve of genus $g$. We take 
$L \in \mathrm{Pic}(C)$ such that 
\begin{align*}
l:=\deg L>2g-2, \ \mbox{ or } L=\Omega_C.
\end{align*}
Let $S$ be the surface 
\begin{align*}
    S=\mathrm{Tot}_C(L) \to C. 
\end{align*}
Here $N^1(S)_{\mathbb Q}$ is one-dimensional, so the positive and negative
quasi-BPS categories agree; we omit the sign as in the numerical Picard rank one case.
Then for $\beta=r[C]$, where $C \hookrightarrow S$ is the zero section, 
and $\chi \in \mathbb{Z}$, we have 
\begin{align*}
    \mathfrak{M}_S(r[C], \chi) \simeq \Hig_{\GL_r}^L(\chi)
\end{align*}
where the right-hand side is the moduli stack of $L$-twisted $\GL_r$-Higgs bundles 
\begin{align*}
    (F, \theta), \ \theta \colon F \to F\otimes L 
\end{align*}
where $F \to C$ is a vector bundle of rank $r$ with $\chi(F)=\chi$. 
Moreover, the map $\mathfrak{M}_S(r[C], \chi) \to \rB$ is identified with the Hitchin fibration 
\begin{align*}
    \Hig_{\GL_r}^L(\chi) \to \rB=\bigoplus_{i=1}^r \Gamma(C, L^{\otimes i}).
\end{align*}

\subsubsection*{Dolbeault geometric Langlands}
In this case, 
Conjecture~\ref{conj:CY3} predicts a $\rB$-linear equivalence  
\begin{align}\label{equiv:L}
    \Coh(\Hig_{\GL_r}^L(w)^{\text{ss}})_{-\chi'} \simeq \LL(\Hig_{\GL_r}^L(\chi))_{w'}, 
\end{align}
where $(\chi', w')$ are given by 
\begin{align*}
    \chi'=\chi+(g-1)r+\frac{l}{2}r(r-1), \ 
    w'=w+(g-1)r+\frac{l}{2}r(r-1).
\end{align*}
For $L=\Omega_C$, 
both sides are unchanged by the replacement of weights $(-\chi', w') \mapsto (-\chi, w)$
(see~\cite[Remark~2.4]{TodaGL2}), 
and (\ref{equiv:L}) is nothing but the main conjecture in~\cite{PTlim} for $G=\GL_r$. 
The equivalence (\ref{equiv:L}) extends the formulation to the $L$-twisted case. 

For $L=\Omega_C$, a version of the equivalence (\ref{equiv:L}) with nilpotent singular supports~\cite{AG} is 
proved over the locus in $\rB$ corresponding to 
spectral curves which are either integral or have at worst type A singularities~\cite{TodaGL2, TodaA}.
In particular, the nilpotent singular support version holds for $r=2$ over the locus in $\rB$ 
corresponding to reduced spectral curves. 

For sufficiently positive $L$, in the forthcoming paper~\cite{PTtwist}, we prove the rank-two case over the full
Hitchin base:
\begin{thm}[\cite{PTtwist}]\label{thm:twisted-DL-GL2}
Suppose that $r=2$ and
$
 \deg L >2g$.
Then, for every $\chi,w\in\mathbb Z$, the $\rB$-linear equivalence
(\ref{equiv:L}) holds over the full Hitchin base $\rB$, including
the locus of non-reduced spectral curves.
\end{thm}

The conjectural equivalence 
\begin{align}\label{equiv:T}
    \mathbb{T}_S(r[C], w)_{-\chi'} \simeq \mathbb{T}_S(r[C], \chi)_{w'}
\end{align}
in the above case 
is nothing but the main conjecture in~\cite{PThiggs}. See~\cite[Section~1.3]{PThiggs} for the cases where 
the equivalence (\ref{equiv:T}) is proved. 
\subsubsection*{Topological K-theory}
For $\deg L>2g-2$, it is proved in~\cite{PThiggs2} that (\ref{equiv:T}) holds 
at the level of topological K-theory 
\begin{align}\label{isom:topK}
    K^{\mathrm{top}}( \mathbb{T}_S(r[C], w)_{-\chi'}) \stackrel{\cong}{\to} K^{\mathrm{top}}( \mathbb{T}_S(r[C], \chi)_{w'}). 
\end{align}

By combining Theorem~\ref{thm:bpt} and Theorem~\ref{thm:sodLS}, we also have an isomorphism 
compatible with the semiorthogonal decompositions in Theorems~\ref{thm:bpt} and~\ref{thm:sodLS}:
\begin{align*}
   K^{\mathrm{top}}(\Coh(\Hig_{\GL_r}^L(w)^{\text{ss}})_{-\chi'}) \simeq K^{\mathrm{top}}(\LL(\Hig_{\GL_r}^L(\chi))_{w'}).
   \end{align*}
\begin{remark}
    The isomorphism (\ref{isom:topK}) is also expected to hold for $L=\Omega_C$. In~\cite[Theorem~1.1]{PThiggs2}, the isomorphism (\ref{isom:topK}) for $L=\Omega_C$ (more precisely, its reduced version) is proved when $(r, \chi, w)\in \mathbb{Z}^3$ is primitive.  
\end{remark}
   \subsubsection*{Semiorthogonal decomposition: simplest case}
   Let $C=\mathbb{P}^1$ and $L=\mathcal{O}_{\mathbb{P}^1}(-1)$. Then any 
   semistable sheaf on $S$ is scheme-theoretically supported 
   on the zero section $C \subset S$. Therefore 
   \begin{align*}
       \mathfrak{M}_S(r[C], \chi)^{\text{ss}}\cong \begin{cases}
       B\GL_r, & r \mid \chi, \\
       \emptyset, & \mbox{otherwise.}
           \end{cases}
   \end{align*}
   By the above isomorphism, one can check that 
   \begin{align*}
       \mathbb{T}_S(r[C], \chi)_w \simeq 
       \begin{cases}
         \mathrm{Vect}, & r=1, \\
         0, & r>1.
       \end{cases}
   \end{align*}
   Here $\mathrm{Vect}=\Coh(\mathrm{pt})$. In particular, 
   the equivalence (\ref{equiv:T}) holds in this case. 

In this case, the semiorthogonal decomposition (\ref{sod:1}) for $w=0$ is 
\begin{align*}
    \Coh(B\GL_r)_{-\chi+r(r+1)/2}
    =\left\langle \mathrm{Vect}_{(\chi_1, \ldots, \chi_r)} \right\rangle.
\end{align*}
Here each component is a copy of $\mathrm{Vect}$, indexed by 
partitions 
\begin{align}\label{part:chi}
    \chi=\chi_1+\cdots+\chi_r, \ \chi_i \in \mathbb{Z}, \ \chi_1>\cdots>\chi_r.
\end{align}
Similarly, the semiorthogonal decomposition (\ref{sod:2}) for $w=0$ is 
\begin{align*}
    \LL(\mathfrak{M}_S(r[C], \chi))_{-r(r+1)/2}=\left\langle \mathrm{Vect}_{(\chi_1, \ldots, \chi_r)} \right\rangle
\end{align*}
where each component is a copy of $\mathrm{Vect}$, indexed by 
partitions (\ref{part:chi}).

\subsection{K3 surfaces}\label{subsec:K3}
\subsubsection*{Dolbeault geometric Langlands for K3 surfaces}
Let $S$ be a smooth projective K3 surface such that 
\begin{align*}
    \mathrm{Pic}(S)=\mathbb Z[\mathcal O_S(H)], \quad H\text{ ample},\quad H^2=2d. 
\end{align*}
In this case, Conjecture~\ref{conj:CY3} predicts a $\rB$-linear equivalence 
\begin{align*}
    \Coh(\mathfrak{M}_S(r[H], w)^{\text{ss}})_{-\chi -r^2 d}\simeq 
    \LL(\mathfrak{M}_S(r[H], \chi))_{w+r^2 d}.
\end{align*}
The conjectural equivalence of quasi-BPS categories in Conjecture~\ref{conj:qbps} is 
\begin{align}\label{equiv:BPSK3}
       \mathbb{T}_S(r[H], w)_{-\chi-r^2 d} \simeq \mathbb{T}_S(r[H], \chi)_{w+r^2 d}.
\end{align}
\subsubsection*{Beauville--Mukai system}
If $r$ and $\chi$ are coprime, then we have 
\begin{align*}
    \mathfrak{M}_S(r[H], \chi)^{\text{ss}}= \mathfrak{M}_S(r[H], \chi)^{\text{st}}.
\end{align*}
Here $(-)^{\mathrm{st}}$ refers to the stable part. 
In this case, it 
is easy to check that
\begin{align*}
     \mathfrak{M}_S(r[H], \chi)^{\text{st}}\simeq 
      \mathcal{M}_S(r[H], \chi)^{\text{st}}\times k[-1]
\end{align*}
where $\mathcal{M}_S(r[H], \chi)^{\text{st}}$ is the classical moduli stack of 
stable sheaves, and the good moduli space map 
\begin{align}\label{Gmgerb}
    \mathcal{M}_S(r[H], \chi)^{\text{st}}\to M_S(r[H], \chi)^{\text{st}}
\end{align}
is a $\mathbb{G}_m$-gerbe
over a holomorphic symplectic manifold $M_S(r[H], \chi)^{\text{st}}$, deformation equivalent to 
the Hilbert scheme of $r^2d +1$ points on $S$.
The map 
\begin{align*}
    \mathfrak{M}_S(r[H], \chi)^{\text{st}} \to \rB=\mathbb{P}(H^0(\mathcal{O}_S(r[H])))\times k[-1]
\end{align*}
is induced by the map 
\begin{align}\label{BMsystem}
M_S(r[H], \chi)^{\text{st}} \to \rB^{\mathrm{cl}}=\mathbb{P}^{r^2 d+1}.
\end{align}
This is a Lagrangian fibration (in the classical sense), and 
is called the \textit{Beauville--Mukai system}. 
\subsubsection*{Dolbeault geometric Langlands for the Beauville--Mukai system}
In the above situation, we have 
\begin{align*}
    \mathbb{T}_S(r[H], \chi)_{w'} &\simeq \Coh(\mathfrak{M}_S(r[H], \chi)^{\text{st}})_{w'} \\
   & \simeq \Coh(M_S(r[H], \chi)^{\text{st}}\times k[-1], \alpha_{\chi}^{w'})
\end{align*}
where $w'$ is given by 
\begin{align*}w'=w+\frac{(rH)^2}{2}=w+r^2d
\end{align*}
and $\alpha_{\chi}$ is the Brauer class corresponding to 
the $\mathbb{G}_m$-gerbe (\ref{Gmgerb}). 
Therefore, if both pairs $(r, \chi)$ and $(r, w)$ are coprime, then 
the equivalence (\ref{equiv:qbps}) follows from a conjectural derived equivalence 
over $\rB^{\mathrm{cl}}$
\begin{align}\label{equiv:BM}
     \Coh(M_S(r[H], w)^{\text{st}}, \alpha_{w}^{-\chi-r^2 d})
     \simeq  \Coh(M_S(r[H], \chi)^{\text{st}}, \alpha_{\chi}^{w+r^2 d}).
\end{align}

When $r=1$, each member $C\in \lvert H \rvert$ is an irreducible planar curve, and 
the moduli space $M_S([H], \chi)^{\text{st}}$ is isomorphic to the 
relative compactified Jacobian of the universal curve $\mathcal{C} \to \rB^{\mathrm{cl}}$:
\begin{align*}
    M_S([H], \chi)^{\text{st}} \cong \overline{J}^{\chi+d}(\cC/\rB^{\mathrm{cl}})
\end{align*}
which is the moduli space of rank-one torsion-free sheaves $F$ on $C \in \lvert H \rvert$
with degree $\chi+d$.
In this case, the equivalence (\ref{equiv:BM}) follows from~\cite{Ardual}; see~\cite[Proposition~2.28]{TodaGL2}.

\subsubsection*{The non-coprime case}
If $r$ and $\chi$ are not coprime, then the good moduli space map 
\begin{align*}
    \mathfrak{M}_S(r[H], \chi)^{\mathrm{ss, cl}} \to M_S(r[H], \chi)^{\mathrm{ss}}
\end{align*}
is not a $\mathbb{G}_m$-gerbe, and $M_S(r[H], \chi)^{\text{ss}}$ is a singular symplectic variety which does not admit a symplectic resolution of singularities except when $d=1, r=2$ (the exceptional case gives O'Grady's 10-dimensional symplectic manifold~\cite{Ogra1}). 

Nevertheless, if $(r, \chi, w) \in \mathbb{Z}^3$ is primitive, it is proved in~\cite[Theorem~1.2]{PTK3} that the dg-category (called a \textit{BPS category})
\begin{align*}\mathbb{T}_S(r[H], \chi)_w
\end{align*}
is smooth, proper and Calabi--Yau (up to removing trivial derived structure, \'etale locally on the good moduli space). Thus it gives a (twisted, \'etale locally)
categorical crepant resolution of $M_S(r[H], \chi)^{\text{ss}}$. 

However, it is not clear what the category $\mathbb{T}_S(r[H], \chi)_w$ actually looks like, e.g. 
whether it is equivalent to a (twisted) derived category of coherent sheaves on a holomorphic symplectic manifold. If the equivalence (\ref{equiv:BPSK3}) holds, then it 
answers this question, at least partially. Indeed, for any pair $(r, \chi)$, not necessarily coprime, the tuple $(r, \chi, 1)$ is primitive. Therefore 
\begin{align*}
    \mathbb{T}_S(r[H], \chi)_1 &\simeq \mathbb{T}_S(r[H], 1-r^2 d)_{-\chi-r^2 d} \\ &\simeq 
    \Coh(M_S(r[H], 1-r^2 d)^{\text{st}}\times k[-1], \alpha_{1-r^2 d}^{-\chi-r^2 d})
\end{align*}
since $r$ and $1-r^2 d$ are coprime. 

\begin{example}\label{example:double}
Let $S$ be a K3 surface, obtained as a double cover $S \to \mathbb{P}^2$ ramified over 
a generic sextic. Let $H$ be an ample divisor on $S$, given by the pull-back of a 
line in $\mathbb{P}^2$. The (twisted) Fourier--Mukai transform of the Beauville--Mukai system 
in this case is discussed in~\cite[Section~5]{Sawon2004}. 
In the notation of~\cite[Section~5.3]{Sawon2004}, we have $Y^k=M_S(2[H], k-4)^{\mathrm{st}}$. 
In~\cite[Section~5.3]{Sawon2004}, the following question is addressed: 
`can one define a dual fibration of $Y^5$, which should agree with $Y^0$ over smooth fibers, in a sensible way?'

If (\ref{equiv:BPSK3}) holds, we have  
\begin{align*}
\mathbb{T}_S(2[H], 1)_{0}  \simeq \mathbb{T}_S(2[H], -4)_5.
\end{align*}
The left-hand side is equivalent to $\Coh(Y^5) \otimes \Coh(k[-1])$. After removing the $k[-1]$-factor, 
the above equivalence answers this question; the BPS category $\mathbb{T}_S(2[H], -4)_5$
is regarded as a categorical dual of $Y^5$. 
\end{example}

\subsection{The projective plane}
\subsubsection*{Moduli stacks of semistable sheaves on $\mathbb{P}^2$}
Let $S=\mathbb{P}^2$. Note that we have 
\begin{align*}
    \mathrm{Pic}(\mathbb{P}^2)=\mathbb{Z}[\mathcal{O}(H)]
\end{align*}
for a hyperplane $H\subset \mathbb{P}^2$. 
Note that $\mathfrak{M}_S(r[H], \chi)^{\mathrm{ss}}$ is smooth
(since $\Ext_S^2(F, F)=0$ for one-dimensional semistable sheaves on $\mathbb{P}^2$), 
and the map to $\rB$ is 
\begin{align*}
    \mathfrak{M}_S(r[H], \chi)^{\mathrm{ss}}\to \rB=\mathbb{P}(H^0(\mathbb{P}^2, \mathcal{O}(rH))). 
\end{align*}
\subsubsection*{Conjectural equivalence for $\mathbb{P}^2$}
The conjectural equivalence in Conjecture~\ref{conj:CY3} is 
\begin{align}\label{equiv:P20}
    \Coh(\mathfrak{M}_S(r[H], w)^{\mathrm{ss}})_{-\chi-r(r-3)/2} \simeq \LL(\mathfrak{M}_S(r[H], \chi))_{w+r(r-3)/2}.
\end{align}
The conjectural equivalence of quasi-BPS categories is 
\begin{align}\label{equiv:P21}
    \mathbb{T}_S(r[H], w)_{-\chi -r(r-3)/2} \simeq \mathbb{T}_S(r[H], \chi)_{w+r(r-3)/2}.
\end{align}
The equivalences (\ref{equiv:P20}) and (\ref{equiv:P21}) are obvious for $r=1$; in this case 
both sides are equivalent to $\Coh(\rB)$. 

\subsubsection*{The case of $r=2$}
For $r=2$, both sides of (\ref{equiv:P21}) are invariant under 
$\chi \mapsto \chi+r$ and $w\mapsto w+r$. Thus, to prove (\ref{equiv:P21}), it is enough to check the equivalence
\begin{align*}
    \mathbb{T}_S(2[H], 0)_0 \simeq \mathbb{T}_S(2[H], 1)_1.
\end{align*}
One can directly check that both sides are equivalent to $\Coh(\rB)$, 
so they are equivalent. However, the equivalence (\ref{equiv:P20}) is 
non-trivial even for $r=2$: both sides admit semiorthogonal decompositions into 
quasi-BPS categories, and we need to compare how these categories are glued on the two sides. 

\subsection{Two $(-1)$-curves}\label{subsec:two-1}
Let 
\begin{align*}
    g\colon S \to \mathbb{P}^1 \times \mathbb{A}^1
\end{align*}
be the blow-up at $(0, 0) \in \mathbb{P}^1 \times \mathbb{A}^1$. 
We have the projection 
\begin{align*}
    f \colon S \to \mathbb{A}^1
    \end{align*}
    whose fiber at $t\in \mathbb{A}^1$ is $\mathbb{P}^1$ 
for $t\neq 0$, and the fiber at $0\in \mathbb{A}^1$ is $C_1 \cup C_2$ where 
$C_i \subset S$ are $(-1)$-curves, $C_1$ is the strict transform of 
$\mathbb{P}^1 \times \{0\}$ and $C_2$ is the exceptional curve of $g$. 

Let $H'\in N^1(\mathbb{P}^1 \times \mathbb{A}^1)$ be the pull-back of the hyperplane class of $\mathbb{P}^1$, and 
set 
\begin{align*}
    H=2g^* H'-C_2 \in N^1(S).
\end{align*}
It is an ample divisor such that $H\cdot C_1=H \cdot C_2=1$. 
By setting $\sigma=\sigma_{0, H}$ and $\beta=[C_1]+[C_2]$, we have 
\begin{align*}
    \mathfrak{M}_S(\beta, 0)^{\sigma\text{-ss}} \cong \mathbb{A}^2/\mathbb{G}_m \times \bgm
\end{align*}
where $\mathbb{G}_m$ acts on $\mathbb{A}^2$ with weights $(1, -1)$, 
and the closed point $0\in \mathbb{A}^2/\mathbb{G}_m$ corresponds to 
$\mathcal{O}_{C_1}(-1) \oplus \mathcal{O}_{C_2}(-1)$. 
It follows that 
\begin{align*}
      \Coh(\mathfrak{M}_S(\beta, 0)^{\sigma\text{-ss}})_{0} \simeq \Coh(\mathbb{A}^2/\mathbb{G}_m).
\end{align*}

Let $\X$ be the stack constructed in \cite[Subsection~5.1]{TodaLimits};
\begin{align*}
    \X=\bigcup_{i\in 1/2+\mathbb{Z}}\X_i, \ \X_i=\mathbb{A}^2/\mathbb{G}_m
\end{align*}
where $\X_i$ and $\X_{i+1}$ are glued by the open substacks isomorphic to $\mathbb{A}^1$. 
Then there is an isomorphism 
\begin{align*}
    \mathfrak{M}_S(\beta, 1) \cong \X \times \bgm
\end{align*}
such that each closed point $p_i \in \X_i$ with $\mathbb{G}_m$-stabilizer for $i\in 1/2 +\mathbb{Z}$ corresponds to 
\begin{align*}\mathcal{O}_{C_1}\left(i-\frac{1}{2}\right) \oplus \mathcal{O}_{C_2}\left(-i-\frac{1}{2}\right)\in \Coh_{\leq 1}^{\heartsuit}(S).
\end{align*}
The element $\delta_N(D_{\sigma, \beta, -1})$ corresponds to $\delta\boxtimes \delta_P(\mathcal{O}_{\bgm}(-1))$, where $\delta=\delta_P(\mathcal{O}_{\mathcal{X}}(-1/2))$ as 
in \cite[Subsection~5.1]{TodaLimits}. It follows that 
\begin{align*}
    \LL(\mathfrak{M}_S(\beta, 1))_{-1}^{\sigma} \simeq 
    \LL(\X)_{\delta}.
\end{align*}

It follows that, for $\chi=1$ and $w=0$ (so that $\chi_{\beta}=0$ and $w_{\beta}=-1$), the 
equivalence 
\begin{align*}
    \Coh(\mathfrak{M}_S(\beta, 0)^{\sigma\text{-ss}})_{0}\simeq \LL(\mathfrak{M}_S(\beta, 1))_{-1}^{\sigma}
\end{align*}
in Conjecture~\ref{conj:CY3} is a direct consequence of the following equivalence, 
see~\cite[Proposition~5.2]{TodaLimits}:
\begin{align*}
    \Coh(\mathbb{A}^2/\mathbb{G}_m) \simeq \LL(\X)_{\delta}.
\end{align*}

\section{Hecke operators}\label{sec:hecke}
In this section, we construct Hecke operators on limit categories and BPS categories, 
generalizing the construction in~\cite[Section~9]{PTlim}. 
\subsection{Moduli stacks of Hecke operators}
Let 
\begin{align*}
    \mathrm{Hecke}_S
\end{align*}
be the derived moduli stack which classifies exact sequences in $\Coh_{\leq 1}^{\heartsuit}(S)$
\begin{align}\label{ses:E}
0\to E_1 \to E_2 \to Q\to 0
\end{align}
where $E_i$ for $i=1, 2$ are pure one-dimensional 
sheaves and $Q$ is a zero-dimensional sheaf. 
It is a quasi-smooth stack; see~\cite[Proposition~3.8]{PoSa}.
We have the evaluation morphisms 
\begin{align}\label{dia:left}
\xymatrix{
    \mathrm{Hecke}_S \ar[r]^-{\ev_2} \ar[d]_-{(\ev_1, \ev_3)} & \mathfrak{M}_S \\
    \mathfrak{M}_S \times \mathfrak{M}_S^0. &
    }
\end{align}
Here 
$\mathfrak{M}_S^0$ is the derived moduli stack of zero-dimensional 
sheaves on $S$. 
For $i=1, 2$, the morphism $\ev_i$ sends (\ref{ses:E}) to $E_i$ and $\ev_3$ sends it to $Q$. By exchanging $\ev_1$ and $\ev_2$, we also obtain the diagram 
\begin{align}\label{dia:right}
\xymatrix{
    \mathrm{Hecke}_S \ar[r]^-{\ev_1} \ar[d]_-{(\ev_3, \ev_2)} & \mathfrak{M}_S \\
  \mathfrak{M}_S^0\times   \mathfrak{M}_S. &
    }
\end{align}

\begin{lemma}\label{lem:qsmooth}
(1) In the diagram (\ref{dia:left}), 
the morphism $(\ev_1, \ev_3)$ is quasi-smooth, and $\ev_2$ is quasi-smooth 
and projective on each connected component. 

(2) In the diagram (\ref{dia:right}), 
the morphism $(\ev_3, \ev_2)$ is quasi-smooth, and $\ev_1$ is quasi-smooth 
and projective on each connected component. 
\end{lemma}
\begin{proof}
The morphism $(\ev_1, \ev_3)$ is quasi-smooth 
by~\cite[Proposition~3.10]{PoSa}.

For a pure one-dimensional sheaf $E$ on $S$ 
and a zero-dimensional sheaf $Q$ on $S$, 
we set
\begin{align*}
    E^{\vee}:=\mathbb{D}_S(E)[1] \in \Coh_{\leq 1}^{\heartsuit}(S), \ 
    Q^{\vee}:=\mathbb{D}_S(Q)[2]\in \Coh^{\heartsuit}_{0}(S)
\end{align*}
Here $\mathbb{D}_S(-)=\cHom(-, \mathcal{O}_S)$ is the derived dual. 
There is an equivalence of derived stacks 
\begin{align*}
    \mathbb{D} \colon \mathrm{Hecke}_S \stackrel{\sim}{\to}
    \mathrm{Hecke}_S
\end{align*}
given by 
\begin{align*}
    \mathbb{D}(0\to E_1 \to E_2 \to Q \to 0)=
    (0\to E_2^{\vee} \to E_1^{\vee} \to Q^{\vee} \to 0).
\end{align*}
It fits into the following commutative diagram: 
\begin{align}\label{dia:heckedual}
    \xymatrix{
\mathrm{Hecke}_S \ar[r]^-{\sim}_-{\mathbb{D}} \ar[d]_-{(\ev_1, \ev_3)} & \mathrm{Hecke}_S \ar[d]^-{(\ev_3, \ev_2)} \\
\mathfrak{M}_S \times \mathfrak{M}_S^0 \ar[r]^-{\sim}_{\mathbb{D}} & 
\mathfrak{M}_S^0 \times \mathfrak{M}_S.
    }
\end{align}
Here the bottom equivalence is given by $(E, Q) \mapsto (Q^{\vee}, E^{\vee})$.
The quasi-smoothness of $(\ev_3, \ev_2)$ follows from the quasi-smoothness of 
$(\ev_1, \ev_3)$ together with the above diagram. 

Since $(\ev_1, \ev_3)$ and $(\ev_3, \ev_2)$ are quasi-smooth, 
and $\mathfrak{M}_S^0$ is a quasi-smooth stack, 
the morphisms $\ev_1, \ev_2$ are also quasi-smooth. 
The morphism $\ev_2$ is also projective on each connected component since it is a relative derived Quot scheme; the fiber of $\ev_2$ at $E_2$ classifies quotients $E_2\twoheadrightarrow Q$. By the diagram (\ref{dia:heckedual}) and the fact that $\ev_2$ is projective on each component, the morphism $\ev_1$ is also projective on each component. 
\end{proof}

\begin{lemma}\label{lem:c1L}
For any map $\nu \colon \bgm \to \mathrm{Hecke}_S$ and $i=1, 2$, we have 
\begin{align*}
    c_1(\nu^* \mathbb{L}_{\ev_i}^{>0})= c_1(\nu^* \mathbb{L}_{\ev_i}^{<0})=0
\end{align*}
and 
  \begin{align}\label{vanish:ev3}
    c_1(\nu^* \mathbb{L}_{(\ev_i, \ev_3)}^{>0})= c_1(\nu^* \mathbb{L}_{(\ev_i, \ev_3)}^{<0})=0.
\end{align}  
\end{lemma}
\begin{proof}
Similarly to (\ref{dsum2}), a map $\nu \colon \bgm \to \mathrm{Hecke}_S$ corresponds to the graded object 
  \begin{align}\label{hecke:graded}
      \bigoplus_{m\in \mathbb{Z}}(0\to E_{1, m} \to E_{2, m} \to Q_m \to 0).
  \end{align}
  Here 
  the sequence for each $m$ is an exact sequence in $\Coh_{\leq 1}^{\heartsuit}(S)$,
  each non-zero $E_{i, m}$ is pure one-dimensional and $Q_m$ is zero-dimensional. 
Since we have 
\begin{align*}
    \nu^* \mathbb{L}_{\ev_2}=\Hom^*\left(\bigoplus_{m\in \mathbb{Z}}E_{1, m}, \bigoplus_{m\in \mathbb{Z}}Q_m  \right)^{\vee}
\end{align*}
we have 
\begin{align*}
    c_1(\nu^* \mathbb{L}_{\ev_2}^{<0})=
    \sum_{m<n}\chi(E_{1, m}, Q_n)\cdot (m-n).
\end{align*}
Since $\chi(E_{1, m}, Q_n)=0$ as $Q_n$ is a zero-dimensional sheaf, 
we obtain $c_1(\nu^* \mathbb{L}_{\ev_2}^{<0})=0$. 
The vanishing $c_1(\nu^* \mathbb{L}_{\ev_2}^{>0})=0$
follows similarly. 
The vanishings $c_1(\nu^* \mathbb{L}_{\ev_1}^{>0})=0$ and 
$c_1(\nu^* \mathbb{L}_{\ev_1}^{<0})=0$ hold by the commutative diagram (\ref{dia:heckedual}). 

The vanishings (\ref{vanish:ev3}) follow similarly from 
\begin{align*}
    \nu^* \mathbb{L}_{(\ev_1, \ev_3)}=\Hom^*\left(\bigoplus_{m\in \mathbb{Z}}Q_m, \bigoplus_{m\in \mathbb{Z}}E_{1, m}  \right)^{\vee}[-1]
\end{align*}
and $\chi(Q_m, E_{1,n})=0$. 
\end{proof}

\subsection{Categorical Hall product}
We have the decomposition into open and closed substacks 
\begin{align*}\mathfrak{M}_S^0=\bigsqcup_{d\geq 0}\mathfrak{M}_S^0(d), \ \mathfrak{M}_S^0(0)=\mathrm{pt},
\end{align*}
where $\mathfrak{M}_S^0(d)$ 
classifies zero-dimensional sheaves of length $d$.  
Let
\begin{align*}
    \mathcal{Q} \in \Coh(S \times \mathfrak{M}_S^0)
\end{align*}
be the universal zero-dimensional sheaf. We consider the following map 
\begin{align*}
    \delta_N \colon \mathbb{Q} \to \mathrm{Pic}(\mathfrak{M}_S^0)_{\mathbb{Q}}^{\mathrm{CL}}, \ 
    r \mapsto (\det p_{\mathfrak{M}*}\mathcal{Q})^r.
\end{align*}
Here $p_{\mathfrak{M}} \colon S \times \mathfrak{M}_S^0 \to \mathfrak{M}_S^0$ is the projection. 
For $d\in \mathbb{Z}_{>0}$, we define 
\begin{align*}
    \mathbb{T}_S(d)_0 :=\LL(\mathfrak{M}_S^0(d))_{\delta=0}.
\end{align*}
This is the quasi-BPS category for zero-dimensional sheaves (with weight zero) on the local surface $X=\mathrm{Tot}_S(\omega_S)$, considered in~\cite{PT0, PT2}. 

Let $\mathfrak{M}_S^0(d_1, d_2)$ be the derived moduli stack which classifies 
exact sequences 
\begin{align*}
    0\to Q_1 \to Q \to Q_2 \to 0
\end{align*}
where $Q_i$ are zero-dimensional sheaves of length $d_i$. 
We have the evaluation morphisms 
\begin{align*}
    \mathfrak{M}_S^0(d_1) \times \mathfrak{M}_S^0(d_2) \stackrel{q}{\leftarrow}
    \mathfrak{M}_S^0(d_1, d_2) \stackrel{p}{\to} \mathfrak{M}_S^0(d_1+d_2).
\end{align*}
Since $p$ is proper and $q$ is quasi-smooth, we have the categorical Hall product 
\begin{align}\label{catHall}
    p_{*}q^* \colon 
    \Coh(\mathfrak{M}_S^0(d_1)) \otimes \Coh(\mathfrak{M}_S^0(d_2)) \to 
    \Coh(\mathfrak{M}_S^0(d_1+d_2)).
\end{align}

On the other hand, we define 
\begin{align*}
    \mathbb{H} :=\bigoplus_{d\geq 0}\mathbb{T}_S(d)_{0}, \ \mathbb{T}_S(0)_0=\mathrm{Vect}.
\end{align*}
\begin{lemma}\label{lem:monoidal}
    The categorical Hall product (\ref{catHall}) restricts to a monoidal structure on $\mathbb{H}$. 
\end{lemma}
\begin{proof}
    The same computation as in Lemma~\ref{lem:c1L} shows that for any 
    $\nu \colon \bgm \to \mathfrak{M}_S^0(d_1, d_2)$ we have 
    \begin{align*}
        c_1(\nu^* \mathbb{L}_{p}^{>0})= c_1(\nu^* \mathbb{L}_{p}^{<0})
        = c_1(\nu^* \mathbb{L}_{q}^{>0})= c_1(\nu^* \mathbb{L}_{q}^{<0})=0.
    \end{align*}
    Therefore the lemma holds by Proposition~\ref{prop:limit-correspondence}. 
\end{proof}

\subsection{Harder--Narasimhan strata under the Hecke correspondence}

The diagram (\ref{dia:left}) restricts to the diagram of open and closed substacks 
\begin{align}\label{dia:left:con}
\xymatrix{
    \mathrm{Hecke}_S(\beta, \chi, \chi+d) \ar[r]^-{\ev_2} \ar[d]_-{(\ev_1, \ev_3)} & \mathfrak{M}_S(\beta, \chi+d) \\
    \mathfrak{M}_S(\beta, \chi) \times \mathfrak{M}_S^0(d). &
    }
\end{align}
For each $\mu \in I_{>0}$, let 
\begin{align*}
    \mathfrak{S}_{> \mu}(\beta, \chi) \subset \mathfrak{M}_S(\beta, \chi)
\end{align*}
be the union of Harder--Narasimhan strata with respect to $\sigma$-stability, with instability value (\ref{id:mu}) greater than $\mu$. We have the following lemma: 
\begin{lemma}\label{lem:HN}
For each $\mu'>0$, there is $\mu \gg 0$ such that 
    \begin{align*}
        (\ev_1, \ev_3)(\ev_2)^{-1}(\mathfrak{S}_{> \mu}(\beta, \chi+d))
        \subset \mathfrak{S}_{> \mu'}(\beta, \chi)\times 
        \mathfrak{M}_S^0(d).
    \end{align*}
\end{lemma}
\begin{proof}
Take a pure one-dimensional sheaf 
$E\in \Coh^{\heartsuit}_{\leq 1}(S)$ such that $\cl(E)=(\beta, \chi+d)$, and consider an exact sequence 
\begin{align*}
    0\to E' \to E \to Q \to 0
\end{align*}
where $Q$ is zero-dimensional of length $d$.
We may assume that $E$ is $\sigma$-unstable, since only such sheaves
occur in $\mathfrak S_{>\mu}(\beta,\chi+d)$ for $\mu>0$.
Let 
\begin{align*}
    0\to E_1 \to E\to E_2 \to 0
\end{align*}
be an exact sequence of pure one-dimensional sheaves. 
Let $Q_1$ be defined by 
\begin{align*}
    Q_1=\operatorname{Im}(E_1 \hookrightarrow E \twoheadrightarrow Q)
\end{align*}
and $Q_2=Q/Q_1$. We have the commutative diagram of exact sequences 
\begin{align*}
    \xymatrix{
0 \ar[r] & E_1 \ar[r] \ar[d] & E \ar[r] \ar[d] & E_2 \ar[r] \ar[d] & 0 \\
0 \ar[r] & Q_1 \ar[r] & Q \ar[r] & Q_2 \ar[r] & 0.
    }
\end{align*}
Here the vertical arrows are surjections. 
By taking the kernels, we obtain the exact sequence 
\begin{align*}
    0\to E_1' \to E' \to E_2' \to 0
\end{align*}
such that $\chi(E_i')=\chi(E_i)-d_i$, where $d_i$ is the length of $Q_i$. 

We apply the above construction to the Harder--Narasimhan filtration of $E$. Let 
\begin{align}\label{HN:E}
    0=E_0\subset E_1 \subset \cdots \subset E_k=E
\end{align}
be the Harder--Narasimhan filtration of $E$. The subquotient 
$F_i=E_i/E_{i-1}$ is $\sigma$-semistable, and we set $\cl(F_i)=(\beta_i, \chi_i)$.
The above construction gives a filtration 
\begin{align}\label{HN:E'}
    0=E_0'\subset E_1' \subset \cdots \subset E_k' =E'
\end{align}
such that $F_i'=E_i'/E_{i-1}'$ satisfies 
\begin{align*}\cl(F_i')=(\beta_i, \chi_i-d_i), \ 
d=d_1+\cdots +d_k, \ d_i \geq 0.
\end{align*}

By \eqref{id:mu}, the instability value of the Harder--Narasimhan
filtration \eqref{HN:E} is
\begin{equation}
\begin{split}
\mu_{\mathrm{HN}}(E)
&=
(H\cdot \beta)
\left(
 \sum_{i=1}^{k}
 (H\cdot\beta_i)
 \mu^{\sigma}(\beta_i,\chi_i)^2
 -
 (H\cdot\beta)
 \mu^{\sigma}(\beta,\chi+d)^2
\right)^{1/2}.
\end{split}
\label{id:mu:hn}
\end{equation}
For simplicity, we set 
\begin{align*}
h_i:=H\cdot\beta_i, \ 
h:=H\cdot\beta, \ 
s_i:=\mu^{\sigma}(\beta_i,\chi_i), \ 
\overline{s}:=\mu^{\sigma}(\beta,\chi+d).
\end{align*}
Then we have
\[
h=\sum_{i=1}^{k}h_i,
\quad
h\overline{s}=\sum_{i=1}^{k}h_i s_i.
\]
We also put $\lambda_i:=s_i-\overline{s}$
and set 
$
V:=\sum_{i=1}^{k}h_i\lambda_i^2>0$. 
Then
\[
\sum_{i=1}^{k}h_i\lambda_i=0,
\quad
\mu_{\mathrm{HN}}(E)=h\sqrt{V}.
\]
By (\ref{weight:nu}), up to multiplication by a common positive number, the tuple
$
(\lambda_1,\ldots,\lambda_k)$
is the one-parameter subgroup associated with the
Harder--Narasimhan filtration of \(E\).

We equip the filtration \eqref{HN:E'} with the same weights
\((\lambda_i)_{i=1}^{k}\). After multiplication by
a common positive integer, these weights define a \(\Theta\)-filtration
\(f'\) of \(E'\). Note that the subquotients \(F_i'\) need not be
semistable. Set
\[
s_i'
:=
\mu^{\sigma}(\beta_i,\chi_i-d_i)
=
s_i-\frac{d_i}{h_i},
\quad
\overline{s}'
:=
\mu^{\sigma}(\beta,\chi)
=
\overline{s}-\frac{d}{h}.
\]
The same computation as the one giving \eqref{id:mu}, applied to
the weighted filtration \(f'\), gives
\begin{align*}
    \mu(f')=h \frac{ \sum_{i=1}^{k}
 h_i(s_i'-\overline{s}')\lambda_i}{\sqrt{\sum_{i=1}^{k}h_i\lambda_i^2}}
 =h\frac{V-\sum_{i=1}^k d_i \lambda_i}{\sqrt{V}}.
\end{align*}
Indeed,
\[
s_i'-\overline{s}'
=
\lambda_i-\frac{d_i}{h_i}+\frac{d}{h},
\]
and the last term does not contribute since
\(\sum_i h_i\lambda_i=0\).

Choose a positive integer \(N\) such that \(NH\) is an ample
Cartier divisor. Since each \(\beta_i\) is a nonzero effective
one-cycle, we have
\[
Nh_i=(NH)\cdot\beta_i\in\mathbb{Z}_{>0},
\]
and hence \(h_i\geq 1/N\). By the Cauchy--Schwarz inequality,
\[
\begin{split}
\left|
 \sum_{i=1}^{k}d_i\lambda_i
\right|
\leq
\left(
 \sum_{i=1}^{k}\frac{d_i^2}{h_i}
\right)^{1/2}
\left(
 \sum_{i=1}^{k}h_i\lambda_i^2
\right)^{1/2}
\leq
d\sqrt{N}\,\sqrt{V},
\end{split}
\]
where we have used
\begin{align*}
d_i\geq 0, \quad 
\sum_{i=1}^{k}d_i=d,
\quad
\sum_{i=1}^{k}d_i^2\leq d^2.
\end{align*}
It follows that
\begin{align*}
\mu(f')
\geq
h\sqrt{V}-hd\sqrt{N}
=
\mu_{\mathrm{HN}}(E)-hd\sqrt{N}.
\end{align*}
Since the Harder--Narasimhan filtration maximizes the numerical
invariant, we obtain
\[
\mu_{\mathrm{HN}}(E')
\geq
\mu_{\mathrm{HN}}(E)-hd\sqrt{N}.
\]

We now take
\[
\mu\geq \mu'+hd\sqrt{N}.
\]
If \(E\in \mathfrak{S}_{>\mu}(\beta,\chi+d)\), then
\[
\mu_{\mathrm{HN}}(E')
>
\mu-hd\sqrt{N}
\geq
\mu'.
\]
Therefore
$
E'\in \mathfrak{S}_{>\mu'}(\beta,\chi)$,
which proves the desired inclusion.

\end{proof}

\subsection{Construction of the Hecke operators}
We construct the Hecke operators using the following lemma: 
\begin{lemma}\label{lem:delpull}
For any $D\in N^1(S)_{\mathbb{Q}}$, we have 
\begin{align*}
    \ev_1^*\delta_N(D)=\ev_2^* \delta_N(D) \in \mathrm{Pic}^{\mathrm{CL}}(\mathrm{Hecke}_S)_{\mathbb{Q}}.
\end{align*}
\end{lemma}
\begin{proof}
    For any map $\nu \colon \bgm \to \mathrm{Hecke}_S$ corresponding
    to the graded object (\ref{hecke:graded}), we have 
    \begin{align*}
        \ev_2^* \delta_N(D)(\nu) &=
        \sum_{m\in \mathbb{Z}} m \cdot (D\cdot [E_{2, m}]) \\
        &
        = \sum_{m\in \mathbb{Z}} m \cdot (D\cdot [E_{1, m}])
        =\ev_1^* \delta_N(D)(\nu).
    \end{align*}
\end{proof}

In the diagram (\ref{dia:left:con}), 
the map $\ev_2$ is projective and $(\ev_1, \ev_3)$ is quasi-smooth. Hence we 
have the functor 
\begin{align}\label{funct:ev23}
    \ev_{2*} (\ev_1, \ev_3)^* \colon \Coh(\mathfrak{M}_S(\beta, \chi)\times
    \mathfrak{M}_S^0(d))\to 
    \Coh(\mathfrak{M}_S(\beta, \chi+d)).
\end{align}
\begin{prop}\label{prop:hecke}
Suppose that $-K_S$ is nef. 
For any $D\in N^1(S)_{\mathbb{Q}}$, 
the functor (\ref{funct:ev23}) 
restricts to the functor 
\begin{align*}
    \ev_{2*}(\ev_1, \ev_3)^* \colon \LL(\mathfrak{M}_S(\beta, \chi))_{\delta_N(D)} \otimes \LL(\mathfrak{M}_S^0(d))_{0}
    \to \LL(\mathfrak{M}_S(\beta, \chi+d))_{\delta_N(D)}.
\end{align*}  
\end{prop}
\begin{proof}
By Proposition~\ref{prop:limit-correspondence} and Lemmas~\ref{lem:c1L} and~\ref{lem:delpull}, 
    the functor (\ref{funct:ev23})
    restricts to the functor 
    \begin{align*}
\widetilde{\LL}(\mathfrak{M}_S(\beta, \chi) \times \mathfrak{M}_S^0(d))_{\delta_N(D)\boxtimes 0}
    \to \widetilde{\LL}(\mathfrak{M}_S(\beta, \chi+d))_{\delta_N(D)}.
    \end{align*}
    By~\cite[Theorem~4.27]{TodaLimits}, we have the equivalence 
    \begin{align*}
        \LL(\mathfrak{M}_S(\beta, \chi))_{\delta_N(D)} \otimes 
\LL(\mathfrak{M}_S^0(d))_{0} \stackrel{\sim}{\to}\LL(\mathfrak{M}_S(\beta, \chi) \times \mathfrak{M}_S^0(d))_{\delta_N(D)\boxtimes 0}.
    \end{align*}
    Take an object 
    \begin{align*}
        E\boxtimes T \in 
        \LL(\mathfrak{M}_S(\beta, \chi))_{\delta_N(D)} \otimes 
\LL(\mathfrak{M}_S^0(d))_{0}.
    \end{align*}
By Corollary~\ref{cor:limit-compactness}, it is enough to show that, 
for $\mu \gg 0$, we have 
\begin{align}\label{vanish:Ltilde}
\Hom(\ev_{2*}(\ev_1, \ev_3)^* (E\boxtimes T), A)=0
\end{align}
where $A\in \widetilde{\LL}(\mathfrak{M}_S(\beta, \chi+d))_{\delta_N(D)}$
is supported on $\mathfrak{S}_{> \mu}(\beta, \chi+d)$.

By Grothendieck duality, we have
\begin{align*}
\cHom\left(\ev_3^*T,\ev_2^!A\right)
\simeq
\ev_2^!A\otimes
\ev_3^*\mathbb{D}_{\mathfrak{M}_S^0(d)}(T).
\end{align*}
Here the left-hand side is the internal Hom. 
Therefore, by adjunction for $\ev_1$ and $\ev_2$,
the vanishing (\ref{vanish:Ltilde}) is equivalent to
\begin{align}\label{vanish:Ltilde2}
\Hom\left(
E,
\ev_{1*}\left(
\ev_2^!A\otimes
\ev_3^*\mathbb{D}_{\mathfrak{M}_S^0(d)}(T)
\right)
\right)=0.
\end{align}

   Since $E \in \LL(\mathfrak{M}_S(\beta, \chi))_{\delta_N(D)}$, 
   by Corollary~\ref{cor:limit-compactness} there is $\mu'>0$ such that 
  \begin{align*}\Hom(E, A')=0, \ \mbox{ for every } A' \in \widetilde{\LL}(\mathfrak{M}_S(\beta, \chi))_{\delta_N(D)}, \ 
  \mathrm{Supp}(A') \subset \mathfrak{S}_{> \mu'}(\beta, \chi).
  \end{align*}
  Take $\mu \gg 0$ as in Lemma~\ref{lem:HN}.
  By~(\ref{adjoint:limit-correspondence}), we have
  \begin{align*}
      \ev_{1*}(\ev_2^! A \otimes \ev_3^* \mathbb{D}_{\mathfrak{M}_S^0(d)}(T)) \in \widetilde{\LL}(\mathfrak{M}_S(\beta, \chi))_{\delta_N(D)}
  \end{align*}
  and it is supported on $\mathfrak{S}_{> \mu'}(\beta, \chi)$ by Lemma~\ref{lem:HN}.
  Therefore the vanishing (\ref{vanish:Ltilde2}) holds. 
\end{proof}

Similarly, the diagram (\ref{dia:right}) induces the functor 
\begin{align}\label{funct:ev231}
    \ev_{1*} (\ev_3, \ev_2)^* \colon 
    \Coh(\mathfrak{M}_S^0(d)\times\mathfrak{M}_S(\beta, \chi))\to 
    \Coh(\mathfrak{M}_S(\beta, \chi-d)).
\end{align}
We also have the following proposition: 
\begin{prop}\label{prop:hecke2}
Suppose that $-K_S$ is nef.
For any $D\in N^1(S)_{\mathbb{Q}}$, 
the functor (\ref{funct:ev231}) 
restricts to the functor 
\begin{align*}
    \ev_{1*}(\ev_3, \ev_2)^* \colon
 \LL(\mathfrak{M}_S^0(d))_{0}  \otimes
    \LL(\mathfrak{M}_S(\beta, \chi))_{\delta_N(D)}
    \to \LL(\mathfrak{M}_S(\beta, \chi-d))_{\delta_N(D)}.
\end{align*}  
\end{prop}
\begin{proof}
    The proposition follows from Proposition~\ref{prop:hecke} together with the diagram (\ref{dia:heckedual}). 
\end{proof}

Together, these results give the following result: 

\begin{thm}\label{thm:hecke1}
Suppose that $-K_S$ is nef. 
For a fixed $\beta$ and $D\in N^1(S)_{\mathbb{Q}}$, the direct sum
\begin{align*}
    \bigoplus_{\chi \in \mathbb{Z}}\LL(\mathfrak{M}_S(\beta, \chi))_{\delta_N(D)}
\end{align*}
admits right and left actions of $\mathbb{H}$.    
\end{thm}
\begin{proof}
    This is a consequence of Proposition~\ref{prop:hecke} and Proposition~\ref{prop:hecke2}. 
    Note that associativity of the module actions follows from the usual associativity of the categorical Hall product by considering stacks of two-step flags, as in~\cite{PoSa}. 
\end{proof}

As a corollary, we obtain Hecke actions on 
quasi-BPS categories: 
\begin{cor}\label{cor:hecke2}
Suppose that $-K_S$ is nef. 
Let $B+iH\in N^1(S)_{\mathbb{C}}$ be rational 
such that $H$ is ample, and set $\sigma=\sigma_{B, H}$. 
Suppose that $(\beta, w) \in \Lambda$ is $\sigma$-primitive. 
Then the direct sum 
\begin{align*}
\bigoplus_{\chi\in \mathbb{Z}}\mathbb{T}^{+}_S(\beta, \chi)_{w_{\beta}}^{\sigma}
\end{align*}
admits right and left actions of $\mathbb{H}$. 
\end{cor}
\begin{proof}
By the semiorthogonal decomposition (\ref{sod:2}), if $(\beta, w)\in \Lambda$ is $\sigma$-primitive, 
we have the equivalence 
\begin{align*}
   j_{!} \colon  \mathbb{T}^{+}_S(\beta, \chi)_{w_{\beta}}^{\sigma}
    \stackrel{\sim}{\to} \LL(\mathfrak{M}_S(\beta, \chi))_{w_{\beta}}^{\sigma}.
\end{align*}
    Then the corollary follows from Theorem~\ref{thm:hecke1}. 
\end{proof}

\begin{remark}\label{rmk:heckelift}
Under the $\sigma$-primitivity assumption of Corollary~\ref{cor:hecke2}, if $S$ is symplectic, then for any object 
\begin{align*}
A\in \LL(\mathfrak{M}_S^0(1))_0 \simeq \Coh(S\times k[-1])
\end{align*}
there is a Hecke action 
\begin{align}\label{funct:A}
    *A \colon \mathbb{T}^{+}_S(\beta, \chi)_{w_{\beta}}^{\sigma} \to \mathbb{T}^{+}_S(\beta, \chi+1)_{w_{\beta}}^{\sigma}.
\end{align}
When $S$ is a symplectic surface, the 
above action should be regarded as a categorical lift of the 
Hecke operator on BPS cohomology, recently constructed in~\cite{DHKSV}. 
The Hecke action on BPS cohomology is used in~\cite{DHKSV} to give a proof of the $\chi$-independence conjecture for BPS cohomology of 
symplectic surfaces. 

On the other hand, the functor (\ref{funct:A}) may not be an equivalence for any $A$. 
Indeed, let $S$ be a K3 surface with $\Pic(S)=\mathbb{Z}[H]$ with 
$H^2=2d\geq 2$. Then the functor (\ref{funct:A}) is not fully faithful 
for any choice of $A$; this follows essentially from the fact that there is no one-cycle 
$Z \subset S$ with $Z \cdot \beta=1$, so that the functor (\ref{funct:A}) never induces 
an equivalence of dg-categories of coherent sheaves on the generic fibers of 
$\mathfrak{M}_S \to \rB$. 
\end{remark}

\subsection{Wilson/Hecke compatibility conjecture}
Assume throughout this subsection that $-K_S$ is nef.
The moduli stack $\mathfrak{M}_S^0(1)$ is a $\mathbb{G}_m$-gerbe over 
\begin{align*}
    S\times_{0, \omega_S, 0}S=S[\omega_S^{-1}[-1]].
\end{align*}
Let $i \colon S \hookrightarrow S[\omega_S^{-1}[-1]]$
be the natural closed immersion. 
For $A \in \Coh(S)$, 
consider the object 
\begin{align*}
    i_{*}A \in \LL(\mathfrak{M}_S^0(1))_0=\Coh(S[\omega_S^{-1}[-1]]).
\end{align*}
We have the left Hecke operator associated with $i_{*}A$
\begin{align*}
    \rH_{A}(-):={i_{*}A}*(-) \colon \LL(\mathfrak{M}_S(\beta, \chi))_w^{\sigma} \to \LL(\mathfrak{M}_S(\beta, \chi-1))_w^{\sigma}.
\end{align*}

On the other hand, using the universal sheaf $\cE$ in (\ref{univ:E}), 
the Wilson operator is the functor 
\begin{align*}
    \rW_A \colon 
    \Coh(\mathfrak{M}_S(\beta, w)^{\sigma\text{-ss}})_{-\chi}
    \to  \Coh(\mathfrak{M}_S(\beta, w)^{\sigma\text{-ss}})_{-\chi+1}
\end{align*}
given by 
\begin{align}\label{formula:W}
   \rW_A(-) :=p_{\mathfrak{M}*}((-)\otimes \cE\otimes p_{S}^* A).
\end{align}
We expect that these functors are compatible under 
the conjectural equivalence. Such a compatibility is expected to play a key role
in proving that a given functor is an equivalence; see~\cite{TodaGL2, TodaA, PTtwist}:
\begin{conj}\label{conj:WH}
    Under an equivalence (\ref{conj:equiv}),
    the functors $\rH_A$ and $\rW_A$ commute for any $A\in \Coh(S)$, i.e. we have the commutative diagram 
    \begin{align*}
        \xymatrix{
          \Coh(\mathfrak{M}_S(\beta, w)^{\sigma\text{-ss}})_{-\chi_{\beta}}
          \ar[r]^-{\Phi} \ar[d]_-{\rW_A} & \LL(\mathfrak{M}_S(\beta, \chi))_{w_{\beta}}^{\sigma} \ar[d]^-{\rH_A} \\
      \Coh(\mathfrak{M}_S(\beta, w)^{\sigma\text{-ss}})_{-\chi_{\beta}+1}
          \ar[r]^-{\Phi} &  \LL(\mathfrak{M}_S(\beta, \chi-1))_{w_{\beta}}^{\sigma}.
        }
    \end{align*}
\end{conj}

\section{Dolbeault geometric Langlands-type conjecture for local surfaces}\label{sec:DLCY3}
In this section, we discuss the Dolbeault geometric Langlands-type conjecture for local surfaces, which are non-compact Calabi--Yau 3-folds. Most of this section is conjectural and speculative. 

\subsection{Categorical DT theory for local surfaces}
So far, our main conjecture has assumed that $-K_S$ is nef. This assumption 
is related to the geometry of moduli of sheaves on the 
Calabi--Yau 3-fold 
\begin{align*}
    X=\mathrm{Tot}_S(\omega_S) \stackrel{\pi_X}{\to} S.
\end{align*}
Let $\mathfrak{M}_X$ be the derived moduli stack of 
compactly supported pure one-dimensional sheaves $E$ on $X$. 
Then we have 
\begin{align*}
    \mathfrak{M}_X=\Omega_{\mathfrak{M}_S}[-1]
\end{align*}
where the right-hand side is the $(-1)$-shifted cotangent, see~\cite{T}. 
We also have the decomposition into open and closed substacks 
\begin{align*}
    \mathfrak{M}_X=\bigsqcup_{(\beta, \chi)}\mathfrak{M}_X(\beta, \chi)
\end{align*}
where each component lies over $\mathfrak{M}_S(\beta, \chi)$.

On the other hand, the stability condition 
$\sigma=(Z, \Coh_{\leq 1}^{\heartsuit}(S))$
on $\Coh_{\leq 1}(S)$
induces the stability condition
\begin{align*}
    \sigma_X=(Z_X, \Coh_{\leq 1}^{\heartsuit}(X))
\end{align*}
with central charge
\begin{align*}
    Z_X \colon K(\Coh_{\leq 1}(X)) \stackrel{\pi_{X*}}{\to}
    K(\Coh_{\leq 1}(S)) \stackrel{Z}{\to} \mathbb{C}.
\end{align*}
\begin{lemma}\label{lem:X}
Suppose that $-K_S$ is nef. Then we have 
\begin{align*}
    \mathfrak{M}_X^{\sigma_X\text{-ss}}=\Omega_{\mathfrak{M}_S^{\sigma\text{-ss}}}[-1].
\end{align*}
\end{lemma}
\begin{proof}
    It is enough to show that, for a $\sigma_X$-semistable sheaf 
    $E\in \Coh_{\leq 1}^{\heartsuit}(X)$, the object
    $F=\pi_{X*}E \in \Coh_{\leq 1}^{\heartsuit}(S)$ is $\sigma$-semistable. Note that $\Coh_{\leq 1}^{\heartsuit}(X)$
    is equivalent to the category of pairs 
    \begin{align*}
    (F, \theta), \ 
    F\in \Coh_{\leq 1}^{\heartsuit}(S), \ \theta \colon F \to F\otimes \omega_S.
    \end{align*}
    Suppose that $F$ is not semistable. Then we have the non-trivial 
    Harder--Narasimhan filtration 
    \begin{align*}
        0=F_0 \subset F_1 \subset \cdots \subset F_k=F
    \end{align*}
    of $F$ with respect to $\sigma$-stability. Then Lemma~\ref{lem:S} implies that 
    each $F_i$ is preserved by $\theta$.
    Thus the filtration lifts to a destabilizing filtration of $E$,
    contradicting the $\sigma_X$-semistability of $E$.
    Hence $F$ is $\sigma$-semistable. 
\end{proof}

Let $\mathcal{Z}_X^{\text{us}}$ be the unstable locus 
\begin{align*}
    \mathcal{Z}_X^{\text{us}} :=\mathfrak{M}_X \setminus \mathfrak{M}_X^{\sigma_X\text{-ss}}.
\end{align*}
We denote by 
\begin{align*}
    \IndCoh_{\mathcal{Z}_X^{\text{us}}}(\mathfrak{M}_S)
    \subset \IndCoh(\mathfrak{M}_S)
\end{align*}
the subcategory of objects with singular supports~\cite{AG} 
contained in $\mathcal{Z}_X^{\text{us}}$.
Then the $\mathbb{G}_m$-equivariant (ind-)DT category for 
$X$ is given by the quotient~\cite{T} 
\begin{align}\label{indDT}
    \mathrm{IndDT}^{\mathbb{G}_m}(\mathfrak{M}_X(\beta, \chi)^{\sigma_X\text{-ss}}):=
    \IndCoh(\mathfrak{M}_S(\beta, \chi))/\IndCoh_{\mathcal{Z}_X^{\text{us}}}(\mathfrak{M}_S(\beta, \chi)).
\end{align}
Here $\mathbb{G}_m$ acts on the fibers of $X \to S$.
We propose the following conjecture:
\begin{conj}\label{CYconj1}
The dg-category (\ref{indDT}) is compactly generated, 
with compact objects
\begin{align*}
     \mathrm{DT}^{\mathbb{G}_m}(\mathfrak{M}_X(\beta, \chi)^{\sigma_X\text{-ss}})=
    \left(\Coh(\mathfrak{M}_S(\beta, \chi))/\Coh_{\mathcal{Z}_X^{\text{us}}}(\mathfrak{M}_S(\beta, \chi))\right)^{\mathrm{idem}}.
\end{align*}
Here $(-)^{\mathrm{idem}}$ stands for the idempotent completion. 

Furthermore, we have the central weight decomposition 
\begin{align*}
     \mathrm{DT}^{\mathbb{G}_m}(\mathfrak{M}_X(\beta, \chi)^{\sigma_X\text{-ss}})=\bigoplus_{w\in \mathbb{Z}}
       \mathrm{DT}^{\mathbb{G}_m}(\mathfrak{M}_X(\beta, \chi)^{\sigma_X\text{-ss}})_w.
\end{align*}
and
there are subcategories
\begin{align}\label{sub:T}
    \mathbb{T}^{\mathbb{G}_m,\pm}_X(\beta, \chi)_w^{\sigma_X} \subset 
     \mathrm{DT}^{\mathbb{G}_m}(\mathfrak{M}_X(\beta, \chi)^{\sigma_X\text{-ss}})_{w}
\end{align}
such that we have semiorthogonal decompositions 
\begin{align}\label{DTSOD}
        \mathrm{DT}^{\mathbb{G}_m}(\mathfrak{M}_X(\beta, w)^{\sigma_X\text{-ss}})_{-\chi_{\beta}}
        =\left\langle \bigotimes_{i=1}^k \mathbb{T}_X^{\mathbb{G}_m,-}(\beta_i, w_i)_{-\chi_{i, \beta_i}}^{\sigma_X}
        \right\rangle.
    \end{align}
    Here the right-hand side ranges over the partitions 
    \begin{align}\label{decom:betachi}
    (\beta, \chi, w)=(\beta_1, \chi_1, w_1)+\cdots+(\beta_k, \chi_k, w_k)
    \end{align}
    such that 
    \begin{align}\label{decom:musigma}
        \mu^{\sigma}(\beta, w)=\mu^{\sigma}(\beta_i, w_i), \ 
        \mu^{\sigma}(\beta_1, \chi_1)>\cdots>\mu^{\sigma}(\beta_k, \chi_k).
    \end{align}
\end{conj}
\begin{remark}\label{rmk:CY3ss}
By Lemma~\ref{lem:X}, if $-K_S$ is nef, we have 
\begin{align*}
     \mathrm{DT}^{\mathbb{G}_m}(\mathfrak{M}_X(\beta, \chi)^{\sigma_X\text{-ss}})_w \simeq \Coh(\mathfrak{M}_S(\beta, \chi)^{\sigma\text{-ss}})_w.
\end{align*}
In this case, the decomposition in Conjecture~\ref{CYconj1} follows from
Theorem~\ref{thm:bpt}, and we take 
\begin{align}\label{qbps:nef}
\mathbb{T}_X^{\mathbb{G}_m,\pm}(\beta, \chi)_w^{\sigma_X}=
\mathbb{T}^{\pm}_S(\beta, \chi)_w^{\sigma}.  
\end{align}
\end{remark}
\subsection{Limit categories for local surfaces}
On the automorphic side, we define
\begin{align}\label{limit:X}
    \IndL^{\mathbb{G}_m}(\mathfrak{M}_X(\beta, \chi))_w^{\sigma_X}:=
    \IndL(\mathfrak{M}_S(\beta, \chi))_w^{\sigma}.
\end{align}
We propose the following conjecture, without the assumption 
that $-K_S$ is nef: 
\begin{conj}\label{CYconj2}
The dg-category (\ref{limit:X}) is compactly generated by 
\begin{align*}
     \LL^{\mathbb{G}_m}(\mathfrak{M}_X(\beta, \chi))_w^{\sigma_X}:=
    \LL(\mathfrak{M}_S(\beta, \chi))_w^{\sigma}
\end{align*}
and it admits a semiorthogonal decomposition 
\begin{align}\label{SOD:CY3L}
    \LL^{\mathbb{G}_m}(\mathfrak{M}_X(\beta, \chi))_{w_{\beta}}^{\sigma_X}&=
       \left\langle \bigotimes_{i=1}^k \mathbb{T}_X^{\mathbb{G}_m,+}(\beta_i, \chi_i)_{w_{i, \beta_i}}^{\sigma_X}
        \right\rangle.
        \end{align}
        Here the
        dg-categories in the semiorthogonal summands are (\ref{sub:T}), and the
        right-hand side ranges over the partitions (\ref{decom:betachi})
        satisfying (\ref{decom:musigma}). 
\end{conj}

\begin{remark}\label{rmk:CY3HN}
 The semiorthogonal decomposition (\ref{SOD:CY3L}) should be associated with the Harder--Narasimhan stratification of $\mathfrak{M}_X$, 
    not $\mathfrak{M}_S$. 
    
    In the case that $-K_S$ is nef, the proof of Lemma~\ref{lem:X} also shows that the Harder--Narasimhan stratification of $\mathfrak{M}_X$ is 
    the pull-back of that of $\mathfrak{M}_S$.
    In this case,  Conjecture~\ref{CYconj2} follows from Theorem~\ref{thm:sodLS}, and 
    the semiorthogonal decompositions (\ref{SOD:CY3L}) should be the same as (\ref{sod:2}).  
\end{remark}

\subsection{Formulation of the conjecture for local surfaces}
We formulate the Dolbeault geometric Langlands-type conjecture
for $X=\mathrm{Tot}_S(\omega_S)$ as follows.
\begin{conj}\label{conj:local-surfaces}
For any smooth quasi-projective surface $S$
(with $-K_S$ not necessarily nef),
there is a $\rB$-linear equivalence
\begin{align}\label{conj:3fold}
  \mathrm{IndDT}^{\mathbb{G}_m}(\mathfrak{M}_X(\beta, w)^{\sigma_X\text{-ss}})_{-\chi_{\beta}}
  \simeq 
    \IndL^{\mathbb{G}_m}(\mathfrak{M}_X(\beta, \chi))_{w_{\beta}}^{\sigma_X}.
\end{align}

Furthermore, the semiorthogonal decompositions (\ref{DTSOD}), (\ref{SOD:CY3L}) 
on the subcategories of compact objects are compatible with the 
above equivalence, hence there is an equivalence 
    \begin{align}\label{qbps:CY3}
          \mathbb{T}^{\mathbb{G}_m,-}_X(\beta, w)_{-\chi_{\beta}}^{\sigma_X} 
          \simeq   \mathbb{T}^{\mathbb{G}_m,+}_X(\beta, \chi)_{w_{\beta}}^{\sigma_X}.
    \end{align}
    \end{conj}

    \begin{remark}
When $-K_S$ is nef, an equivalence (\ref{conj:3fold}) is equivalent to Conjecture~\ref{conj:CY3}. 
 In particular
    in this case, an equivalence (\ref{qbps:CY3}) should be the same as the equivalence in Conjecture~\ref{conj:qbps}. 
    \end{remark}
\begin{remark}
     Note that the categories in (\ref{sub:T}) are not yet defined,
    so an equivalence (\ref{qbps:CY3}) remains speculative. If they exist, then under suitable
    conditions they are expected to categorify
    BPS invariants~\cite{DM, BDIKP} on the Calabi--Yau 3-fold $X$.
\end{remark}

   \subsection{Wilson/Hecke compatibility for local surfaces}
   Finally, we formulate a conjecture on the Wilson/Hecke compatibility of a conjectural equivalence (\ref{conj:3fold}). 
   Let $\mathfrak{M}_X^0(d)$ be the derived moduli stack of zero-dimensional sheaves on $X$ of length $d$. 
   We define 
   \begin{align*}
       \LL^{\mathbb{G}_m}(\mathfrak{M}_X^0(d))_0:=\LL(\mathfrak{M}_S^0(d))_0.
   \end{align*}
Then the direct sum
\begin{align*}
    \mathbb{H}_X^{\mathbb{G}_m}:=\bigoplus_{d\geq 0}
    \LL^{\mathbb{G}_m}(\mathfrak{M}_X^0(d))_0
\end{align*}
    admits a monoidal structure by the categorical Hall product, 
    and left and right actions on 
    \begin{align*}
        \bigoplus_{\chi \in \mathbb{Z}}\IndL^{\mathbb{G}_m}(\mathfrak{M}_X(\beta, \chi))_{w}^{\sigma_X}.
    \end{align*}

    On the other hand, we define 
\begin{align}\label{def:WS}
\mathbb{W}_S:=
  \bigoplus_{d\geq 0}\Coh(S^d/\mathfrak{S}_d).
  \end{align}
  It admits a monoidal structure given by the push-forward along 
     \begin{align*}
         S^{d_1}/\mathfrak{S}_{d_1}\times S^{d_2}/\mathfrak{S}_{d_2} 
         =S^{d_1+d_2}/(\mathfrak{S}_{d_1}\times \mathfrak{S}_{d_2}) 
         \to S^{d_1+d_2}/\mathfrak{S}_{d_1+d_2}.
     \end{align*}
The universal sheaf $\cE$ in (\ref{univ:E}) defines the 
object 
\begin{align*}
    \cE^{\boxtimes d} \in \Coh(S^d/\mathfrak{S}_d \times \mathfrak{M}_S)
\end{align*}
which defines the functor 
\begin{align*}
    \Coh(S^d/\mathfrak{S}_d)\otimes \Coh(\mathfrak{M}_S)_{*}
    \to \Coh(\mathfrak{M}_S)_{*+d}
\end{align*}
given by 
\begin{align*}
    (-)\mapsto p_{\mathfrak{M}*}((-)\otimes \cE^{\boxtimes d})
\end{align*}
By replacing $\cE$ with $\cE^{\vee}=\mathbb{D}(\cE)[1]$, 
we also obtain 
\begin{align*}
  \Coh(\mathfrak{M}_S)_{*} \otimes 
    \Coh(S^d/\mathfrak{S}_d)
    \to \Coh(\mathfrak{M}_S)_{*-d}.
\end{align*}
These functors descend to the singular support quotients
and define, respectively, the left action (using $\cE$) and the right action
(using $\cE^{\vee}$) of $\mathbb{W}_S$ on 
\begin{align}\label{indDT:chi}
    \bigoplus_{\chi\in \mathbb{Z}}\mathrm{IndDT}^{\mathbb{G}_m}(\mathfrak{M}_X(\beta, w)^{\sigma_X\text{-ss}})_{-\chi}.
\end{align}

Consider the following 
\begin{align}\label{def:WX}
\mathbb{W}_X^{\mathbb{G}_m}:=
  \bigoplus_{d\geq 0}\mathrm{MF}^{\mathbb{G}_m}(X^d/\mathfrak{S}_d, 0).
  \end{align}
  It admits a monoidal structure given by the push-forward along 
     \begin{align*}
         X^{d_1}/\mathfrak{S}_{d_1}\times X^{d_2}/\mathfrak{S}_{d_2} 
         =X^{d_1+d_2}/(\mathfrak{S}_{d_1}\times \mathfrak{S}_{d_2}) 
         \to X^{d_1+d_2}/\mathfrak{S}_{d_1+d_2}.
     \end{align*}
       The following conjecture generalizes the conjecture in~\cite{PT0}
  for $X=\mathbb{C}^3$:
  \begin{conj}\label{conj:monoidal}
  There is a monoidal equivalence 
  \begin{align}\label{mon:equiv}
 \mathbb{W}_X^{\mathbb{G}_m} \simeq 
  \mathbb{H}_X^{\mathbb{G}_m}.
  \end{align}
  \end{conj}

The projection $\pi_X \colon X \to S$ defines the morphism 
\begin{align*}
    \pi_X^{(d)} \colon X^d/\mathfrak{S}_d \to S^d/\mathfrak{S}_d.
\end{align*}
The pull-back along the above morphism 
determines the monoidal functor 
\begin{align}\label{equiv:WX}
    \pi_X^* \colon \mathbb{W}_S \to \mathbb{W}_X^{\mathbb{G}_m}. 
\end{align}

\begin{conj}\label{conj:final}
The $\mathbb{W}_S$-action on (\ref{indDT:chi}) lifts to a monoidal action of 
$\mathbb{W}_X^{\mathbb{G}_m}$ through the functor (\ref{equiv:WX})
preserving the compact objects of (\ref{indDT:chi}), 
and 
the equivalences (\ref{mon:equiv}) and (\ref{conj:3fold})
are compatible with the above monoidal actions; namely, we have the commutative diagram 
\begin{align*}
\xymatrix@C=3em@R=1em{
\mathbb{W}_X^{\mathbb{G}_m}
  \ar[r]^-{\sim}
  \ar@{}[d]|{\curvearrowright}
&
\mathbb{H}_X^{\mathbb{G}_m}
  \ar@{}[d]|{\curvearrowright}
\\
\displaystyle\bigoplus_{\chi\in \mathbb{Z}}
\mathrm{DT}^{\mathbb{G}_m}
(\mathfrak{M}_X(\beta,w)^{\sigma_X\text{-ss}})_{-\chi_{\beta}}
  \ar[r]^-{\sim}
&
\displaystyle\bigoplus_{\chi\in \mathbb{Z}}
\LL^{\mathbb{G}_m}
(\mathfrak{M}_X(\beta,\chi))_{w_{\beta}}^{\sigma_X}.
}
\end{align*}
    
\end{conj}

\newcommand{\etalchar}[1]{$^{#1}$}
\providecommand{\bysame}{\leavevmode\hbox to3em{\hrulefill}\thinspace}
\providecommand{\MR}{\relax\ifhmode\unskip\space\fi MR }
\providecommand{\MRhref}[2]{%
  \href{http://www.ams.org/mathscinet-getitem?mr=#1}{#2}
}
\providecommand{\href}[2]{#2}

\end{document}